\documentclass[12pt]{article}
\usepackage{cite}
\usepackage{amsmath,amssymb,amsthm,mathrsfs}
\usepackage{titlesec,hyperref}
\usepackage{color}
\usepackage{fancyhdr}
\usepackage{cases}
\usepackage[margin=2.5cm]{geometry}
\usepackage{enumerate}
\usepackage[integrals]{wasysym}
\usepackage{bm}
\usepackage{graphicx}
\usepackage{caption}
\usepackage{float}
\usepackage{subfigure}
\usepackage{indentfirst}
\newcommand{\dd}{\mathop{}\!\mathrm{d}}
\let\del\partial

\newcommand{\supp}{\text{supp}\,}
\newcommand{\MM}{ T}
\newcommand{\MMM}{{\bar T}}
\newcommand{\RR}{\mathring R}
\newcommand{\R}{\mathbb R}

\newcommand{\ZZ}{\mathbb Z}
\newcommand{\NN}{\mathbb N}
\newcommand{\RRR}{\mathring{\bar R}}
\newcommand{\ootimes}{\mathbin{\mathring{\otimes}}}

\let\div\relax
\DeclareMathOperator{\div}{div}

\newcommand{\vin}[0]{v^{\textup{in}}}
\newcommand{\tin}[0]{\theta^{\textup{in}}}

\newcommand{\pex}[0]{\textsf{\textit{p}}}
\newcommand{\vex}[0]{\textsf{\textit{v}}}
\newcommand{\vv}[0]{\bar v}
\newcommand{\tb}[0]{\bar \theta}

\newcommand{\ppp}[0]{\bar{\bar p}}
\newcommand{\pp}[0]{\bar p}

\newcommand{\Rnash}{R_{q+1}^\textup{Nash}}
\newcommand{\Rtransport}{R_{q+1}^\textup{trans}}
\newcommand{\Rosc}{R_{q+1}^\textup{osc}}
\newcommand{\Mtrans}{T_{q+1}^\textup{trans}}
\newcommand{\Mosc}{T_{q+1}^\textup{osc}}
\newcommand{\Mnash}{T_{q+1}^\textup{Nash}}

\newcommand{\TT}[0]{\mathsf{T}}
\newcommand{\TTT}[0]{\mathbb{T}}

\newcommand{\matd}[1]{D_{t,#1}}
\def\dashint{\,\ThisStyle{\ensurestackMath{%
  \stackinset{c}{.2\LMpt}{c}{.5\LMpt}{\SavedStyle-}{\SavedStyle\phantom{\int}}}%
  \setbox0=\hbox{$\SavedStyle\int\,$}\kern-\wd0}\int}
\def\ddashint{\,\ThisStyle{\ensurestackMath{%
  \stackinset{c}{.2\LMpt}{c}{.5\LMpt+.2\LMex}{\SavedStyle-}{%
    \stackinset{c}{.2\LMpt}{c}{.5\LMpt-.2\LMex}{\SavedStyle-}{%
      \SavedStyle\phantom{\int}}}}\setbox0=\hbox{$\SavedStyle\int\,$}\kern-\wd0}\int}
\usepackage{scalerel,stackengine}
\stackMath
\newcommand\widecheck[1]{%
\savestack{\tmpbox}{\stretchto{%
  \scaleto{%
    \scalerel*[\widthof{\ensuremath{#1}}]{\kern-.6pt\bigwedge\kern-.6pt}%
    {\rule[-\textheight/2]{1ex}{\textheight}}
  }{\textheight}%
}{0.5ex}}%
\stackon[1pt]{#1}{\scalebox{-1}{\tmpbox}}%
}

\newcommand{\coloneq}{\mathrel{\mathop:}=}

\newcommand{\ii}{\textup{i}}
\newcommand{\ee}{\textup{e}}

\DeclareMathOperator{\curl}{curl}

\newtheorem{thm}{Theorem}[section]

\newtheorem{lem}[thm]{Lemma}
\newtheorem{prop}[thm]{Proposition}
\newtheorem{con}[thm]{Conjecture}
\theoremstyle{definition}
\newtheorem{defn}[thm]{Definition}

\theoremstyle{remark}
\newtheorem{rem}[thm]{Remark}
\numberwithin{equation}{section}

\allowdisplaybreaks

\begin{document}
\title{On Onsager's type conjecture for the inviscid Boussinesq equations}

\author{
Changxing $\mbox{Miao}$ \footnote{Institute  of Applied Physics and Computational Mathematics, Beijing 100191, P.R. China, email: miao changxing@iapcm.ac.cn},~
Yao $\mbox{Nie}$ \footnote{School of Mathematical Sciences and LPMC, Nankai University, Tianjin, 300071, P.R.China, email address: nieyao@nankai.edu.cn}
,~
Weikui $\mbox{Ye}$ \footnote{School of Mathematical Sciences, South China Normal University, Guangzhou, 510631, P.R. China, email:
904817751@qq.com}
}
\date{}

\maketitle

\begin{abstract}
In this paper, we investigate  the Cauchy problem for the three dimensional inviscid Boussinesq system in the periodic setting.  For $1\le p\le \infty$, we show that the threshold regularity exponent for $L^p$-norm conservation of temperature of this system is  $1/3$, consistent with  Onsager exponent.  More precisely, for $1\le p\le\infty$, every weak solution $(v,\theta)\in C_tC^{\beta}_x$ to the inviscid Boussinesq equations satisfies that $\|\theta(t)\|_{L^p(\mathbb{T}^3)}=\|\theta_0\|_{L^p(\mathbb{T}^3)}$ if $\beta>\frac{1}{3}$, while if $\beta<\frac{1}{3}$, there exist infinitely many weak solutions $(v,\theta)\in C_tC^{\beta}_x$ such that the $L^p$-norm of
temperature is not conserved.  As a byproduct, we are able to construct many weak solutions in $C_tC^{\beta}_x$ for $\beta<\frac{1}{3}$ displaying wild behavior, such as fast kinetic energy dissipation and high oscillation of velocity. Moreover,  we also show that  if a weak solution $(v, \theta)$ of this system  has at least one interval of regularity, then this weak solution $(v,\theta)$ is not unique in  $C_tC^{\beta}_x$ for  $\beta<\frac{1}{3}$.
\end{abstract}

\section{Introduction}
The fact that the dynamics of almost non-stationary astrophysical plasma flows feature a diversity of length scales and speed indicates the widespread existence of turbulence in the universe.  Atmospheric turbulence may result from   small-scale irregular air motions influenced by wind and ocean circulation acts as a turbulent motion on a wide range of scales \cite{YK}.  In certain ranges of scales in the atmosphere and in the ocean,  hydrodynamics is dominated by the interplay  between gravity and the earth's rotation with density variations about a reference state \cite{Majda}.
The Boussinesq system describes the convection phenomena  on these scales in the ocean or atmosphere and has been known as one of the most  geophysical models \cite{A,Majda}. In the periodic setting $\TTT^3=\R^3/\ZZ^3$, the inviscid Boussinesq  system is described by
\begin{equation}
\left\{ \begin{alignedat}{-1}
&\del_t v+(v\cdot\nabla) v  +\nabla p   =  \theta e_3, &{\rm in}\quad \TTT^3\times(0,T),
 \\
&  \del_t \theta+v\cdot\nabla\theta  =0, &{\rm in}\quad \TTT^3\times(0,T),
 \\
&  \nabla \cdot v  = 0, &{\rm in}\quad \TTT^3\times(0,T),
\end{alignedat}\right.  \label{e:B}\tag{B}
\end{equation}
associated with initial data $(v, \theta)\big|_{t=0} = (\vin,  \tin)$, where $v: \TTT^3\times [0,T]\to\R^3$ denotes the velocity of the incompressible fluid, $p:\TTT^3\times [0,T]\to\R$ the pressure field, $\theta:\TTT^3\times [0,T]\to\R$  the temperature, $e_3=(0,0,1)^{T}$ and $\theta e_3$ represents the buoyancy force.

From the viewpoint of mathematics, the system \eqref{e:B} has been attracted considerable attention
in the past years since it is closely related to the incompressible Euler equations. More precisely,
the system \eqref{e:B} boils down to the classical Euler equations when $\theta\equiv 0$ and the 2D inviscid Boussinesq system coincides with  the 3D axisymmetric swirling Euler equations  away from the symmetry axis \cite{MB}. The local well-posedness of smooth solution for the system \eqref{e:B} has been proved (see e.g. \cite{CKN, CN, MP}). The global existence of smooth solution for the 3D
Euler equations remains outstandingly open. Therefore, the question of the global well-posedness of smooth solution for  the multidimensional inviscid Boussinesq system is challenging.

To understand the turbulence phenomena in fluid mechanics, the study of weak solutions to the hydrodynamic equations becomes more and more important in the past several decades. For instance, the physicist and chemist Lars Onsager \cite{Onsa} proposed a celebrated  conjecture for the Euler equations.
\begin{con}[Onsager's conjecture]  Consider weak solutions of the  Euler equations in $C_tC^{\beta}(\TTT^3)$.
 \begin{enumerate}
    \item $[\rm{Positive\,\,\,part}]$ If $\beta>1/3$,  the energy $\int_{\mathbb T^3} |v|^2 \dd x$ is conserved in time.
    \item $[\rm{Negative\,\,\,part}]$ If  $\beta\le 1/3$, there exist at least one solution that does not conserve energy.
\end{enumerate}
 \end{con}
 The positive part of the conjecture was proven by  Constantin, E and Titi in~\cite{CET}. So far, the sharpest result on energy conservation of weak solutions was obtained by  {Cheskidov} et al.. More specifically, they proved that  the kinetic energy
is conserved   in the class $L^3_t B^{1/3}_{3, c(\NN)}\cap C_t L^2_x$ in \cite{CCFS}.

 The negative part of the conjecture has been solved through several decades of unremitting efforts of many researchers.  The first result in this direction was obtained by Scheffer \cite{S93}, in which the author constructed  a nontrivial weak solution with compact support in space-time $\R^2\times\R$. Later, Shnirelman \cite{Sh97} gave a simpler and more transparent example of a weak solution with  compact support in time in the class of $L^2(\TTT^2\times\R)$. These weak solutions are  discontinuous, unbounded $L^2$-fields. Therefore, these examples cannot show whether there exists a weak solution with monotonically decreasing kinetic energy or not. The first result on the existence of a dissipative solution   was  demonstrated by  Shnirelman \cite{Sh00} in the class $L^2(\TTT^3\times[T_1, T_2])$. After nearly a decade,  De Lellis and Sz\'{e}kelyhidi applied  either the
Baire category theorem (e.g., \cite{DM97}) or the general framework of convex integration (e.g., \cite{MS,S06}) to the multidimensional incompressible Euler equations and recover Scheffer and Shnirelman's counterexamples  with bounded velocity and pressure in the seminal paper \cite{DS09}. However, the Baire category approach  does not work for obtaining continuous solutions. As a result, the method of convex integration has been developed as an alternative in the follow-up research on Onsager's conjecture. The authors in \cite{DS13} abandoned this  essentially ``soft'' functional analytic technique in favor of a more ``hard'' PDE approach
based on Schauder's estimates and oscillatory integrals to show the existence of continuous periodic dissipative solutions. Subsequently,  
they \cite{DS14} constructed periodic
weak solutions in the H\"{o}lder space $C^{1/10-\epsilon}(\TTT^3\times[0,1])$,  which dissipate the total kinetic energy with velocity.  Choffrut \cite{Ch13} extended this result of  \cite{DS14} by establishing optimal h-principles in two spatial dimensions. Isett \cite{Ise17} constructed   H\"{o}lder continuous periodic solutions with compact support in time in $C^{1/5-\epsilon}_{t,x}$  and Buckmaster et al. \cite{BDIS15} gave a shorter proof that includes a result on the existence of anomalous dissipative solutions in $C^{1/5-\epsilon}_{t,x}$. Furthermore,
Buckmaster \cite{B15} proved the existence of non-trivial H\"{o}lder continuous solutions belonging to $C^{1/3-\epsilon}$ for almost
every time. Making use of more subtle modifications, the authors in \cite{BDS} show the existence of continuous periodic weak dissipative solutions in $L^1_t(C^{1/3-\epsilon})$.  Using  ``Mikado flows'' introduced in \cite{DaS17}  as an alternative to Beltrami flows to build the waves in the construction and developing  a new ``gluing approximation'' technique, Isett \cite{Ise18} constructed weak solutions in $C_t(C^{1/3-\epsilon})$ that have  nonempty, compact temporal support and thus fail to conserve the total kinetic energy, which finally proves the negative part of  Onsager's conjecture up to $\beta=\frac{1}{3}$.


The study on the Onsager's conjecture of the incompressible Euler equations motivates us to consider whether the weak solutions of the inviscid Boussinesq system \eqref{e:B} satisfy \eqref{vL2} and whether  temperature conserves $L^p$-norm or not.  Assume that $(v,\theta)$ is a smooth solution of the system \eqref{e:B} on $\TTT^3\times [0,T]$ with smooth initial data $(\vin, \tin)$, one easily obtains that
\begin{equation}\label{vL2}
\|v( t)\|^2_{L^2(\TTT^3)}=\|\vin\|^2_{L^2(\TTT^3)} +\int_0^t\int_{\TTT^3}\theta(x,s) v_3(x,s)\dd x,\quad \forall s\in [0,T],
\end{equation}
and for each $p\ge 1$,
 \begin{equation}\label{tLp}
\|\theta( t)\|_{L^p(\TTT^3)}=\|\tin\|_{L^p(\TTT^3)}, \quad \forall t\in [0,T].
\end{equation}
 Similar to Onsager's conjecture of the Euler equations, a natural question is
\vskip 2mm
(Q): \emph{What is the threshold H\"{o}lder exponent for a weak solution $(v ,\theta)\in C_tC^{\beta}_x$ of the system \eqref{e:B} to satisfy \eqref{vL2} and \eqref{tLp}?}
\vskip 2mm
  To the best of our knowledge, Tao and Zhang \cite{TZ17, TZ18} firstly considered this problem.
 They constructed the H\"{o}lder continuous solution $\small{(v, \theta)\in C^{1/28-\epsilon}_{t,x}\times  C^{1/25-\epsilon}_{t,x}}$ with space-time compact support on $\R^2\times\R$ and  pointed out that similar results also hold on $\R^3\times\R$ in  \cite{TZ17}.  These results imply that there are multiple solutions in some H\"{o}lder spaces starting from the identically zero flow, and so breaking \eqref{vL2}--\eqref{tLp}.  Subsequently, the authors \cite{TZ18}  showed the existence of  continuous solutions with prescribed kinetic energy and some other prescribed properties on torus $\TTT^3$.  

In this paper, we  address this problem (Q) on  tours  $\TTT^3$. From the proof in \cite{CET} by Constantin, E and Titi, one can directly show that $\theta$ conserves $L^2$-norm if $(v, \theta)\in C_tC^{\beta}$ with $\beta >\frac{1}{3}$. A rough viewpoint of derivative allocation
 \[\int_{\TTT^3}(v\cdot\nabla \theta)\theta \dd x\approx  \int_{\TTT^3}\nabla^{\frac{1}{3}}v(\nabla^{\frac{1}{3}}\theta)^2 \dd x\]
 may be used to explain the index $\beta>\frac{1}{3}$. This motivates us to expect that $\theta$ could conserve $L^p$-norm  if $(v, \theta)\in C_tC^{\beta}_x(\TTT^3)$ with $\beta >\frac{1}{p+1}$. Interestingly, we show that the threshold exponent for question (Q) is still $\frac{1}{3}$ even for the $L^p$-norm conservation of temperature, which indicates that this rough viewpoint above is not always effective. More precisely,
 for any $C^{\infty}$ initial data, we show that the corresponding weak solution $(v, \theta)\in C_tC^{\beta}_x(\TTT^3)$ for $\beta>\frac{1}{3}$ obeys  \eqref{vL2}--\eqref{tLp}, while there exist infinitely many weak solutions in  $ C_tC^{\beta}_x(\TTT^3)$ for $\beta<\frac{1}{3}$ breaking \eqref{vL2}--\eqref{tLp}. Here, the weak solutions are defined as follows.
\begin{defn}\label{def} Let $(\vin, \tin)\in  C^{\infty}(\TTT^3)$ and $\vin$ be divergence-free in the sense of distributions. We say that a pair $(v, \theta)\in {C(\mathbb \TTT^3 \times [0,T])}$ is a \emph{weak solution} on $[0,T]$ to the system~\eqref{e:B} with initial data $(\vin, \tin) $ if \begin{itemize}
     \item [(1)] For a.e. $t\in [0,T]$, $v(\cdot, t)$ is divergence-free in the sense of distributions;
     \item [(2)]For all divergence-free test functions $\phi\in  C^\infty_c\big([0,T); \big(C^{\infty}(\TTT^3)\big)^3\big)  $ ,
     \begin{align}\nonumber
-\int_{\mathbb T^3} \vin \phi(x,0)\dd x-\int_0^T \int_{\mathbb T^3} \del_s\,\phi\cdot v\dd x \dd s=\int_0^T \int_{\mathbb T^3} (\nabla \phi : v\otimes v+\theta e_3\cdot\phi)  \dd x \dd s,
\end{align}
and for all test functions $\varphi\in C^\infty_c([0,T); C^{\infty}(\TTT^3)) $,
\begin{align}\label{weaktheta}
 -\int_{\mathbb T^3} \tin\varphi(x,0)\dd x-\int_0^T \int_{\mathbb T^3} \del_s\,\varphi \theta\dd x \dd s=\int_0^T \int_{\mathbb T^3} v\cdot\nabla\varphi\,\theta\dd x \dd t.
\end{align}
   \end{itemize}
\end{defn}
Now we are in position to state our result.

\begin{thm}\label{main00}
Consider weak solutions $(v,\theta)$ of the system \eqref{e:B} in $C_tC^{\beta}(\TTT^3)$.
 \begin{enumerate}
    \item $[\rm{Positive\,\,\,part}]$ If $\beta>1/3$,  then $\|\theta(t)\|_{L^p(\TTT^3)}(1\leq p\leq \infty)$ is conserved in time.
    \item $[\rm{Negative\,\,\,part}]$ If  $\beta<1/3$,  there exist many  solutions that do not conserve $L^p$-norm.
\end{enumerate}
\end{thm}
\begin{rem}Following the work of Constantin, E and Titi in \cite{CET}, one immediately deduces that the equality \eqref{vL2} holds and $\theta$ conserves $L^2$-norm  for weak solutions  $(v, \theta)\in C_tC^{\beta}(\TTT^3)$ with $\beta>\frac{1}{3}$.  However, for $1\le p\le\infty$, the critical H\"{o}lder  exponent  that ensures $\theta$  the $L^p$-norm conservation is not straightforward. Motivated by \cite{CCFS}, we show that  $(v, \theta)$ satisfies \eqref{vL2}--\eqref{tLp} when it belongs to some Besov spaces that are weaker than $C_tC^{\beta}(\TTT^3)$ with $\beta>\frac{1}{3}$.
\end{rem}
The proof of  Theorem \ref{main00} can be concluded by Theorem \ref{main} and Theorem \ref{t:main0}.
\begin{thm}\label{main} Let $(v, \theta)$ be a weak solution of the system \eqref{e:B}. Suppose that $1\le p\le\infty$, $\tfrac{1}{p}+\tfrac{1}{p'}=1$, $v\in L^{\infty}([0,T], B^{\frac{1}{3}}_{\infty, c(\NN)}(\TTT^3))$ and
\begin{itemize}
  \item $\theta\in C([0,T], L^p(\TTT^3))\cap L^2([0,T];B^{\frac{1}{3}}_{p', c(\NN)}(\TTT^3))\,\, \text{for} \,\,1\le p<2,$
\item $\theta\in C([0,T], L^p(\TTT^3))\cap L^p([0,T]; B^{\frac{1}{3}}_{p, c(\NN)}(\TTT^3))\,\, \text{for}\,\, p\ge2,$
\end{itemize}
then the weak solution  $(v,\theta)$ satisfies \eqref{vL2} and \eqref{tLp}.
\end{thm}

The above positive result immediately shows that the weak solutions $(v, \theta)$ in $C_tC^{\beta}_x(\TTT^3)$, $\beta>\tfrac{1}{3},$ satisfy \eqref{vL2} and \eqref{tLp}.  The negative part of Theorem \ref{main00} can be concluded by the following result.
\begin{thm}\label{t:main0}
For any given initial data $(\vin, \tin)\in C^\infty(\TTT^3)$ and  $0<\beta<\tfrac{1}{3}$,  there exist infinitely many weak solutions $(v, \theta)$ on $[0,T]$ with $(v, \theta)|_{t=0}=(\vin, \tin)$  such that  $(v, \theta)\in C^{\beta}_{t,x}$ and violates \eqref{vL2}--\eqref{tLp}.
\end{thm}
\begin{rem}For $\beta<\frac{1}{3}$, Theorem \ref{t:main0} immediately show the non-uniqueness of a weak solution $(v , \theta)$ of the
system~\eqref{e:B} in the class $C^{\beta}_{t,x}$ as long as $(v , \theta)$ has at least one interval of regularity.
\end{rem}
Theorem \ref{t:main0} can be summarized in the following proposition.
\begin{prop}\label{t:main} 
Let  $(v^{(1)}, \theta^{(1)})\in C([0, T_1]; C^\infty(\TTT^3))$ \footnote{Here and throughout the paper, for three dimensional space-time vector function $f$ and  space-time scalar function $g$, we denote $(f , g)\in C(I; (X)^3\times X)$ by $(f, g)\in C(I;X)$. For three dimensional space vector function $u$ and space scalar function $v$, we denote $(u , v)\in  (X)^3\times X$ by $(u, v)\in X$.} and $(v^{(2)}, \theta^{(2)})\in C([0, T_2]; C^\infty(\TTT^3))$ solve the system \eqref{e:B} with mean-free initial data $(v^{(1)}(0, x)$ $ \theta^{(1)}(0,x))$, $(v^{(2)}(0, x), \theta^{(2)}(0, x))\in C^\infty(\TTT^3)$, respectively.  Fixed $ \widetilde{T}\le \tfrac{1}{4}\min\{T_1, T_2\}$ and $0<\beta<\frac{1}{3}$, there exists a weak solution $(v, \theta)\in C^{\beta}_{t,x}$ of the Cauchy problem for the system \eqref{e:B}  satisfying
\begin{align}\label{TW}
(v, \theta)\equiv(v^{(1)}, \theta^{(1)})\,\,\text{\rm on}\,\,[0, 2 \widetilde{T}],\quad \text{ and}\quad(v, \theta)\equiv(v^{(2)}, \theta^{(2)})\,\,\text{\rm on}\,\,[3\widetilde{T}, T_2].
\end{align}
\end{prop}
\begin{rem} For the sake of simplicity, we give this result under the assumption of mean-free initial data. In fact, Proposition \ref{t:main} can be proved without the mean-free condition.  For general  initial data $(v^{(1)}(0, x)$ $ \theta^{(1)}(0,x))$, $(v^{(2)}(0, x), \theta^{(2)}(0, x))\in C^\infty(\TTT^3)$, by transformation
\begin{align*}
(\widetilde{v}^{(1)}, \widetilde{\theta}^{(1)})=\big(v^{(1)}-\int_{\TTT^3}v^{(1)}(0, x)\dd x, \,\,\theta^{(1)}-\int_{\TTT^3}\theta^{(1)}(0,x)\dd x\big),\\
(\widetilde{v}^{(2)}, \widetilde{\theta}^{(2)})=\big(v^{(2)}-\int_{\TTT^3}v^{(2)}(0, x)\dd x, \,\,\theta^{(2)}-\int_{\TTT^3}\theta^{(2)}(0,x)\dd x\big),
\end{align*}
one immediately shows from the proof of Proposition \ref{t:main} that for $0<\beta<\frac{1}{3}$, {there exists a weak solution $({v}, {\theta})\in C^{\beta}_{t,x}$ of the system \eqref{e:B} satisfying \eqref{TW}}.
\end{rem}
\noindent \textbf{Proposition \ref{t:main} implies Theorem \ref{t:main0}.}
\vskip 1mm
 Given initial data $(\vin, \tin)\in C^\infty(\TTT^3)$, there exists a unique local-in-time smooth solution $(v^{(1)},\theta^{(1)})\in C([0,T_1];C^\infty(\TTT^3))$ by the local well-posedness of the Cauchy problem for the Boussinesq equations \eqref{e:B} (e.g. \cite{MP}). Let $(v^{(2)}, \theta^{(2)})$ be the shear flows such that
\begin{equation}\label{shear flow}
v^{(2)}(t,x)=(0,0,Af(x_1)t), \quad\theta^{(2)}(t,x)=Af(x_1),
\end{equation}
      where $A$ is a constant and $f(x_1)$ is a smooth function with zero mean. $(v^{(2)}, \theta^{(2)})$ solves the equations \eqref{e:B} on $\TTT^3\times[0,\infty)$. Therefore, we obtain from Proposition \ref{t:main} that for fixed $\beta<\tfrac{1}{3}$ and $T>0$, there exists a weak solution $(v, \theta)\in C^\beta_{t,x}$ of the system \eqref{e:B} on $[0,T]$ such that
      \begin{align*}
(v, \theta)\equiv(v^{(1)}, \theta^{(1)})\,\,\text{on}\,\,\big[0, \tfrac{1}{2}\min\{T_1, T\}\big],\quad \text{and}\quad(v, \theta)\equiv(v^{(2)}, \theta^{(2)})\,\,\text{on}\,\,\big[\tfrac{3}{4}\min\{T_1, T\}, T\big].
\end{align*}
Therefore, the solution $(v, \theta)$ satisfies $(v,\theta)|_{t=0}=(\vin, \tin)$, and for every $t\in\big[\tfrac{3}{4}\min\{T_1, T\}, T\big]$,
\begin{align*}
&\|v(t)\|^2_{L^2(\TTT^3)}=A^2 t^2\int_{\TTT^3}f^2(x_1)\dd x, \quad\|\theta(t)\|_{L^p(\TTT^3)}=A\|f(x_1)\|_{L^p(\TTT^3)},
\end{align*}
and
\begin{align*}
&\int_{\TTT^3}v_3(t,x)\theta(t,x)\dd x=A^2 t \int_{\TTT^3}f^2(x_1)\dd x.
\end{align*}
These equalities immediately show that there exist infinitely many $A$ and $f$ such that \eqref{vL2} and \eqref{tLp} do not hold on $\big[\tfrac{3}{4}\min\{T_1, T\}, T\big]$. Hence,
we conclude Theorem \ref{t:main0}.
\vskip 1mm
\noindent \textbf{Some remarks.}
\begin{enumerate}

  \item 
  Tao and Zhang  found an interesting phenomenon that there exist some weak solutions $(v,\theta)\in C(\TTT^3\times[0,1])$ such that for the small initial velocity, the oscillation of velocity after some time could be as large as possible for enough heat in the system~\eqref{e:B}, see \cite[Theorem 1.4]{TZ18}. In fact, for any $0<\beta<\frac{1}{3}$, we can construct many weak solutions $(v,\theta)\in C^\beta_{t,x}$ which exhibit such behavior.   More precisely,  taking $f(x_1)=\cos(2\pi x_1)$ in~\eqref{shear flow},  for a given initial data $(\vin, \tin)\in C^\infty(\TTT^3)$, $0<\beta<\tfrac{1}{3}$ and $T>0$, we can construct a weak solution $(v, \theta)\in C^\beta_{t,x}$ such that $(v,\theta)|_{t=0}=(\vin, \tin)$ and
\begin{equation}\label{v}
v(t,x)=(0,0, A\cos(2\pi x_1)t), \,\,\theta(t,x)=A\cos(2\pi x_1), \quad\forall t\in \big[\tfrac{T}{2}, T\big].
\end{equation}
Particularly, when  $\|(\vin, \tin)\|_{H^3(\TTT^3)}\ll 1$,  taking $A$ large enough, then the weak solution $(v, \theta)$ in the class $C^\beta_{t,x}$  of the system \eqref{e:B} demonstrates the same behavior as in  \cite[Theorem 1.4]{TZ18}.  
\item Interestingly, for $0<\beta<\frac{1}{3}$, weak solutions of the system \eqref{e:B} displaying the opposite  behavior in the class $C^{\beta}_{t,x}$  also exist.  Indeed, in the process of proving Theorem  \ref{t:main0}, choosing
     $$v^{(1)}=(0,0, A\cos(2\pi x_1)t), \,\,\theta^{(1)}=A\cos(2\pi x_1),\,\, (v^{(2)}, \theta ^{(2)})=0,$$
   then  for $0<\beta<\frac{1}{3}$ and $\widetilde{T}\ll 1$,  we can construct a global weak solution in $C^{\beta}_{t,x}$ satisfying
    \begin{align*}
    &v(t,x)=(0,0, A\cos(2\pi x_1)t), \,\,\theta(t,x)=A\cos(2\pi x_1) \,\,\text{on}\,\,[0, 2{\widetilde{T}}],\\
    &v(t,x)=\theta(t,x)=0\,\,\text{on}\,\,[ 3{\widetilde{T}}, \infty).
    \end{align*}
 These properties show that  even though the kinetic energy and the oscillation of velocity are large at initial time, the weak solution could dissipate the kinetic energy of this system as fast as possible.
\end{enumerate}

\noindent Now let us introduce main ideas in the proofs of Theorem \ref{main} and Theorem \ref{t:main0}.
For the proof of Theorem \ref{main}, Poisson's summation formula allows us to apply the argument in \cite{CCFS} to the framework of the torus
 and and we introduce a smooth approximation of $|z|$ to show $L^p$-norm conservation of temperature  for $1\le p\le \infty$. Theorem~\ref{t:main0}  can be reduced  to
a suitable iterative procedure, summarized in Proposition~\ref{p:main-prop} below.
In \cite{TZ18}, the authors employed a multi-step iteration scheme. More precisely, due to  no known analogue for Beltrami flows,  they decomposed the stress error into blocks by the geometric lemma and  removed one component of stress error in each step by putting  several plane waves which oscillate along the same direction with different frequency. In this paper, we build waves based on intermittent cuboid flows, which are constructed  by collecting these flows supported in disjoint cuboids that point in multiple directions and essentially are Mikado flows given in \cite{DaS17}.
Different from the flows built in \cite{TZ18},  a intermittent cuboid flow on its own does not produce unacceptable error terms over a sufficiently long time scale, even when composed with the Lagrangian flow of the slow velocity field. To fully utilize this  superiority  and avoid the interference of waves with different oscillation directions, using the gluing approximation technique introduced in \cite{Ise18}, we  construct a gluing Boussinesq flow $(\vv_q, \tb_q, \pp_q, \RRR_q, \MMM_q)$, which satisfy the same  estimates of  $(v_q, b_q, p_q, \RR_q, \MM_q)$ up to a small error. Benefiting from the stresses $(\RRR_q, \MMM_q)$ supported in disjoint short time intervals,
we can decompose $(\RRR_q, \MMM_q)$ as a sum $(\sum_i\RRR_{q,i}, \sum_i\MMM_{q,i})$ and   cancel each $(\RRR_{q,i}, \MMM_{q,i})$ by adopting a single intermittent cuboid flow that  is Lie-advected by $\vv_q$ and prevent the distinct flows from  interacting with each other. Therefore, the interference between distinct deformed intermittent cuboid flows built is controlled over a sufficiently long time scale. This is the core of improving the regularity of weak solution from continuous space  to $C^{1/3-}$ H\"{o}lder space.

\vskip 2mm
\noindent \textbf{Organization of the paper.}
\vskip 1mm
In Section \ref{Preliminaries}, we recall the  Littlewood-Paley theory, the definition of Besov spaces on torus and some related preliminary lemmas. Section \ref{proof of positive} contains the proof of Theorem \ref{main}. In Section \ref{proof of negative}, we focus on proving Proposition \ref{t:main}. In  Section 4.1, we reduce Proposition~\ref{t:main} to Proposition 4.1. The proof of Proposition 4.1 occupies Sections~4.2-4.5. To be more precise, in Section 4.2, we implement mollification and gluing procedure to obtain the glued velocity and temperature, of which the corresponding Reynolds and temperature stresses are supported in disjoint time intervals. In Section 4.3, we introduce the deformed  intermittent cuboid flows to construct the perturbation of the glued velocity and temperature. We establish the H\"{o}lder estimates for the perturbation, the new Reynolds and temperature stresses in Section 4.3 and Section 4.4, respectively. Collecting these estimates together,  we verify that the inductive bounds stated in Proposition 4.1 are propagated from step  $q$ to $q+1$ and thus complete the proof of Proposition~4.1 in Section 4.5.
\section{Preliminaries}\label{Preliminaries}
For the convenience of the readers,  we review briefly the so-called Littlewood-Paley  theory introduced  in \cite{Ba, C, MWZ12}. Suppose that $\chi$ and $\varphi$ are two smooth radial functions with values in $[0,1]$, satisfying
\begin{align*}
&\text{supp}\,\chi\subset\big\{\xi\in\mathbb{R}^3\big{|}|\xi|\le \tfrac{4}{3}\big\}; \qquad  \quad \chi(\xi)\equiv1 \quad \text{for}\,\,\, |\xi|\le\tfrac{3}{4};\\
&\text{supp}\,\varphi\subset\big\{\xi\in\mathbb{R}^3\big{|}\tfrac{3}{4}\le|\xi|\le\tfrac{8}{3}\big\}; \quad\varphi(\xi)\equiv1 \quad \text{for}\,\,\,\tfrac{6}{7}\le |\xi|\le \tfrac{12}{7},
\end{align*}
and
\[\chi(\xi)+\sum_{j\ge 0}\varphi(2^{-j}\xi)=1, \quad \forall\xi\in \mathbb{R}^3.\]

Let $D(\TTT^3)$ be the set of all infinitely differentiable functions on $\TTT^3$. Then $D'(\TTT^3)$ is defined to be the topological dual of $D(\TTT^3)$. Formally, every distribution $f\in D'(\TTT^3)$ can be represented by its Fourier series
\begin{align*}
f=\sum_{k\in\ZZ^d}\widehat{f}(k)e^{2\pi \ii k\cdot x} \quad(\text {convergence in } D'(\TTT^3) ).
\end{align*}
For $j\le -2,$ we define $\Delta_{j}f=0$.
\begin{align*}
\Delta_{-1}f=\sum_{k\in\ZZ^d}\chi(k)\widehat{f}(k)e^{ 2\pi \ii k\cdot x},
\end{align*}
and for $j\ge 0$,
\begin{align*}
\Delta_{j}f=\sum_{k\in\ZZ^d}\varphi(2^{-j}k)\widehat{f}(k)e^{ 2\pi \ii k\cdot x}:=\sum_{k\in\ZZ^d}\varphi_j(k)\widehat{f}(k)e^{ 2\pi\ii k\cdot x}.
\end{align*}
Furthermore, let  nonhomogeneous low-frequency cut-off operator $S_Q$ be defined by
\begin{equation}\label{defSQf}
S_Q f=\sum_{j=-1}^{Q-1}\Delta_j f=\sum_{k\in\ZZ^d}\chi(2^{-Q}k)\widehat{f}(k)e^{2\pi\ii  k\cdot x}.
\end{equation}
Let us denote by $\mathcal{F}^{-1}$ the inverse Fourier transform on $\R^3$. In the following, we denote $h=\mathcal{F}^{-1}\varphi$ and $\widetilde{h}=\mathcal{F}^{-1}\chi$. In fact, one can deduce that for $f\in D'(\TTT^3)$,
\begin{align}
&\Delta_{-1}u=\int_{\mathbb{R}^3}\widetilde{h}(y)f(x-y)\dd y,\label{D1f}\\
&\Delta_j f=2^{3j}\int_{\mathbb{R}^3} h(2^jy)f(x-y)\,{\rm d}y,\qquad\forall j\ge0,\label{Djf}\\
&S_Q f=\sum_{j=-1}^{Q-1}\Delta_j u=2^{3Q}\int_{\mathbb{R}^3}\widetilde{h}(2^Qy)f(x-y)\,{\rm d}y \label{SQ}.
\end{align}
Indeed, we obtain from \eqref{defSQf} that
\begin{align*}
S_Q f=&\int_{\TTT^3}f(x-y)\sum_{k\in\ZZ^d}\chi(2^{-Q}k)e^{2\pi\ii  k\cdot y}\dd y\\
=&\int_{\TTT^3}f(x-y)\sum_{m\in\ZZ^d}\mathcal{F}^{-1}\big(\chi(2^{-Q}\cdot)\big)(y+m)\dd y,
\end{align*}
where the second equality follows by Poisson's summation formula. Whence,
\begin{align*}
S_Q f=&\sum_{m\in\ZZ^d}2^{3Q}\int_{\TTT^3}\widetilde{h}(2^Q (y+m))f(x-y)\dd y\\
=&\sum_{m\in\ZZ^d}2^{3Q}\int_{\TTT^3+m}\widetilde{h}(2^Q y)f(x-y-m)\dd y\\
=&\sum_{m\in\ZZ^d}2^{3Q}\int_{\TTT^3+m}\widetilde{h}(2^Q y)f(x-y)\dd y\\
=&2^{3Q}\int_{\mathbb{R}^3}\widetilde{h}(2^Q y)f(x-y)\dd y.
\end{align*}
Therefore, we deduce \eqref{SQ}. In the same way as in deriving \eqref{SQ}, one can obtain \eqref{D1f} and \eqref{Djf}. For any functions whose Fourier transforms are supported in an annulus, there is a nice property
called Bernstein's inequality in $\mathbb{R}^3$ \cite{Ba, MWZ12}. By making use of \eqref{D1f}--\eqref{SQ}, one can show the  Bernstein's inequality for  $\TTT^3$-periodic functions.
\begin{prop}[Bernstein's inequality on torus]\label{Bernstein}Let $\mathcal{C}$ be an annulus and $\mathcal{B}$ be a ball, $1\le p\le q\le\infty$. Then there exists a constant $C$ such that for any $k\in \NN_{+}$ and any $\TTT^3$-periodic function $u\in L^p(\TTT^3)$, we have
\begin{align*}
\sup_{|\alpha|=k}\|\partial^\alpha u\|_{L^q(\TTT^3)}&\le C^{k+1}\lambda^{k+3(\frac{1}{p}-\frac{1}{q})}\|u\|_{L^p(\TTT^3)},\qquad\text{\rm supp}\,\,\hat{u}\subset \lambda \mathcal{B},\\
C^{-(k+1)}\lambda^k\|u\|_{L^p(\TTT^3)}&\le\sup_{|\alpha|=k}\|\partial^\alpha u\|_{L^p(\TTT^3)}\le C^{k+1}\lambda^{k}\|u\|_{L^p(\TTT^3)},\qquad\text{\rm supp}\,\,\hat{u}\subset \lambda \mathcal{C}.
\end{align*}
\end{prop}
\begin{defn}[Nonhomogeneous Besov spaces $B^s_{p,q}(\TTT^3)$]\label{def.Be}
Let $s\in\mathbb{R}$ and $1\le p,q\le\infty$.  The nonhomogeneous Besov space $ B^s_{p,q}(\TTT^3)$ consists of all $\TTT^3$-periodic functions $u\in D'(\TTT^3)$  such that
\begin{equation}\nonumber
\|u\|_{ B^s_{p,q}(\TTT^3)}\overset{\text{def}}{=}\Big{\|}\left(2^{js}\|\Delta_j u\|_{L^p(\TTT^3)}\right)_{j\in\ZZ}\Big{\|}_{\ell^q(\ZZ)}<\infty.
\end{equation}
\end{defn}

\begin{defn}[Besov spaces $B^{s}_{p, c(\NN)}(\TTT^3)$] For $s\in \mathbb{R}$ and $1\le p\le\infty$, we define the space $B^{s}_{p, c(\NN)}$ as follows:
\[B^{s}_{p, c(\NN)}:=\big\{u\in B^{s}_{p, \infty}(\TTT^3): 2^{sj}\|\Delta_j u\|_{L^p(\TTT^3)}\to0,\,\, {\rm as}\,\, j\to\infty\big\}.\]
\end{defn}
\noindent Employing \eqref{D1f}--\eqref{SQ} and the commutator estimate \cite[Lemma 2.100]{Ba}, one immediately show the following result on $\TTT^3$.
\begin{lem}[A priori estimates for the transport equations in Besov spaces]\label{trans}
Let $1\le p\le p_1\le\infty$ and $\sigma\in [-1-d\min\{\frac{1}{p_1}, 1-\frac{1}{p}\}, 1+\frac{d}{p_1}].$ Let $v$ be a divergence-free vector field such that $\nabla v\in L^1_T( B^{\frac{d}{p_1}}_{p_1, \infty}(\TTT^3)\cap L^\infty(\TTT^3))$. There exists a constant $C$ depending on $p, p_1,\sigma$ such that for all solutions $f\in L^\infty_T(B^{\sigma}_{p,\infty}(\TTT^3))$ of the transport equation
\[\partial_t f+v\cdot\nabla f=g,\,\,\,\,f(0,x)=f_0(x),\]
with initial data $f_0\in   B^{\sigma}_{p,\infty}(\TTT^3)$ and $g\in L^1_{T}(B^{\sigma}_{p,\infty}(\TTT^3))$, we have, for $t\in[0, T]$,
\begin{align*}
\|f\|_{ L^\infty_T(B^{\sigma}_{p,\infty}(\TTT^3))}\le e^{CV_{p_1}(t)}\Big(\|f_0\|_{  B^{\sigma}_{p,\infty}(\TTT^3)}+\int_0^t  e^{-CV_{p_1}(\tau)}\|g(\tau)\|_{ B^{\sigma}_{p,\infty}(\TTT^3)}\dd\tau\Big),
\end{align*}
where $V_{p_1}(t)=\int_0^t\|\nabla v\|_{( B^{\frac{d}{p_1}}_{p_1, \infty}(\TTT^3)\cap L^\infty(\TTT^3))}\dd s$.
\end{lem}
\begin{defn}[The operators $\mathcal{R}$ and $\mathcal R_{\vex}$]\label{def.R}
 For vector field $u\in C^{\infty}(\TTT^3, \R^3)$, the operator $\mathcal{R}$ introduced in \cite{BDIS15} is defined by
\begin{align*}
    \mathcal R u = -(-\Delta)^{-1} (\nabla u + \nabla u^\TT ) - \frac12(-\Delta)^{-2} \nabla^2 \nabla\cdot u  + \frac12 (-\Delta)^{-1} (\nabla\cdot u ){\rm Id}_{3\times 3}.
\end{align*}
It is a matrix-valued right inverse of the divergence operator for mean-free vector fields, in the sense that
\[ \div\mathcal R u= u - \dashint_{\mathbb T^3} u. \]
In addition, $\mathcal Ru$ is traceless and symmetric.

For mean-free scalar function $f\in C^\infty(\TTT^3, \R)$, we define the operator $R_{\vex}$ by
$$\mathcal R_{\vex} f:={\nabla}{\Delta}^{-1} f.$$
It is a vector-valued right inverse of the divergence operator for mean-free scalar  functions, in the sense that
\[ \div\mathcal R_{\vex} f= f - \dashint_{\mathbb T^3} f. \]
\end{defn}
\begin{lem}[\cite{BDSV, TZ18}]\label{l:non-stationary-phase} {If $a\in C^\infty(\mathbb T^3;\mathbb R^3)$}, $b\in C^\infty(\mathbb T^3;\mathbb R)$ and $\Phi\in  C^\infty (\mathbb T^3 ; \mathbb R^3)$ satisfying  $|\nabla \Phi|  \sim 1$, then
\begin{align*}
\|\mathcal R( a \ee^{\ii k\cdot \Phi}) \|_{C^\alpha}
    \lesssim  \frac{\|a\|_{C^0}}{|k|^{1-\alpha}} + \frac{\|a\|_{C^{N+\alpha}} + \|a\|_{C^0} \|\Phi\|_{C^{N+\alpha}}} {|k|^{N-\alpha}},\\
    \|\mathcal R_{\vex}( b \ee^{\ii k\cdot \Phi}) \|_{C^\alpha}
    \lesssim  \frac{\|b\|_{C^0}}{|k|^{1-\alpha}} + \frac{\| b\|_{C^{N+\alpha}} + \|b\|_{C^0} \|\Phi\|_{C^{N+\alpha}}} {|k|^{N-\alpha}}.
\end{align*}
\end{lem}
\vskip 3mm
\noindent\textbf{{Notation}}\,\, In the following, $\alpha\in(0,1)$, $N\in \NN$ and  $\sigma$ is a multi-index. Give a function $f$ on $\TTT^3\times[0,T]$, we denote the supremum norm by
\[\|f\|_{0}:=\sup_{\TTT^3\times[0,T]}|f|.\]
For the space derivatives $\partial^{\sigma}$, we define the H\"{o}lder seminorms by
\begin{align*}
[f]_{N}:=\max_{|\sigma|=N}\|\partial^{\sigma}f\|_0, \quad [f]_{N+\alpha}:=\max_{|\sigma|=N}\sup_{\substack{x\neq y\in \TTT^3\\t\in[0,T]}}\frac{|(\partial^{\sigma}f)(x,t)-(\partial^{\sigma}f)(x,t)|}{|x-y|^{\alpha}},
\end{align*}
and define the H\"{o}lder norms by
\begin{align*}
\|f\|_{N}:=\sum_{j=0}^N[f]_j, \quad [f]_{N+\alpha}:=\|f\|_N+[f]_{N+\alpha}.
\end{align*}

Let  $\widetilde\phi\in C^\infty(\R;[0,\infty))$ be an even non-negative bump function satisfying $\supp \,\,\widetilde\phi\subset[-1,1]$
For each $\epsilon>0$, we define  two sequences of  mollifiers as follows:
\begin{align}
     \phi(t)
        &\coloneq \frac{\widetilde\phi(t)}{\int_{[-1,1]}\widetilde\phi (\tau)\dd \tau },
    \quad\,\,\, \phi_{\epsilon}(t)
        \coloneq \frac1{\epsilon} \phi\left(\frac t\epsilon\right)\label{e:defn-mollifier-t},
     \\
    \psi(x)
        &\coloneq\frac{\widetilde\phi(|x|)}{\int_{B_1(0)}\widetilde\psi (|y|)\dd y},
    \quad\psi_\epsilon(x)
            \coloneq \frac1{\epsilon^3} \psi\left(\frac{x}\epsilon\right). \label{e:defn-mollifier-x}
\end{align}
Throughout this paper, we use the notation $x\lesssim y$ to denote $x\le Cy$ for a universal constant that may change from line to line. The symbol $\lesssim_N$ will imply that the constant in the inequality depends on $N$. The notations  $x\ll y$  and $x\gg y$ mean that $x$ is much smaller than $y$ and $x$ is much larger than $y$, respectively.
\section{The proof of Theorem \ref{main}}\label{proof of positive}
In this section, we prove that the weak solutions $(v, \theta)$ of the system \eqref{e:B} satisfy \eqref{vL2} and \eqref{tLp} under the assumptions in Theorem \ref{main}. Following the arguments in \cite{CCFS},  one can directly show that $(v ,\theta)$ satisfies \eqref{vL2}. Next, we are focused on proving that such $(v ,\theta)$ also obeys~\eqref{tLp}.

Let $S_Q$ be the frequency localization operator defined in \eqref{SQ} and $\phi_{\epsilon}$ be a standard mollifier defined in  \eqref{e:defn-mollifier-t}.  Fixed  $t\in (0, T]$,  let $\epsilon<\eta<\frac{t}{2}$, we denote
$$(S_\epsilon f)(\tau)=\int_0^T \phi_{\epsilon}(\tau-s)f(s)\dd s.$$
Now, for  $p\in [2,\infty)$, putting  the test function
$$\varphi(x,s)=\int_{\eta}^{t-\eta}\phi_{\epsilon}(\tau-s)S_Q(S_\epsilon S_Q\theta|S_\epsilon S_Q\theta|^{p-2})(x,\tau)\dd\tau$$
into the weak formulation \eqref{weaktheta}, we have
\begin{equation}\nonumber
\begin{aligned}
&-\int_{\TTT^3}\int_0^T \int_{\eta}^{t-\eta}\del_s(\phi_{\epsilon}(\tau-s))(S_\epsilon S_Q\theta|S_\epsilon S_Q\theta|^{p-2})(\tau)\dd\tau S_Q\theta(s)\dd s\dd x\\
=&\int_{\eta}^{t-\eta} \langle S_{\epsilon}S_Q(v\theta), \nabla (S_\epsilon S_Q\theta |S_\epsilon S_Q\theta|^{p-2})\rangle \dd\tau.
\end{aligned}
\end{equation}
We rewrite the left term of this equality as
\begin{align*}
&\int_{\TTT^3}\int_{\eta}^{t-\eta}\int_0^T \del_{\tau}(\phi_{\epsilon}(\tau-s))S_Q\theta(s)\dd s(S_\epsilon S_Q\theta|S_\epsilon S_Q\theta|^{p-2})(\tau)\dd\tau \dd x\\
=&\int_{\TTT^3}\int_{\eta}^{t-\eta}(\del_{\tau}S_{\epsilon}S_Q\theta(s))(S_\epsilon S_Q\theta|S_\epsilon S_Q\theta|^{p-2})(\tau)\dd\tau \dd x\\
=&\frac{1}{p}\int_{\TTT^3}\int_{\eta}^{t-\eta} \del_{\tau}|S_\epsilon S_Q\theta(\tau) |^p\dd\tau\dd x,
\end{align*}
from which it follows that
\begin{equation}\label{energy-conserve}
\begin{aligned}
\frac{1}{p}\int_{\TTT^3}\int_{\eta}^{t-\eta} \del_{\tau}|S_\epsilon S_Q\theta(\tau) |^p\dd\tau\dd x=\int_{\eta}^{t-\eta} \langle S_{\epsilon}S_Q(v\theta), \nabla (S_\epsilon S_Q\theta |S_\epsilon S_Q\theta|^{p-2})\rangle \dd\tau
\end{aligned}
\end{equation}
Now we tackle with the term on the right-hand side of \eqref{energy-conserve}.  First of all, we decompose the term $S_Q(v\theta)$ by \eqref{SQ} as follows:
\begin{equation}\label{decomposition}
S_Q(v\theta)=r_Q(v,\theta)(x)-(v-S_Qv)(\theta-S_Q\theta)+S_QvS_Q\theta,
\end{equation}
where
\begin{equation}\nonumber
r_Q(v,\theta)(x)=2^{3Q}\int_{\mathbb{R}^3}\widetilde{h}(2^Qy)(v(x-y)-v(x))(\theta(x-y)-\theta(x))\dd y.
\end{equation}
Thanks to $\div S_Q v=0$, we have
\begin{align*}
&\int_{\eta}^{t-\eta} \langle S_{\epsilon}S_Q(v\theta), \nabla (S_\epsilon S_Q\theta |S_\epsilon S_Q\theta|^{p-2})\rangle \dd\tau\\
=&\int_{\eta}^{t-\eta} \langle S_{\epsilon}r_Q(v,\theta), \nabla (S_\epsilon S_Q\theta |S_\epsilon S_Q\theta|^{p-2})\rangle \dd\tau\\
&-\int_{\eta}^{t-\eta} \langle S_{\epsilon}((v-S_Q)(\theta-S_Q\theta)), \nabla (S_\epsilon S_Q\theta |S_\epsilon S_Q\theta|^{p-2})\rangle \dd\tau\\
&+\int_{\eta}^{t-\eta} \langle S_{\epsilon}(S_QvS_Q\theta)-S_QvS_{\epsilon} S_Q\theta, \nabla (S_\epsilon S_Q\theta |S_\epsilon S_Q\theta|^{p-2})\rangle \dd\tau\\
:=&\rm I+ II+III.
\end{align*}
For $u\in \mathcal{D}'(\TTT^3)$, we denote
\[d_{j,p}(u)=2^{\frac{1}{3}j}\|\Delta_j u\|_{L^p_x},\qquad d_p(u)=\{d_{j,p}(u)\}_{j\ge-1},\]
and
\begin{equation}\label{K}
K(j)=\left\{ \begin{alignedat}{-1}
&2^{\frac{2}{3}j},\quad j\le0,
 \\
 & 2^{-\frac{1}{3}j},\quad j\ge 1.
\end{alignedat}\right.
\end{equation}
By Berntein's inequality, one gets that
\begin{align}
\|v(x-y)-v(x)\|_{L^{\infty}_x}
\lesssim& \sum_{j\le Q}|y|2^j\|\Delta_j v\|_{L^{\infty}_x}+\sum_{j> Q}\|\Delta_j v\|_{L^{\infty}_x}\nonumber\\
=& 2^{\frac{2}{3}Q}|y|\sum_{j\le Q}2^{\frac{2}{3}(j-Q)}d_{j,\infty}(v)+2^{-\frac{1}{3}Q}\sum_{j> Q}2^{\frac{1}{3}(Q-j)}d_{j,\infty}(v)\nonumber\\
\le &C(2^{\frac{2}{3}Q}|y|+2^{-\frac{1}{3}Q})(K\ast d_{\infty}(v))(Q). \label{v-infty}
\end{align}

Similarly, we can infer from the above inequality that
\begin{equation}\label{t-p}
\begin{aligned}
\|\theta(x-y)-\theta(x)\|_{L^{p}_x}
\le C(2^{\frac{2}{3}Q}|y|+2^{-\frac{1}{3}Q})(K\ast d_{p}(\theta))(Q).
\end{aligned}
\end{equation}
Therefore, the term $\rm I$ can be bounded as
\begin{align}
|{\rm I}|\le& C\int_{\eta}^{t-\eta}S_{\epsilon}\big(2^{3Q}\int_{\mathbb{R}^3}|{\widetilde{h}}|(2^Qy)(2^{\frac{4}{3}Q}|y|^2+2^{-\frac{2}{3}Q})\dd y[K\ast (d_{\infty}(v))](Q)[K\ast d_{p}(\theta)](Q)\big)\nonumber\\
&\times\|\nabla (S_{\epsilon}S_Q\theta |S_{\epsilon}S_Q \theta|^{p-2})\|_{L^{\frac{p}{p-1}}_x}\dd \tau\nonumber\\
\le&C2^{-\frac{2}{3}Q}\int_{\eta}^{t-\eta}S_{\epsilon}[K\ast (d_{\infty}(v))](Q)[K\ast(d_{p}(\theta))](Q)\|\nabla (S_{\epsilon}S_Q\theta |S_{\epsilon}S_Q \theta|^{p-2})\|_{L^{\frac{p}{p-1}}_x}\dd \tau.\label{I}
\end{align}
Taking advantage of Bernstein's inequality, we have
\begin{align}
&\|\nabla (S_{\epsilon}S_Q\theta |S_{\epsilon}S_Q \theta|^{p-2})\|_{L^{\frac{p}{p-1}}_x}\le(p-1)\|\nabla S_{\epsilon}S_Q\theta\|_{L^p_x}\|S_{\epsilon}S_Q\theta\|^{p-2}_{L^p_x}\nonumber\\
\le&C2^{\frac{2}{3}Q}\sum_{j\le Q}2^{\frac{2}{3}(j-Q)}2^{\frac{1}{3}j}\|\Delta_jS_{\epsilon}\theta\|_{L^p}\|S_{\epsilon}\theta\|^{p-2}_{L^p_x}\nonumber\\
\le& C2^{\frac{2}{3}Q}S_{\epsilon}(K\ast d_p(\theta)(Q))\|S_{\epsilon}\theta\|^{p-2}_{L^p_x}.\label{nablasp-2}
\end{align}
Putting this inequality into \eqref{I} yields
\begin{align*}
|\rm I|\le& C\int_{\eta}^{t-\eta}S_{\epsilon}[K\ast (d_{\infty}(v))(Q)][K\ast (d_{p}(\theta))](Q)S_{\epsilon}[K\ast d_p(\theta)(Q)]\|S_{\epsilon}\theta\|^{p-2}_{L^p_x}\dd \tau.
\end{align*}
With the aid of \eqref{nablasp-2}, one can infer that
\begin{align*}
|{\rm II}|\le& \int_{\eta}^{t-\eta}\|v-S_Q v\|_{L^\infty_x}\|\theta-S_Q\theta\|_{L^p_x}\|\nabla (S_{\epsilon}S_Q\theta |S_{\epsilon}S_Q \theta|^{p-2})\|_{L^{\frac{p}{p-1}}_x}\dd\tau\\
\le& C\int_{\eta}^{t-\eta}2^{\frac{2}{3}Q}\|v-S_Q v\|_{L^\infty_x}\|\theta-S_Q\theta\|_{L^p_x}S_{\epsilon}(K\ast d_p(\theta)(Q))\|S_{\epsilon}\theta\|^{p-2}_{L^p_x}\dd\tau\\
\le& C\int_{\eta}^{t-\eta}\big(\sum_{j>Q}2^{\frac{1}{3}(Q-j)}2^{\frac{1}{3}j}\|\Delta_j v\|_{L^\infty_x}\big)\times\big(\sum_{j>Q}2^{\frac{1}{3}(Q-j)}2^{\frac{1}{3}j}\|\Delta_j \theta\|_{L^p_x}\big)\\
&\times S_{\epsilon}(K\ast d_p(\theta)(Q))\|S_{\epsilon}\theta\|^{p-2}_{L^p_x}\dd\tau\\
\le& C\int_{\eta}^{t-\eta}(K\ast d_{\infty}(v))(Q)(K\ast d_{p}(\theta))(Q)S_{\epsilon}(K\ast d_p(\theta)(Q))\|S_{\epsilon}\theta\|^{p-2}_{L^p_x}\dd\tau.
\end{align*}
Thanks to $S_Qv, S_Q\theta\in C^{\infty }(\TTT^3)$, $S_Qv\in L^\infty_TL^\infty_x$ and $S_Q\theta\in L^{p}_TL^p_x$, using the fact that $S_{\epsilon} f(\tau)\to f(\tau), a.e. \tau\in [\eta, t-\eta]$ for a locally integrable function $f(\tau)$, it is easy to check that $|{\rm III}|\to 0$ as $\epsilon\to 0$.  Therefore,
plugging the estimates of $\rm I$ and $\rm II$ into \eqref{energy-conserve}, and letting $\epsilon \to 0$, we have
\begin{equation}\label{SQtheta}
\begin{aligned}
\Big|\tfrac{1}{p}\int_{\TTT^3}\int_{\eta}^{t-\eta} \del_{\tau}| S_Q\theta(\tau) |^p\dd\tau\dd x\Big|
 \le
 &C\int_{\eta}^{t-\eta}[K\ast d_{\infty}(v)](Q)[K\ast d_{p}(\theta)]^2(Q)\|\theta\|^{p-1}_{L^p_x}\dd \tau.
\end{aligned}
\end{equation}
The fact that $v\in L^\infty([0,T]; B^{\frac{1}{3}}_{\infty,c(\NN)})$ and $\theta\in L^{p}([0,T]; B^{\frac{1}{3}}_{p,c(\NN)})$ yields
\begin{align*}
\lim_{Q\to\infty}d_{\infty}(v)(Q)=0,\,\,\text{and}\,\,\lim_{Q\to\infty}d_{p}(\theta)(Q)=0, \,\,{\rm a.e.}\,\, \,\,t\in [0,T].
\end{align*}
Therefore, for $K(j)$ defined in \eqref{K}, one deduces that
\begin{equation}\label{KQ}
\lim_{Q\to\infty}K\ast d_{p}(\theta)(Q)=\lim_{Q\to\infty}K\ast d_{\infty}(v)(Q)=0, \,\,{\rm a.e.}\,\, \, \,t\in[0,T],
\end{equation}
where we  have used  the fact \cite[(16)]{CCFS} that
\begin{align*}
0\le&\lim \sup_{Q\to\infty}K\ast d_{\infty}(v)(Q)+\lim \sup_{Q\to\infty}K\ast d_{p}(\theta)(Q)\\
&\le\lim \sup_{Q\to\infty}d_{\infty}(v)(Q)+\lim \sup_{Q\to\infty}d_{p}(\theta)(Q)=0.
\end{align*}
Taking $\eta\to 0$ in \eqref{SQtheta}, and then taking $Q\to\infty$,  by Lebesgue dominated convergence theorem and \eqref{KQ},then the right-hand side of \eqref{SQtheta} vanishes. Therefore,  we obtain that $\|\theta(t)\|_{L^p}=\|\theta_0\|_{L^p}$ for $2\le p<\infty$.

For the case of $p=\infty$, due to $\theta_0\in L^\infty(\TTT^3)\subset L^p(\TTT^3)$ and $\theta\in L^\infty([0,T]; B^{\frac{1}{3}}_{\infty, c(\NN)})\subset L^{p+1}([0,T]; B^{\frac{1}{3}}_{p, c(\NN)})$ for $2\le p<\infty$, we conclude that for $t\in [0,T)$ and for any $2\le p<\infty$,
\[\|\theta(t)\|_{L^p(\TTT^3)}=\|\theta_0\|_{L^p(\TTT^3)}\le \|\theta_0\|_{L^{\infty}(\TTT^3)}.\]
Using  $\displaystyle\lim_{p\to\infty}\|f\|_{L^p(\TTT^3)}=\|f\|_{L^{\infty}(\TTT^3)}$, we obtain that
\[\|\theta(t)\|_{L^{\infty}(\TTT^3)}=\|\theta_0\|_{L^{\infty}(\TTT^3)}.\]

For $p\in[1,2)$, we choose the test function $\varphi(x,s)$ as follows:
$$\varphi(x,s)=\int_{\eta}^{t-\eta}\phi_{\epsilon}(\tau-s)S_Q[\varphi^{p-1}_{\delta}(S_{\epsilon}S_Q\theta)\varphi'_{\delta}(S_{\epsilon}S_Q\theta)](x,\tau)\dd\tau,$$ where $\varphi_{\delta}(z)$ is a smooth approximation of $|z|$. More precisely,
 \begin{equation}\label{def.Phi}
\varphi_{\delta}(z)=\left\{ \aligned
    & -z-\tfrac{\delta}{2}, \quad z\le -{\delta},\\
     &\tfrac{1}{2\delta}z^2, \quad\quad \,\,\,\, |z|\le {\delta},\\
     &z-\tfrac{\delta}{2},\qquad\,\, z\ge {\delta}.
\endaligned
\right.
\quad\text{and}\quad
\varphi'_{\delta}(z)=\left\{ \aligned
    & -1, \quad z\le -{\delta},\\
     &\tfrac{z}{\delta}, \quad\,\,\,\,\,\,\,   |z|\le {\delta},\\
     &1,\qquad\,\, z\ge {\delta}.
\endaligned
\right.
\end{equation}
In the same way as leading to \eqref{energy-conserve}, we have
\begin{align}
&\frac{1}{p}\int_{\TTT^3}\int_{\eta}^{t-\eta}\partial_{\tau}\big[(\varphi_{\delta}(S_{\epsilon}S_Q\theta))^{p}\big]\dd\tau\dd x\nonumber\\
=&\int_{\eta}^{t-\eta} \langle S_{\epsilon}S_Q(v\theta), \nabla \big[(\varphi_{\delta}(S_{\epsilon}S_Q\theta))^{p-1}\varphi'_{\delta}(S_{\epsilon}S_Q\theta)\big]\rangle \dd\tau.\label{p-2}
\end{align}

By  \eqref{decomposition} and $\div S_Q v=0$, we decompose  the term on the right-hand side of \eqref{p-2} as
\begin{align*}
&\int_{\eta}^{t-\eta} \langle S_{\epsilon}S_Q(v\theta), \nabla \big[(\varphi_{\delta}(S_{\epsilon}S_Q\theta))^{p-1}\varphi'_{\delta}(S_{\epsilon}S_Q\theta)\big]\rangle \dd\tau\\
=&\int_{\eta}^{t-\eta} \langle S_{\epsilon}r_Q(v,\theta), \nabla \big[(\varphi_{\delta}(S_{\epsilon}S_Q\theta))^{p-1}\varphi'_{\delta}(S_{\epsilon}S_Q\theta)\big]\rangle \dd\tau\\
&-\int_{\eta}^{t-\eta} \langle S_{\epsilon}((v-S_Q)(\theta-S_Q\theta)), \nabla \big[(\varphi_{\delta}(S_{\epsilon}S_Q\theta))^{p-1}\varphi'_{\delta}(S_{\epsilon}S_Q\theta)\big]\rangle \dd\tau\\
&+\int_{\eta}^{t-\eta} \langle S_{\epsilon}(S_QvS_Q\theta)-S_QvS_{\epsilon} S_Q\theta, \nabla \big[(\varphi_{\delta}(S_{\epsilon}S_Q\theta))^{p-1}\varphi'_{\delta}(S_{\epsilon}S_Q\theta)\big]\rangle \dd\tau\\
:=&\rm IV+ V+VI.
\end{align*}
By making use of the arguments on estimating $\rm I$, $\rm II$ and $\rm III$ for the case $p\in[2,\infty)$, it suffices to estimate $\|\nabla \big((\varphi_{\delta}(S_{\epsilon}S_Q\theta))^{p-1}\varphi'_{\delta}(S_{\epsilon}S_Q\theta)\big)\|_{L^{\frac{p}{p-1}}_x}$.
Straightforward calculations entail that
\begin{align*}
\nabla \big[(\varphi_{\delta}(S_{\epsilon}S_Q\theta))^{p-1}\varphi'_{\delta}(S_{\epsilon}S_Q\theta)\big]
=&(p-1) (\varphi_{\delta}(S_{\epsilon}S_Q\theta))^{p-2}(\varphi'_{\delta})^2(S_{\epsilon}S_Q\theta)\nabla S_{\epsilon}S_Q\theta\\
&+(\varphi_{\delta}(S_{\epsilon}S_Q\theta))^{p-1}\varphi''_{\delta}(S_{\epsilon}S_Q\theta)\nabla S_{\epsilon}S_Q\theta.
\end{align*}
We deduce from the definition of $\varphi_{\delta}$ that
\begin{align*}
&\|(\varphi_{\delta}(S_{\epsilon}S_Q\theta))^{p-2}(\varphi'_{\delta})^2(S_{\epsilon}S_Q\theta))\nabla S_{\epsilon}S_Q\theta\|_{L^{\frac{p}{p-1}}_x}\\
\le &C\delta^{-p}\|\chi_{|S_{\epsilon}S_Q\theta|\le\delta}|S_{\epsilon}S_Q\theta|^{2(p-1)}\nabla S_{\epsilon}S_Q\theta\|_{L^{\frac{p}{p-1}}_x}\\
&+C\|\chi_{|S_{\epsilon}S_Q\theta|>\delta}|S_{\epsilon}S_Q\theta|^{p-2}\nabla S_{\epsilon}S_Q\theta\|_{L^{\frac{p}{p-1}}_x}\\
\le&C\delta^{p-2}\|\nabla S_{\epsilon}S_Q\theta\|_{L^{\frac{p}{p-1}}_x}.
\end{align*}
Similarly, one has
\begin{align*}
\|(\varphi_{\delta}(S_{\epsilon}S_Q\theta))^{p-1}(\varphi''_{\delta})(S_{\epsilon}S_Q\theta)\nabla S_{\epsilon}S_Q\theta\|_{L^{\frac{p}{p-1}}_x}
\le&C\delta^{p-2}\|\nabla S_{\epsilon}S_Q\theta(\tau)\|_{L^{\frac{p}{p-1}}_x}.
\end{align*}
Collecting the above two estimates together shows that
\begin{equation}\label{testp-2}
\|\nabla \big((\varphi_{\delta}(S_{\epsilon}S_Q\theta))^{p-1}\varphi'_{\delta}(S_{\epsilon}S_Q\theta)\big)\|_{L^{\frac{p}{p-1}}_x}\le C\delta^{p-2}\|\nabla S_{\epsilon}S_Q\theta(\tau)\|_{L^{\frac{p}{p-1}}_x}.
\end{equation}
Using the estimates in the process of bounding $\rm I, II $ and $\rm III$ and \eqref{testp-2}, we deduce from \eqref{p-2}  that
\begin{align*}
&\Big|\frac{1}{p}\int_{\TTT^3}\int_{\eta}^{t-\eta}\partial_{\tau}\big[(\varphi_{\delta}(S_{\epsilon}S_Q\theta))^{p}\big]\dd\tau\dd x\Big|\\
\le&C\delta^{p-2}\int_{\eta}^{t-\eta}S_{\epsilon}[K\ast d_{\infty}(v)](Q)[K\ast d_{p}(\theta)](Q) S_{\epsilon}[K\ast d_{\frac{p}{p-1}}(\theta)](Q)\dd \tau\\
&+C\delta^{p-2}\int_{\eta}^{t-\eta}(K\ast d_{\infty}(v))(Q)(K\ast d_{p}(\theta))(Q)S_{\epsilon}[K\ast d_{\frac{p}{p-1}}(\theta)](Q)\dd\tau\\
&+\Big|\int_{\eta}^{t-\eta} \langle S_{\epsilon}(S_QvS_Q\theta)-S_QvS_{\epsilon} S_Q\theta, \nabla [(\varphi_{\delta}(S_{\epsilon}S_Q\theta))^{p-1}\varphi'_{\delta}(S_{\epsilon}S_Q\theta)]\rangle \dd\tau\Big|.
\end{align*}
Taking $\epsilon\to 0$, and then letting $\eta\to 0$ in  this inequality, we obtain that
\begin{align}
&\Big|\tfrac{1}{p}\|\varphi_{\delta}(S_Q\theta)(t)\|^p_{L^p}-\tfrac{1}{p}\|\varphi_{\delta}(S_Q\theta_0)\|^p_{L^p}\Big|\nonumber\\
\le& C\delta^{p-2}\int_0^t(K\ast d_{\infty}(v))(Q)(K\ast d_{p}(\theta))(Q)(K\ast d_{\frac{p}{p-1}}(\theta))(Q)\dd\tau.\label{SQp-2}
\end{align}
Since  $v\in L^{\infty}([0,T];B^{\frac{1}{3}}_{\infty, c(\NN)})$ and $\theta\in L^{2}([0,T];B^{\frac{1}{3}}_{p, c(\NN)}\cap B^{\frac{1}{3}}_{{\frac{p}{p-1}}, c(\NN)})$, we have as in proving~\eqref{KQ} that
 \begin{equation}\nonumber
\lim_{Q\to\infty}K\ast d_{{\frac{p}{p-1}}}(\theta)(Q)=\lim_{Q\to\infty}K\ast d_{p}(\theta)(Q)=\lim_{Q\to\infty}K\ast d_{\infty}(v)(Q)=0, \,{\rm a.e.} \,\,t\in[0,T].
\end{equation}
Therefore, letting $Q\to\infty$ in \eqref{SQp-2} then gives  $\|\varphi_{\delta}(\theta(t))\|_{L^p}=\|\varphi_{\delta}(\theta_0)\|_{L^p}$. Finally, taking $\delta\to 0$, we conclude that $\|\theta(t)\|_{L^p}=\|\theta_0\|_{L^p}$ for $1\le p<2$. In conclusion, we complete the proof of Theorem \ref{main}.
\section{Proof of Proposition \ref{t:main}}\label{proof of negative}
In this section, we reduce  Proposition~\ref{t:main} to an  iteration proposition---Proposition~\ref{p:main-prop} below.  In Section 4.1,  we show that Proposition \ref{t:main} will be achieved by  applying Proposition~\ref{p:main-prop} inductively. The proof of Proposition \ref{p:main-prop} occupies the rest of this section, mainly consisting of mollification, gluing and constructing the perturbation procedure.
 \subsection{Iteration proposition}
 For the sake of proving  Proposition \ref{t:main} in the framework of convex integration,  we consider  the  following relaxation  of  Boussinesq equations.
 \begin{equation}
\left\{ \begin{alignedat}{-1}
&\del_t v_q-\div (v_q\otimes v_q)  +\nabla p_q   =e_3\theta_q + \div \RR_q,
 \\
 &\del_t \theta_q+v_q\cdot\nabla \theta_q     = \div \MM_q,
 \\
  &\nabla \cdot v_q = 0,
  \\ &(v_q, \theta_q) |_{t=0}=(\vin, \tin):=(v^{(1)}(0, x), \theta^{(1)}(0, x)),
\end{alignedat}\right.  \label{e:subsol-B}
\end{equation}
where
\begin{equation}\label{p_q}
\int_{\TTT^3}p_q\dd x=0,
\end{equation}
 the \emph{Reynolds stress} $\RR_q$  is a symmetric trace-free $3\times3$  matrix and the \emph{temperature stress} $\MM_q$ is a vector field.

Before giving the iteration proposition,  we introduce all parameters needed  in the inductive procedure.   For all $q\ge 1$ and given  $0<\beta<1/3$, we define $$b_0=\min\big\{1+\tfrac{1-3\beta}{2\beta},2\big\}.$$
For the positive parameters  $a\gg 1$ and  $b\in(1, b_0)$, we define
\begin{align}\label{lambdaq}
    \lambda_q \coloneq  \left\lceil a^{b^q}\right\rceil,\qquad   \delta_q \coloneq \lambda_2^{3\beta}\lambda_q^{-2\beta}.
\end{align}
here $\lceil x\rceil$ denotes the ceiling function.
Given $b$ and $\beta$, there exists a $\alpha_0>0$ depending on $b, \beta$ such that for any $0<\alpha<\alpha_0(b, \beta)$,
\begin{equation}\label{alpha0}
\frac{\delta^{1/2}_{q+1}\delta^{1/2}_q\lambda_q}{\lambda_{q+1}}\lesssim\frac{\delta_{q+2}}{\lambda^{8\alpha}_{q+1}}.
\end{equation}
Let $0<\widetilde{T}\le\tfrac{1}{4}$ and $M$ be a universal constant defined in \eqref{M}. In the following, we  require that
\begin{align}
   \alpha < \min\big\{\alpha_0,\tfrac{\beta (b-1)}{3}\big\},\quad  a > \max\Big\{50^{\beta/\alpha },  M^{\frac{2}{\alpha}}, \tfrac{3}{\widetilde{T}}\Big\}. \label{e:params0}
\end{align}
{These inequalities imply that  $\lambda_q^{3\alpha} \le \frac{\delta_q^{3/2}}{\delta_{q+1}^{3/2}}  \le  \frac{\lambda_{q+1}}{\lambda_q}$, as in \cite{BDSV, KMY}.}
And we define $\ell_q$ and $\tau_q$ by
\begin{align}
    \ell_q \coloneq \frac{\delta_{q+1}^{1/2}}{\delta_q^{1/2} \lambda_q^{1+3\alpha/2}}  \in \Big(\frac12\lambda_q^{-1-\frac{3\alpha}2-(b-1)\beta},\lambda_q^{-1-\frac{3\alpha}2}\Big) \subset \Big(\lambda^{-\frac{13}{10}}_q,\lambda_q^{-1}\Big)\label{e:ell},\,\,\,
     \tau_q \coloneq \frac{\ell_q^{2\alpha}}{\delta_q^{1/2} \lambda_q}.
\end{align}
This definition entails  that $M\le \ell^{-\frac{\alpha}{2}}_q$ and $\ell^{-1}_q<\lambda_{q+1}$. Finally, we define  $N_0(b,\beta)\in \ZZ_{+}$ such that
\begin{equation}\label{lambdaN}
 \frac{1}{\lambda_{q+1}^{N_0-\alpha} \ell_q^{N_0+\alpha}} \leq \frac{1}{\lambda_{q+1}^{1-\alpha}}, \quad \forall q\ge 1.
\end{equation}

Now we are in position to state the   iteration proposition.
\begin{prop}
\label{p:main-prop}
Let $b\in(1,b_0)$ and $\alpha$ satisfy \eqref{e:params0}.
There exists  $a_0$ such that for any $\alpha<\alpha_0(\beta,b)$ and $a>a_0(\beta, b,\alpha)$, the following holds.
If $(v_q, \theta_q, p_q,\RR_q, \MM_q)$ obeys equations \eqref{e:subsol-B} and \eqref{p_q} with
\begin{align}
     &\|(v_q, \theta_q)\|_{0} \le\sum_{i=1}^q\delta^{1/2}_i ,
    \label{e:vq-C0}
    \\
    &\|(v_q, \theta_q)\|_{1} \le \delta_{q}^{1/2} \lambda_q ,
    \label{e:vq-C1}
    \\
    &\|(\RR_q, \MM_q)\|_{0} \le \delta_{q+1}\lambda_q^{-3\alpha} ,
    \label{e:RR_q-C0}
\\
&(v_q, \theta_q)=(v^{(1)}, \theta^{(1)}) \,\,\text{\rm on}\,\, [0, 2\widetilde{T}+\tau_q],\label{e:initial1}\\
 &(v_q, \theta_q)= (v^{(2)}, \theta^{(2)})\,\,\text{\rm on}\,\,[3\widetilde{T}-2\tau_q, T_2],\label{e:initial}
\end{align}
where  $T_1$, $T_2$ $\widetilde{T}$, $(v^{(1)}, \theta^{(1)})$ and $ (v^{(2)}, \theta^{(2)})$ are consistent with those in Proposition \ref{t:main}, then there exist smooth functions $(v_{q+1}, \theta_{q+1}, p_{q+1}, \RR_{q+1}, \MM_{q+1})$ satisfying \eqref{e:subsol-B}, \eqref{p_q},
\eqref{e:vq-C0}--\eqref{e:initial}
with $q$ replaced by $q+1$ and
\begin{align}
        \big\|(v_{q+1}, \theta_{q+1}) - (v_q, \theta_q)\big\|_{0} +\frac1{\lambda_{q+1}} \big\|(v_{q+1}, \theta_{q+1})-(v_q, \theta_q)\big\|_1 &\le   \delta_{q+1}^{1/2}.
        \label{e:velocity-diff}
\end{align}
\end{prop}
\noindent \textbf{Proposition \ref{p:main-prop} implies Proposition \ref{t:main}.}\quad Without loss of generality, we assume that $T_1=T_2=1$. To implement the iterative procedure, we firstly  construct $(v_1, b_1, p_1, \RR_1, \MM_1)$ as follows. {Let $\chi(t)$ be a smooth cut-off function  such that $\text{supp}\chi(t)=[-1, 2]$, $\chi|_{[0, 2\widetilde{T}+\tau_1]}=1$} and $\chi|_{[3\widetilde{T}-2\tau_1, 2]}=0$.  We construct  $(v_1, \theta_1)$ by gluing $(v^{(1)}, \theta^{(1)})$ and $(v^{(2)}, \theta^{(2)})$ as follows:
\begin{align*}
(v_1, \theta_1)=(\chi v^{(1)}+(1-\chi)v^{(2)}, \chi \theta^{(1)}+(1-\chi)\theta^{(2)}) \,\,\text{on} \,\,[0,1].
\end{align*}
Since $\nabla \cdot v^{(1)}=\nabla \cdot v^{(2)}=0$ and $\chi$ only depends on time variable, $v_1$ is divergence-free. For sufficient large $a$,  it is easy to verify that $(v_1, \theta_1)$  satisfies
\begin{align}
&\|(v_1,~\theta_1)\|_0\le\delta^{1/2}_2\lambda^{-4\alpha}_1\le \delta^{\frac{1}{2}}_1,\label{v1b1}\\
&\|(v_1,~\theta_1)\|_{1}\le\lambda_1.\nonumber
\end{align}
Therefore, $(v_1, \theta_1)$ satisfies \eqref{e:vq-C0} and \eqref{e:vq-C1} for $q=1$. By the definition of $(v_1, \theta_1)$, we have
\begin{align}
    \RR_1 \coloneq&
        \del_t \chi \mathcal R(v^{(1)}-v^{(2)} ) - \chi(1-\chi)(v^{(1)}-v^{(2)}  )\ootimes (v^{(1)}-v^{(2)}  ), \label{RR_1}\\
        \MM_1 \coloneq&
        \del_t \chi \mathcal R_{\vex}(\theta^{(1)}-\theta^{(2)}  ) - \chi(1-\chi)(v^{(1)}-v^{(2)}   )(\theta^{(1)}-\theta^{(2)} ), \label{MM_1}\\
    p_1 \coloneq & \chi \big(p^{(1)}-\int_{\TTT^3}p^{(1)}\dd x\big)+(1-\chi)\big(p^{(2)}-\int_{\TTT^3}p^{(2)}\dd x\big)\nonumber\\
     &- \chi(1-\chi)\big( |v^{(1)}-v^{(2)}  |^2  - \int_{\mathbb T^3} |v^{(1)}-v^{(2)} |^2 \dd x\big)\nonumber,
\end{align}
here and below, $v\ootimes u:=v\otimes u-\frac{1}{3}(v\cdot u)\rm{Id}_{3\times3}$. By the estimate \eqref{v1b1},
 we can obtain from \eqref{RR_1} and \eqref{MM_1}  that
\begin{align*}
\|(\RR_1, \MM_1)\|_{0}\le C\delta^{1/2}_2\lambda^{-4\alpha}_1+C\delta_2\lambda^{-8\alpha}_1\le \delta_2\lambda^{-3\alpha}_1.
\end{align*}
Noting the support of $\chi(t)$, we deduce that $(\RR_1, \MM_1)\equiv 0$ on $[0, 2\widetilde{T}+\tau_1]\cup [3\widetilde{T}-2\tau_1, 1]$, which implies that
\[(v_1, \theta_1)=(v^{(1)}, \theta^{(1)}) \,\,\text{on}\,\, [0, 2\widetilde{T}+\tau_1],\quad (v_1, \theta_1)= (v^{(2)}, \theta^{(2)})\,\,\text{on}\,\,[3\widetilde{T}-2\tau_1, 1].\]
Using Proposition \ref{p:main-prop}, we can get a sequence of solutions $\{(v_q, \theta_q, p_q,\RR_q,\MM_q)\}$ to the system \eqref{e:subsol-B} satisfying \eqref{p_q} and \eqref{e:vq-C1}--\eqref{e:initial}. First of all, $(v_q, \theta_q)$ converges uniformly  to some continuous functions $(v, \theta)$ by \eqref{e:velocity-diff}. From \eqref{p_q}, \eqref{e:vq-C0} and \eqref{e:RR_q-C0}, one deduces that $p_q$ converges to $p$ in $L^m$ for $m<\infty$. Moreover, since $\|(\RR_q, \MM_q)\|_{0}\rightarrow 0$ as $q\to\infty$, the triple $(v, \theta, p)$ obeys the Boussinesq system \eqref{e:B} in a weak sense.  By a standard interpolation argument with $C^0$ and $C^1$ estimates in \eqref{e:velocity-diff}, we conclude that $(v_q, \theta_q)$ is a Cauchy sequence in $C^0_tC^{\beta'}_x$ for all $\beta'<\beta$, which implies that  $(v, \theta)\in C^0_tC^{\beta'}_x$. On the other hand, notice  that $(v,\theta)$ obeys the equations \eqref{e:B}, one obtains that $(v, \theta)\in C^{\beta''}_tC^{0}_x$ for all $\beta''<\beta'$ as in \cite{BDSV}. Therefore, there exists a weak solution $(v, \theta)\in C^{\beta''}(\TTT^3\times[0, 1])$ solving the equations \eqref{e:B}. Furthermore,
\begin{align*}
 (v, b)\equiv (v^{(1)}, b^{(1)}) \,\,\,\text {on}\,\,\,[0, 2\widetilde{T}] \,\,\,\text{and} \,\,\,(v, b)\equiv (v^{(2)}, b^{(2)}) \,\,\, \text {on}\,\,\, [3\widetilde{T}, 1].
 \end{align*}
Hence, we complete the proof of Theorem \ref{t:main} by Proposition \ref{p:main-prop}.

The rest of this section is devoted to the proof of Proposition \ref{p:main-prop}. More precisely,
\begin{enumerate}
  \item [$\bullet$]In Section 4.2, we complete the two stages: \\
  mollification: $(v_q, \theta_q, p_q, \RR_q,\MM_q)\mapsto (v_{\ell_q}, \theta_{\ell_q}, p_{\ell_q}, \RR_{\ell_q}, \MM_{\ell_q})$,\\
  gluing:  $(v_{\ell_q}, \theta_{\ell_q}, p_{\ell_q}, \RR_{\ell_q}, \MM_{\ell_q})\mapsto (\vv_q, \tb_q, \pp_q, \RRR_q, \MMM_q)$;
  \item [$\bullet$]In Section 4.3, we construct the perturbation to produce a new solution at $q+1$ level:\\
  perturbation: $ (\vv_q, \tb_q, \pp_q, \RRR_q, \MMM_q)\mapsto (v_{q+1}, \theta_{q+1}, p_{q+1}, \RR_{q+1},\MM_{q+1})$;\\
  And we give the H\"{o}lder estimates for the  perturbation parts.
  \item [$\bullet$]Section 4.4 contains the H\"{o}lder estimates for the new Reynolds stress  and temperature stress  $(\RR_{q+1}, \MM_{q+1})$ and we further show Proposition \ref{p:main-prop} in Section~4.5.
\end{enumerate}
\subsection{Mollification and gluing procedure}
\subsubsection{Mollification}
This step plays an important role  in a convex integration scheme,
which can be used  to  solve the loss of derivative problem. For $\ell_q$ in \eqref{e:ell}, we define functions $(v_{\ell_q}, \theta_{\ell_q}, p_{\ell_q}, \RR_{\ell_q}, \MM_{\ell_q})$ by the spatial mollifier defined in \eqref{e:defn-mollifier-x} as follows.
\begin{align*}
   & v_{\ell_q} \coloneq v_q * \psi_{\ell_q}, \quad\theta_{\ell_q} \coloneq \theta_q * \psi_{\ell_q},  \quad  p_{\ell_q}\coloneq p_q *\psi_{\ell_q} -|v_q|^2 + |v_{\ell_q}|^2,\\
 &   \RR_{\ell_q} \coloneq \RR_q * \psi_{\ell_q}  - (v_q \ootimes v_q) * \psi_{\ell_q}  + v_{\ell_q} \ootimes v_{\ell_q} , \\
 &\MM_{\ell_q} \coloneq \MM_q * \psi_{\ell_q}  - (v_q \theta_q) * \psi_{\ell_q}  + v_{\ell_q}\theta_{\ell_q}.
\end{align*}
They obey the following equations
\begin{equation}
\left\{ \begin{alignedat}{-1}
&\del_t v_{\ell_q} +\div (v_{\ell_q}\otimes v_{\ell_q})  +\nabla p_{\ell_q}   =  \div \RR_{\ell_q} +\theta_{\ell_q}e_3,
\\
&\del_t\theta_{\ell_q}+v_{\ell_q}\cdot\nabla\theta_{\ell_q}=\div\MM_{\ell_q},\\
 & \nabla \cdot v_{\ell_q} = 0,
  \\
  &( v_{\ell_q}, \theta_{\ell_q})\big|_{t=0}= (\vin*\psi_{\ell_{q}}, \tin*\psi_{\ell_{q}}).
\end{alignedat}\right.  \label{e:mollified-euler}
\end{equation}
By the standard method as in \cite[Propistion 2.2]{BDSV} and $\lambda^{-3\alpha}_q\le \ell^{\frac{9}{4}\alpha}_q$,  we easily conclude the following bounds.
\begin{prop}[ Estimates for mollified functions]\label{p:estimates-for-mollified}
\begin{align}
&\|(v_{\ell_q}-v_{q}, \theta_{\ell_q}-\theta_{q})\|_{0} \lesssim \delta_{q+1}^{1 / 2}  \lambda_{q}^{-\alpha}, \label{e:v_ell-vq}
\\
&\|(v_{\ell_q}, \theta_{\ell_q})\|_{N+1} \lesssim \delta_{q}^{1 / 2} \lambda_{q} \ell_q^{-N}, && \forall N \geq 0, \label{e:v_ell-CN+1}
\\
&\|(\RR_{\ell_q}, \MM_{\ell_q})\|_{N+\alpha} \lesssim  \delta_{q+1} \ell_q^{-N+\frac{5}{4}\alpha}, && \forall N \geq 0 .\label{e:R_ell}
\end{align}
\end{prop}
\subsubsection{Classical exact  flows}
We define $
t_0=2\widetilde{T}, t_i\coloneq t_0+ i\tau_q, i\in\NN$
and
\begin{align*}
  i_{\max}=\sup\{i\ge 1, \tau_q\le3\widetilde{T}-t_i<2\tau_q\}+1.
\end{align*}
Let
\begin{equation}\nonumber
 (\vex_0, {\bm{\theta}_0}, \pex_0)=(\vex_{i_{\max}}, {\bm\theta_{i_{\max}}}, \pex_{i_{\max}}):=(v_q, \theta_q, p_q).
 \end{equation}
Note that
$$(v_q, \theta_q)=(v^{(1)}, \theta^{(1)})\,\,{\rm on}\,\,[0, t_0+\tau_q], \,\,\,(v_q, \theta_q)=(v^{(2)}, \theta^{(2)})\,\,{\rm on}\,\,[3\widetilde{T}-2\tau_q, 1]. $$
This clearly entails  that
\begin{equation}\label{def-v0}
(\vex_0, {\bm{\theta}_0})=(v^{(1)}, \theta^{(1)})\,\,{\rm on}\,\,[0, t_0+\tau_q], \,\,\,(\vex_{i_{\max}}, {\bm\theta_{i_{\max}}})=(v^{(2)}, \theta^{(2)})\,\,{\rm on}\,\,[3\widetilde{T}-2\tau_q, 1],
\end{equation}
which imply that $(\vex_0, {\bm{\theta}_0})$ and $(\vex_{i_{\max}}, {\bm{\theta}}_{i_{\max}})$ are the exact solutions to the Boussinesq equations~\eqref{e:B} on $[0, t_0+\tau_q]$ and $[3\widetilde{T}-2\tau_q, 1]$, respectively. For $1\le i\le i_{\max}-1$, let $(\vex_i, \bm{\theta}_i, \pex_i)$ be the unique exact solution to the following Boussinesq equations
 \begin{equation}
 \left\{ \begin{alignedat}{-1}
&\del_t \vex_i +(\vex_i\cdot\nabla) \vex_i  +\nabla \pex_i   =  \bm\theta_i e_3,
\\
&\del_t  \bm\theta_i +(\vex_i\cdot\nabla)  \bm\theta_i  =0,
\\
&  \nabla \cdot \vex_i  = 0,
  \\
& (\vex_i, \bm\theta_i)\big|_{t=t_i} = (v_{\ell_q}(\cdot,t_i), \theta_{\ell_q}(\cdot,t_i)).
\end{alignedat}\right.  \label{e:exact-B}
\end{equation}
defined over their own maximal interval of existence.
{The fact that the lifespan of smooth solutions to the above system is greater than $[t_{i}-\tau_q, t_{i}+\tau_q]$ can be guaranteed by the following proposition.}
\begin{prop}[Estimates for exact solutions to Boussinesq equations]
\label{p:exact-B}
Let $ N\ge1$ and $\alpha\in(0,1)$. Assume that initial data $(V_0, \Theta_0)\in C^{ \infty}$, then  there exists a universal positive constant $c$ such that the system \eqref{e:B} possesses a unique smooth solution $(V,\Theta, P)$ to \eqref{e:B} on $[-T, T]$, where $T= \frac c{\|V_0\|_{{1+\alpha}}+1}$, and satisfies the following derivative estimates,
\begin{equation}\label{est-N+alpha}
 \|V\|_{C^0([- T, T];C^{N+\alpha})} +  \|\Theta\|_{C^0([- T,T];C^{N+\alpha})} \lesssim  \| V_0 \|_{{N+\alpha}}+\|\Theta_0 \|_{{N+\alpha}} .
 \end{equation}
\end{prop}
\begin{proof}We focus on the proof of the a priori estimate for smooth solutions. The existence of a unique solution is standard (see e.g. \cite{MB})
We denote  $\matd{V}$ by material derivative. For any multi-index $\sigma$ with $|\sigma|=N$, we have
\[\matd{V}\partial^\sigma V+[\partial^\sigma, V\cdot\nabla]V+\nabla\partial^\sigma P=\partial^\sigma(\Theta e_3),\]
and
\[\matd{V}\partial^\sigma \Theta +[\partial^\sigma, V\cdot\nabla]\Theta=0.\]
Therefore, by commutator estimates we have
\begin{equation}\nonumber
\|\matd{V}\partial^\sigma V\|_{\alpha}\lesssim \|V\|_{1+\alpha}\|V\|_{N+\alpha}+\|\Theta\|_{\alpha+N}.
\end{equation}
This estimate shows that for $|t|\le T= \frac c{\|V_0\|_{{1+\alpha}}+1}$,
\begin{equation}\label{VN}
\| V(t)\|_{N+\alpha}\lesssim \|V_0\|_{N+\alpha}+\int_0^{|t|}\big(\|V(\tau)\|_{1+\alpha}\|V(\tau)\|_{N+\alpha}+\|\Theta(\tau)\|_{N+\alpha}\big)\dd \tau.
\end{equation}
Similarly, we can bound $\Theta$ by
\begin{equation}\label{TN}
\| \Theta(t)\|_{N+\alpha}\lesssim \| \Theta_0\|_{N+\alpha}+\int_0^{|t|}\big(\|V(\tau)\|_{1+\alpha}\|\Theta(\tau)\|_{N+\alpha}+\|V(\tau)\|_{N+\alpha}\|\Theta(\tau)\|_{1+\alpha}\big)\dd \tau.
\end{equation}
For $N=1$, by Gr\"{o}nwall's inequality and the continuity argument, we obtain that for $|t|\le T$, we have
\[\|V\|_{1+\alpha}+\|\Theta\|_{1+\alpha}\lesssim \|V_0\|_{1+\alpha}+\|\Theta_0\|_{1+\alpha}.\]
Combining this estimate with \eqref{VN} and \eqref{TN} leads to \eqref{est-N+alpha} for all $N\ge 1$.
\end{proof}
This proposition together with \eqref{e:v_ell-CN+1} immediately shows that  $(\vex_i, \bm\theta_i, \pex_i)$ is well-defined on $[t_{i}-\tau_q, t_{i}+\tau_q]$ provided that the parameter $a$ is chosen sufficiently large. Furthermore,  we can obtain the following stability estimates.
\begin{prop}[Stability]
\label{p:stability}For $i\ge 0$, $|t-t_i|\le \tau_q$, and $0\le N\le N_0+1$, we have
\begin{align}
    \|(\vex_i -v_{\ell_q}, \bm\theta_i -\theta_{\ell_q}) \|_{{N+\alpha}} & \lesssim \tau_q \delta_{q+1} \ell_q^{-N-1+\frac{5}{4}\alpha}\label{e:stability-v},\\
    \|\nabla \pex_i -\nabla p_{\ell_q} \|_{{N+\alpha}} &\lesssim\delta_{q+1} \ell_q^{-N-1+\frac{5}{4}\alpha},\label{e:stability-p}\\
    \|(\matd{v_{\ell_q}}(\vex_i - v_{\ell_q}),  \matd{v_{\ell_q}}(\bm\theta_i - \theta_{\ell_q}))\|_{{N+\alpha}} & \lesssim  \delta_{q+1} \ell_q^{-N-1+\frac{5}{4}\alpha},\label{e:stability-matd}\\
  \|(\mathcal{R}(\vex_i -v_{\ell_q}), \mathcal{R}_{\vex}(\bm\theta_i -\theta_{\ell_q})) \|_{{\alpha}} &\lesssim \tau_q \delta_{q+1} \ell_q^{\frac{5}{4}\alpha}\label{e:stability-v-1}.
\end{align}
\end{prop}
\begin{proof}
From \eqref{e:mollified-euler} and \eqref{e:exact-B}, we have
 \begin{align}
& \matd{v_{\ell_q}} (\vex_i - v_{\ell_q}) + (\vex_i-v_{\ell_q})\cdot\nabla \vex_i + \nabla (\pex_i - p_{\ell_q}) = (\bm\theta_i-\theta_{\ell_q})- \div \RR_{\ell_q},\label{diff-1}\\
&\matd{v_{\ell_q}} (\bm\theta_i - \theta_{\ell_q}) +(\vex_i-v_{\ell_q})\cdot\nabla \bm\theta_i  =-\div\MM_{\ell_q}\label{diff-2}.
 \end{align}
 Noting that $\div \vex_i=\div v_{\ell_q}=0$, we obtain from \eqref{diff-1} that for $N\ge0$,
\begin{align}
 \|\nabla (\pex_i - p_{\ell_q})\|_{N+\alpha}\lesssim& \|\vex_i - v_{\ell_q}\|_{N+\alpha}\|(\vex_i, v_{\ell_q})\|_{1+\alpha}+\|\vex_i - v_{\ell_q}\|_{\alpha}\|(\vex_i, v_{\ell_q})\|_{N+\alpha}\notag\\
&+ \| \RR_{\ell_q}\|_{N+1+\alpha}+\|\bm\theta_i - \theta_{\ell_q}\|_{{N+\alpha}}.\label{pex-p}
\end{align}
First of all, we show  the estimate \eqref{e:stability-v}. For $1\le i\le i_{\max}-1$, by Lemma \ref{trans}, \eqref{e:v_ell-CN+1} and \eqref{e:R_ell}, due to $(\bm\theta_i - \theta_{\ell_q})(t_i)=0$, one obtains  that for $t\in [t_{i-1}, t_{i+1}]$,
\begin{align*}
\|(\bm\theta_i - \theta_{\ell_q})(t)\|_{\alpha}
\le &\exp\big(\int_{t_{i-1}}^{t_{i+1}}\|\nabla v_{\ell_q}\|_{0}\dd s\big)\big(\int_{t_{i-1}}^{t}(\|(\vex_i-v_{\ell_q})\cdot\nabla \bm\theta_i \|_{\alpha}+\|\div\MM_{\ell_q}\|_{\alpha})\dd s\big)\\
\lesssim&\int_{t_{i-1}}^t\|\vex_i-v_{\ell_q}\|_{\alpha}\|\bm\theta_i\|_{1+\alpha}\dd s+\tau_q\delta_{q+1}\ell^{-1+\frac{5}{4}\alpha}_q.
\end{align*}
Using \eqref{e:v_ell-CN+1}, \eqref{e:R_ell}, \eqref{pex-p} and Lemma \ref{trans} leads to
\begin{align*}
\|(\vex_i - v_{\ell_q})(t)\|_{\alpha}
\lesssim &\int_{t_{i-1}}^{t}\big(\|\vex_i - v_{\ell_q}\|_{\alpha}\|(\vex_i, v_{\ell_q})\|_{1+\alpha}+\|\bm\theta_i - \theta_{\ell_q}\|_{{\alpha}}+\| \div\RR_{\ell_q}\|_{\alpha}\big)\dd s\\
\lesssim&\int_{t_{i-1}}^t\big(\|\vex_i - v_{\ell_q}\|_{\alpha}\|(\vex_i, v_{\ell_q})\|_{1+\alpha}+\|\bm\theta_i - \theta_{\ell_q}\|_{{\alpha}}\big)\dd s+\tau_q\delta_{q+1}\ell^{-1+\frac{5}{4}\alpha}_q.
\end{align*}
Collecting the above two inequalities together, by making use of Gr\"{o}nwall's inequality and~\eqref{e:v_ell-CN+1}, we obtain
\begin{align}
\|(\vex_i - v_{\ell_q})(t)\|_{\alpha}+\|(\bm\theta_i - \theta_{\ell_q})(t)\|_{\alpha}
\lesssim&\tau_q\delta_{q+1}\ell^{-1+\frac{5\alpha}{4}}_q\exp\Big(\int_{t_{i-1}}^t\|(\vex_i, v_{\ell_q}, \bm\theta_i )\|_{1+\alpha}\dd s\Big)\nonumber\\
\lesssim&\tau_q\delta_{q+1}\ell^{-1+\frac{5\alpha}{4}}_q\exp(\tau_q\delta^{1/2}_q\lambda_q\ell^{-\alpha}_q)\nonumber\\
\lesssim&\tau_q\delta_{q+1}\ell^{-1+\frac{5}{4}\alpha}_q.\label{vex-vq-alpha}
\end{align}
For $N\ge 1$ and $|\sigma|=N$, one deduces from \eqref{diff-1} and \eqref{diff-2} that
\begin{align}
& \matd{v_{\ell_q}} \partial^\sigma(\vex_i - v_{\ell_q}) +[\partial^\sigma, v_{\ell_q} \cdot\nabla](\vex_i - v_{\ell_q}) + \partial^\sigma\big((\vex_i-v_{\ell_q})\cdot\nabla \vex_i\big)\nonumber\\
&\quad\quad \quad\quad\quad\quad\quad\quad\quad\quad\quad\quad\quad\quad+ \partial^\sigma\nabla (\pex_i - p_{\ell_q}) = \partial^\sigma(\bm\theta_i-\theta_{\ell_q})- \partial^\sigma\div \RR_{\ell_q},\label{diff-1N}\\
&\matd{v_{\ell_q}} \partial^\sigma(\bm\theta_i - \theta_{\ell_q}) +[\partial^\sigma, v_{\ell_q} \cdot\nabla](\bm\theta_i - \theta_{\ell_q})+\partial^\sigma\big((\vex_i-v_{\ell_q})\cdot\nabla \bm\theta_i \big) =-\partial^\sigma\div\MM_{\ell_q}\label{diff-2N}.
 \end{align}
Using Lemma \ref{trans} and \eqref{pex-p}, we obtain that for $t\in [t_{i-1}, t_{i+1}]$,
\begin{align*}
&\|\partial^\sigma(\vex_i - v_{\ell_q})(t)\|_{\alpha}+\|\partial^\sigma(\bm\theta_i - \theta_{\ell_q})(t)\|_{\alpha}\\
\lesssim&\int_{t_{i-1}}^t\big(\|\vex_i - v_{\ell_q}\|_{N+\alpha}\|(\vex_i, v_{\ell_q},\bm\theta_i)\|_{1+\alpha}+\|\vex_i - v_{\ell_q}\|_{\alpha}\|(\vex_i, v_{\ell_q},\bm\theta_i)\|_{N+1+\alpha}\\
&+\|\bm\theta_i - \theta_{\ell_q}\|_{N+\alpha}(\| v_{\ell_q}\|_{1+\alpha}+1)+\|(\RR_{\ell_q}, \MM_{\ell_q})\|_{N+1+\alpha}\big)\dd s.
 \end{align*}
By Gr\"{o}nwall's inequality, combining with \eqref{e:v_ell-CN+1}, \eqref{e:R_ell} and \eqref{vex-vq-alpha}, we infer from this inequality that
\begin{align*}
&\|(\vex_i - v_{\ell_q})(t)\|_{N+\alpha}+\|(\bm\theta_i - \theta_{\ell_q})(t)\|_{N+\alpha}\\
\lesssim&\int_{t_{i-1}}^{t_{i+1}}\big(\|\vex_i - v_{\ell_q}\|_{\alpha}\|(\vex_i, v_{\ell_q},\bm\theta_i)\|_{N+1+\alpha}+\|(\RR_{\ell_q}, \MM_{\ell_q})\|_{N+1+\alpha}\big)\dd s\\
&\times\exp\big(\int_{t_{i-1}}^t(\|(\vex_i, v_{\ell_q},\bm\theta_i)\|_{1+\alpha}+1)\dd s\big)\\
\lesssim& (\tau^2_q \delta_{q+1} \delta_{q}^{1 / 2} \lambda_{q}\ell^{-N-1+\frac{1}{4}\alpha}_q +\tau_q\delta_{q+1} \ell_q^{-N-1+\frac{5}{4}\alpha})\exp(\tau_q\delta_{q}^{1 / 2} \lambda_{q}\ell^{-\alpha}_q)\\
\lesssim&\tau_q\delta_{q+1} \ell_q^{-N-1+\frac{5}{4}\alpha}.
\end{align*}
Hence, we have proved \eqref{e:stability-v} for $1\le i\le i_{\max}-1$.  For $i=0$ or $i=i_{\max}$, by \eqref{def-v0}, it is easy to verify that
for $t$ satisfying $|t-t_0|\le\tau_q$, $0\le N\le N_0+1$ and large enough $a$,
\begin{align*}
&\|\vex_0 - v_{\ell_q}(t)\|_{N+\alpha}+\|\bm\theta_0 - \theta_{\ell_q}(t)\|_{N+\alpha}\\
=&\|(v^{(1)} - v^{(1)}\ast\ell_q)(t)\|_{N+\alpha}+\|(\theta^{(1)} - \theta^{(1)}\ast\ell_q)(t)\|_{N+\alpha}\\
\lesssim&\|(v^{(1)}, \theta^{(1)})\|_{N+1+\alpha}\ell_q\lesssim\tau_q\delta_{q+1}\ell^{-1-N+\frac{5}{4}\alpha}_q,
\end{align*}
and for every $t$ with $|t-t_{i_{\max}}|\le\tau_q$, $0\le N\le N_0+1$ and large enough $a$,
\begin{align*}
&\|(\vex_{i_{\max}} - v_{\ell_q})(t)\|_{N+\alpha}+\|(\bm\theta_{i_{\max}} - \theta_{\ell_q})(t)\|_{N+\alpha}\\
=&\|(v^{(2)} - v^{(2)}\ast\ell_q)(t)\|_{N+\alpha}+\|(\theta^{(2)} - \theta^{(2)}\ast\ell_q)(t)\|_{N+\alpha}\\
\lesssim&\|(v^{(2)}, \theta^{(2)})\|_{N+1+\alpha}\ell_q\lesssim\tau_q\delta_{q+1}\ell^{-1-N+\frac{5}{4}\alpha}_q.
\end{align*}
Hence, we complete the proof of \eqref{e:stability-v}. Applying  \eqref{e:v_ell-CN+1},  \eqref{e:R_ell} and \eqref{e:stability-v} to \eqref{pex-p}, we deduce \eqref{e:stability-p}. Using  \eqref{e:v_ell-CN+1}--\eqref{e:R_ell} and \eqref{e:stability-v}--\eqref{e:stability-p}, from \eqref{diff-1}--\eqref{diff-2}, one concludes \eqref{e:stability-matd}.

Next, we aim to  prove \eqref{e:stability-v-1}. For $i=0$, one immediately gets from \eqref{def-v0} that for $|t-t_0|\le \tau_q$,
\begin{align}
&\|\mathcal{R}(\vex_0 - v_{\ell_q})(t)\|_{\alpha}+\|\mathcal{R}_{\vex}(\bm\theta_0 - \theta_{\ell_q})(t)\|_{\alpha}\notag\\
=&\|\mathcal{R}v^{(1)} - \mathcal{R}v^{(1)}\ast\ell_q)(t)\|_{\alpha}+\|\mathcal{R}_{\vex}\theta^{(1)} - \mathcal{R}_{\vex}\theta^{(1)}\ast\ell_q)(t)\|_{\alpha}\notag\\
\lesssim&\|(v^{(1)}, \theta^{(1)})\|_{2+\alpha}\ell^2_q\lesssim\ell^2_q,\label{v0-vq}
\end{align}
 and  for  $i=i_{\max}$ and each $t$ with $|t-t_{i_{\max}}|\le\tau_q$,
\begin{align}
&\|\mathcal{R}(\vex_{i_{\max}} - v_{\ell_q})(t)\|_{\alpha}+\|\mathcal{R}_{\vex}(\bm\theta_{i_{\max}} - \theta_{\ell_q})(t)\|_{\alpha}\notag\\
=&\|\mathcal{R}v^{(2)} - \mathcal{R}v^{(2)}\ast\ell_q)(t)\|_{\alpha}+\|\mathcal{R}_{\vex}\theta^{(2)} - \mathcal{R}_{\vex}\theta^{(2)}\ast\ell_q)(t)\|_{\alpha}\notag\\
\lesssim&\|(v^{(2)}, \theta^{(2)})\|_{2+\alpha}\ell^2_q\lesssim\ell^2_q.\label{vmax-vq}
\end{align}
For $1\le i\le i_{\max}-1$, by {Lemma \ref{trans}} and $(\vex_i -v_{\ell_q}, \bm\theta_i -\theta_{\ell_q})(t_i)=0$, we have, for $t\in [t_i-\tau_q, t_i+\tau_q]$,
\begin{align*}
&\|(\vex_i -v_{\ell_q}, \bm\theta_i -\theta_{\ell_q})(t)\|_{B^{-1+\alpha}_{\infty, \infty}}\\
\lesssim&\int_{t_i-\tau_q}^{t}\big(\|\vex_i - v_{\ell_q}\|_{B^{-1+\alpha}_{\infty, \infty}}\|(\vex_i, v_{\ell_q}, \bm\theta_i)\|_{1}+\|\bm\theta_i - \theta_{\ell_q}\|_{B^{-1+\alpha}_{\infty, \infty}}+\|(\RR_{\ell_q}, \MM_{\ell_q})\|_{B^{\alpha}_{\infty, \infty}}\big)\dd s.
\end{align*}
Taking advantage of  Gr\"{o}nwall's inequality and \eqref{e:R_ell},  we obtain that
\begin{align}
&\|(\vex_i -v_{\ell_q}, \bm\theta_i -\theta_{\ell_q})(t)\|_{B^{-1+\alpha}_{\infty, \infty}}\nonumber\\
\lesssim&\int_{t_i-\tau_q}^{t_i+\tau_q}\|(\RR_{\ell_q}, \MM_{\ell_q})\|_{B^{\alpha}_{\infty, \infty}}\dd s
\times\exp\big(\int_{t_i-\tau_q}^{t_i-\tau_q}\|(\vex_i, v_{\ell_q}, \bm\theta_i)\|_{1}+1\dd s\big)\nonumber\\
\lesssim&\tau_q\delta_{q+1}\ell_q^{\frac{5}{4}\alpha}.\label{vi-vq}
\end{align}
Thanks to $\|(\vex_i -v_{\ell_q}, \bm\theta_i -\theta_{\ell_q})(t)\|_{B^{-1+\alpha}_{\infty, \infty}}\sim \|\mathcal{R}(\vex_i -v_{\ell_q}, \bm\theta_i -\theta_{\ell_q})(t)\|_{\alpha}$, collecting \eqref{v0-vq}--\eqref{vi-vq} together yields \eqref{e:stability-v-1}.
\end{proof}

\subsubsection{Gluing procedure}
Now, we construct a new triple $(\vv_q, \tb_q, \ppp_q)$ by gluing these exact flows $(\vex_i, \bm\theta_i, \pex_i)$. More precisely, we define intervals $I_i, J_i$ by
\begin{align*}
 &I_i \coloneq [ t_i + \tfrac{\tau_q}3,\ t_i+\tfrac{2\tau _q}3], \,\, 0\le i\le i_{\max}-1,\\
 &J_i \coloneq (t_i - \tfrac{\tau_q}3,\ t_i+\tfrac{\tau _q}3),\,\, 1\le i\le i_{\max}-1,\\
 &J_0=[0, t_0+\tfrac{\tau_q}3),\,\,J_{i_{\max}}=(t_{i_{\max}-1}+\tfrac{2\tau _q}3, 1].
\end{align*}
Let $\{ \chi_i\}_{i=0}^{i_{\textup{max}}}$  be a partition of unity on $[0,1]$ such that for all $N\ge 0$,
\begin{align*}
       &\text{supp} \chi_0=[-1, t_0+\tfrac{2\tau_q}{3}], \qquad \,\,\,\chi_0 |_{J_0} =1, \qquad\|  \del_t^N \chi_0\|_{0} \lesssim \tau_q^{-N},\\
           & \text{supp} \chi_i=I_{i-1} \cup J_i \cup I_{i},  \qquad\,\,\,\,\chi_i |_{J_i} =1, \qquad\, \|  \del_t^N \chi_i\|_{0} \lesssim \tau_q^{-N},\ \   1\le i\le i_{\max}-1,\\
       &\text{supp} \chi_{i_{\max}}= [t_{i_{\max}-1}+\tfrac{\tau_q}{3}, 2],  \,\,\,\,\chi_i |_{J_{i_{\max}}} =1, \quad\|  \del_t^N \chi_{i_{\max}}\|_{0} \lesssim \tau_q^{-N}.
\end{align*}
We define the glued {velocity, temperature  pressure }$(\vv_q, \tb_q, \ppp_q)$ by
\begin{align*}
    \vv_q(x,t) &\coloneq \sum_{i=0}^{i_{\textup{max}}} \chi_i(t) \vex_i(x,t) ,
    \\
     \tb_q(x,t) &\coloneq \sum_{i=0}^{i_{\textup{max}}} \chi_i(t) \bm\theta_i(x,t) ,
    \\
    \ppp_q(x,t) &\coloneq \sum_{i=0}^{i_{\textup{max}}} \chi_i(t) \pex_i(x,t).
\end{align*}
From the gluing procedure,  it is easy to verify  that $(\vv_q, \tb_q, \pp_q)$ solves the following equations on $[0,1]$:
\begin{equation}
\left\{ \begin{alignedat}{-1}
&\del_t \vv_q+\div (\vv_q\otimes \vv_q)  +\nabla \pp_q   =e_3\tb_q + \div \RRR_q,
 \\
 &\del_t \tb_q+\vv_q\cdot\nabla \tb_q     = \div \MMM_q,
 \\
  &\nabla \cdot \vv_q = 0,
\end{alignedat}\right.  \label{e:glu-B}
\end{equation}
where
\begin{align}
    \RRR_q &\coloneq
        \sum_{i=0}^{i_{\max}}\del_t \chi_i \mathcal R(\vex_i-\vex_{i+1} ) - \sum_{i=0}^{i_{\textup{max}}} \chi_i(1-\chi_i)(\vex_i-\vex_{i+1} )\ootimes (\vex_i-\vex_{i+1} )\notag,
\\
 \MMM_q &\coloneq
        \sum_{i=0}^{i_{\textup{max}}} \del_t \chi_i \mathcal R_{\vex}(\bm\theta_i-\bm\theta_{i+1} ) - \sum_{i=0}^{i_{\textup{max}}} \chi_i(1-\chi_i)(\vex_i-\vex_{i+1} ) (\bm\theta_i-\bm\theta_{i+1} )\notag,
\\
    \pp_q &  \coloneq \ppp_q -\int_{\TTT^3}\ppp_q \dd x- \sum_{i=0}^{i_{\textup{max}}} \chi_i(1-\chi_i)\big( |\vex_i - \vex_{i+1}|^2  - \int_{\mathbb T^3} |\vex_i - \vex_{i+1}|^2 \dd x\big).\label{ppp}
\end{align}
Particularly, $(\vv_q, \tb_q, \pp_q)$  is an exact solution to the inviscid Boussinesq system \eqref{e:B} on $J_i$.
As in the proof of \cite[Proposition 4.4]{BDSV}, we get the following estimates for  $(\RRR_q, \MMM_q)$.
\begin{prop}[Estimates for $\RRR_q$ and $\MMM_q$]\label{p:estimate-RRRq}For all $0\le N\le N_0+1$, we have
    \begin{align}
        \|(\RRR_q, \MMM_q)\|_{N+\alpha}
        &\lesssim \delta_{q+1} \ell_q^{-N+\frac{5}{4}\alpha}, \label{e:RRR_q-N+alpha-bd}
        \\
        \| (\matd{\vv_q} \RRR_q, \matd{\vv_q}\MMM_q)\|_{N+\alpha }
        &\lesssim  \delta_{q+1} \delta_{q}^{1/2}\lambda_q \ell_q^{-N-\frac{3}{4}\alpha}.\label{e:matd-RRR_q}
    \end{align}
\end{prop}
Since $(\vv_q, \tb_q)$ is defined by $(\vex_i, \bm\theta_i)$, one can obtain the following estimates for $(\vv_q, \tb_q)$ from the stability estimates in Proposition \ref{p:stability}.
\begin{prop}[Estimates for $\vv_q$ and $\tb_q$]  For  all $0\le N\le N_0+1$,
\begin{align}
&\|(\vv_{q}-v_{\ell_q}, \tb_{q}-\theta_{\ell_q})\|_{\alpha}  \le \delta_{q+1}^{1 / 2} \ell_q^\alpha,\label{e:stability-vv_q} \\
&\|(\vv_{q}, \tb_{q})\|_{0}  \le \sum_{i=1}^q\delta^{1 / 2}_i+\delta_{q+1}^{1 / 2} \ell_q^\alpha, \label{e:stability-vv_q-N} \\
&\|(\vv_{q}, \tb_q)\|_{1+N}  \lesssim \delta_{q}^{1 / 2} \lambda_{q} \ell_q^{-N}.\label{e:vv_q-bound}
\end{align}
Furthermore, we also have
\begin{align}\label{vvqtbq}
(\vv_{q}, \tb_q)\equiv (v^{(1)}, b^{(1)}) \,\,{\rm on}\,\,[0, 2\widetilde{T}+\tfrac{\tau_q}{3}],\quad (\vv_{q}, \tb_q)\equiv (v^{(2)}, b^{(2)}) \,\,\text{\rm on}\,\,[3\widetilde{T}-\tfrac{\tau_q}{3}, 1].
\end{align}
\end{prop}
\begin{proof}With aid of \eqref{e:stability-v}, we have
\begin{equation}\nonumber
\begin{aligned}
\|(\vv_{q}-v_{\ell_q}, \tb_{q}-\theta_{\ell_q})\|_{\alpha}&\le \|(\sum_{i=0}^{i_{\max}}\chi_i(t)(\vex_i-v_{\ell_q}), \sum_{i=0}^{i_{\max}}\chi_i(t)(\bm\theta_i-\theta_{\ell_q}))\|_{\alpha}\\
& \lesssim \tau_{q} \delta_{q+1} \ell_q^{-1+\alpha}\le\delta^{1/2}_{q+1}\ell_q^{\alpha}.
\end{aligned}
\end{equation}
Combining this inequality with \eqref{e:vq-C0}, one can deduce \eqref{e:stability-vv_q-N}. By \eqref{e:v_ell-CN+1}, we obtain that
\begin{equation}\nonumber
\begin{aligned}
\|(\vv_{q}, \tb_{q})\|_{N+1}\lesssim\|(\sum_{i=0}^{i_{\max}}\chi_i(t)\vex_i, \sum_{i=0}^{i_{\max}}\chi_i(t)\bm\theta_i)\|_{N+1}
\lesssim\|(v_{\ell_q}, \theta_{\ell_q})\|_{N+1}\lesssim \delta_{q}^{1 / 2} \lambda_{q} \ell_q^{-N}.
\end{aligned}
\end{equation}
By the definition of $J_0$ and $J_{i_{\max}}$, it is easy to verify that for $t\in [0, 2\widetilde{T}+\tfrac{\tau_q}{3}]$,
$$(\vv_q, \tb_q)=(\vex_0, \bm\theta_0)=(v_q, \theta_q)=(v^{(1)}, \theta^{(1)}),$$
and for $t\in [3\widetilde{T}-2\tau_q, 1]$,
$$(\vv_q, \tb_q)=(\vex_{i_{\max}},\bm\theta_{i_{\max}})=(v_q, \theta_q)=(v^{(2)}, \theta^{(2)}).$$
Thus, we conclude~\eqref{vvqtbq}.
\end{proof}

\subsection{Construction and estimates for the perturbation}
In the previous subsection, we have completed the stage: $(v_q, \theta_q)\mapsto (\vv_q, \tb_q)$. This section will be dedicated to construct the perturbation $(w_{q+1}, d_{q+1})$ of $(\vv_q, \tb_q)$ to complete this iteration. More precisely,
\[(\vv_q, \tb_q)\mapsto (v_{q+1}, \theta_{q+1}):=(\vv_q+w_{q+1}, \tb_q+d_{q+1}).\]

The main idea in the construction of the perturbation is to use 2D spatial intermittent building blocks, which are essentially consistent with Mikado flows. This highly oscillatory perturbation is designed  to cancel the  error $(\RRR_q, \MMM_q)$ by the low frequency of the term $\div(w_{q+1}\otimes w_{q+1})$. To achieve this goal,
the following geometric lemma   plays a central role in the choices of the oscillation directions of these building blocks.

\begin{lem}[Geometric Lemma\cite{2Beekie}]\label{first S}Let $B_{\sigma}({\rm Id})$ denote the ball of radius $\sigma$ centered at $\rm Id$ in the space of $3\times3$ symmetric matrices.
There exists a set $\Lambda_v\subset\mathbb{S}^2\cap\mathbb{Q}^3$ that consists of vectors $k$ with associated orthonormal basis $(k,\bar{k},\bar{\bar{k}}),~\epsilon_v>0$ and smooth function $a_{v,k}:B_{\epsilon_v}(\rm Id)\rightarrow\mathbb{R}$ such that, for every $R_v\in B_{\epsilon_v}(\rm Id)$ we have the following identity:
$$R_v=\sum_{k\in\Lambda_v}a^2_{v,k}(R_v)\bar{k}\otimes\bar{k}.$$
Furthermore, for $0\le N\le N_0+1$, there exists a constant $M$ such that
\begin{equation}\label{M}
\|a_{v,k}\|_{C^N(\bar{B}_{1/2}({\rm Id}))}\le M.
\end{equation}
\end{lem}
We define $\Lambda_{\theta}=\{e_1, e_2, e_3\}$ and  choose $\Lambda_v$ in Lemma \ref{first S} such that
$\Lambda_{\theta}\cap \Lambda_v=\emptyset$. With respect to the choice of $\Lambda_v$, readers can refer to \cite{2Beekie} for more details. Since the set $\Lambda_{v}\cup\Lambda_{\theta} \subseteq\mathbb{S}^2\cap\mathbb{Q}^3$ is finite, there exists a universal number $N_{\Lambda}$ such that
\begin{equation}\label{N}
    \{N_{\Lambda}\bar{k}, N_{\Lambda}\bar{\bar{k}}\}\subseteq N_{\Lambda}\mathbb{S}^2\cap \mathbb{N}^3, \quad \forall k\in \Lambda_{v}\cup\Lambda_{\theta}.
\end{equation}

\vskip 3mm
\noindent \textbf{Step 1:  Deformed intermittent cuboid flows}\\
Before introducing the deformed intermittent cuboid flows, we give spatial building blocks, which are essentially given in\cite{2Beekie, BCV}.
\vskip2mm
 Assume that $\phi:\mathbb{R}\rightarrow\mathbb{R}$ is a smooth mean-free cut-off function with a compact support which is much smaller than the interval $[0, 1]$,  We normalize $\phi$ such that
\begin{equation}\label{phinorm}
\int_{\mathbb{R}}\phi^2\dd x=1.
\end{equation}
Then we periodize $\phi$ so that the resulting function is a periodic function defined on $\mathbb{R}/\mathbb{Z}=\mathbb{T}$ and we still denote it by $\phi$. Now we construct spatial building blocks by $\phi$:
{\begin{lem}[Spatial building blocks]\label{SPB}Let $N_{\Lambda}$ be defined by \eqref{N} and the $\mathbb{T}$-periodic function $\phi$ support on $[0, \lambda^{-1}_{1}]$, where $\lambda_1$ is defined in \eqref{lambdaq}.  For each $k\in \Lambda_{v}\cup\Lambda_{\theta}$, let $\{k, \bar{k}, \bar{\bar{k}}\}$ be an orthonormal basis. Then for  $\lambda_1$ large enough,  there exist suitable shifts $\{x^k\}_{k\in\Lambda_{v}\cup \Lambda_{\theta}}:=\{(x^k_1, x^k_2, x^k_3)\}_{k\in\Lambda_{v}\cup \Lambda_{\theta}}$  such that
\begin{equation}
\left\{ \begin{alignedat}{-1}\label{WkWk'1}
&W^{\theta}_{(k)}(x)W^{\theta}_{{(k')}}(x)=0, \quad \forall k\neq k'\in \Lambda_{\theta},\\
&W^v_{(k)}(x)W^{v}_{{(k')}}(x)=0, \quad \forall k\neq k'\in \Lambda_v,\\
&W^v_{(k)}(x)W^{\theta}_{{(k')}}(x)=0, \quad \forall k\in \Lambda_v, k'\in\Lambda_{\theta}.
\end{alignedat}\right.
\end{equation}
That is,  the set of the supports of  $\big\{W^v_{(k)}(x)\big\}$ and $\big\{W^{\theta}_{(k)}(x)\big\}$ are mutual disjoint, where the spatial building blocks  $\big\{W^v_{(k)}(x)\big\}$ and $\big\{W^{\theta}_{(k)}(x)\big\}$ are defined as follows:
\begin{align*}
W^{\theta}_{(k)}(x)=\phi(N_{\Lambda}\bar{k}\cdot(x-x^k))\phi(N_{\Lambda}\bar{\bar{k}}\cdot(x-x^k)), \,k\in \Lambda_{\theta},\\
W^v_{(k)}(x)=\phi(N_{\Lambda}k\cdot(x-x^k))\phi(N_{\Lambda}\bar{\bar{k}}\cdot(x-x^k)), \,k\in \Lambda_{v}.
\end{align*}
\end{lem}
\begin{proof} Since $\phi(x)$ is $\mathbb{T}$-periodic and \eqref{N}, for each $k\in \Lambda_v\cup \Lambda_{\theta}$,  the function
\[W_{(k)}(x):=\phi(N_{\Lambda}\bar{k}\cdot x)\phi(N_{\Lambda}\bar{\bar{k}}\cdot x)\]
is $\mathbb{T}^3$-periodic. The construction of  $W_{(k)}(x)$ reveals that it supports on finite cuboids of which the length side are parallel to  $k$, the width and height are equivalent to  $\lambda^{-1}_1$. Note that $\Lambda_v{\cup}\Lambda_{\theta}\subseteq \mathbb{Q}^3$ is finite,  we take $\lambda_1$ large enough such that these finite cuboids  are mutually disjoint by  choosing  suitably shifts  $\{x^k\}_{k\in\Lambda_v{\cup}\Lambda_{\theta}}$, which shows  \eqref{WkWk'1}.
\end{proof}}
For the case when  $W_{(k)}(x)$ supports on pipes,  readers can refer to \cite{BV-19} for the proof. For the sake of dealing with the transport error of the perturbation to achieve the optimum regularity $\beta$, we will construct the deformed intermittent cuboid flows by the spatial building blocks in  Lagrangian coordinates $\Phi_i$, the so-called inverse flow map defined as follows.
\vskip 3mm
\noindent\textbf{Inverse flow map.}\quad
Let $\Phi_i$ be  the solution to the glued vector transport equations
    \begin{align}
         (\partial_t+\vv_q \cdot \nabla) \Phi_i = 0, \quad \Phi_i\big|_{t=t_i} = x. \label{e:phi_i-defn}
    \end{align}
We define  time cutoffs $\eta_i\in C^{\infty}_c(\R)$  by
\begin{equation}\label{etai}
{\rm supp} \,\,\eta_i=I_i+[-\tfrac{\tau_q}{6}, \tfrac{\tau_q}{6}], \quad \eta_i|_{I_i}\equiv 1,\quad\|\partial^N_t\eta_i\|_0\lesssim \tau^{-N}_q.
\end{equation}
By the estimates for $\vv_q$ in \eqref{e:vv_q-bound}, we have the following stability and derivative estimates of $\Phi_i$ in a short time. Readers can refer to \cite[Proposition 4.10]{KMY} for the details.
\begin{prop}[Estimates for $\Phi_i$\cite{KMY}]\label{p:estimates-for-inverse-flow-map} For $a\gg 1$, $0\le N\le N_0+1$ and every $t\in {\rm supp} \,\,\eta_i$,
\begin{align}
 \|\nabla\Phi_i-{\rm Id}_{3\times3}\|_0 &\le\frac1{10},\label{e:nabla-phi-i-minus-I3x3}
 \\
 \| (\nabla \Phi_i)^{-1}\|_N+ \|  \nabla \Phi_i\|_N &\le \ell_q^{-N}, \label{e:nabla-phi-i-CN}
\\
\|\matd{\vv_q} \nabla \Phi_i\|_N &\lesssim \delta_q^{1/2} \lambda_q \ell_q^{-N}. \label{e:nabla-phi-i-matd}
\end{align}
\end{prop}
Noting that $\Phi_i(x,t)$ is a measure preserving map,  $\phi(N_{\Lambda}k\cdot(\Phi_i-x^k))$,  $\phi(N_{\Lambda}\bar{k}\cdot(\Phi_i-x^k))$ and  $\phi(N_{\Lambda}\bar{\bar{k}}\cdot(\Phi_i-x^k))$ are  mean-free functions on $\TTT^3$. Furthermore, we follow from
\[\Phi_i(x+\ZZ^3,t)=\Phi_i(x,t)+\ZZ^3\]
that
$\phi(N_{\Lambda}k\cdot(\Phi_i-x_k))$, $\phi(N_{\Lambda}\bar{k}\cdot(\Phi_i-x_k))$ and  $\phi(N_{\Lambda}\bar{\bar{k}}\cdot(\Phi_i-x_k))$  are $\TTT^3$-periodic functions.
 Consequently, we conclude that $W^v_{(k)}(\Phi_i(x,t))$ and  $W^{\theta}_{(k)}(\Phi_i(x,t))$ are $\TTT^3$-periodic function with zero mean. Moreover, we have the following property.
{\begin{lem}For $t\in \text{supp}\, \eta_i(t)$, we have
\begin{equation}
\left\{ \begin{alignedat}{-1}\label{WkWk'}
&W^{\theta}_{(k)}(\Phi_i(x,t))W^{\theta}_{{(k')}}(\Phi_i(x,t))=0, \quad \forall k\neq k'\in \Lambda_{\theta},\\
&W^v_{(k)}(\Phi_i(x,t))W^{v}_{{(k')}}(\Phi_i(x,t))=0, \quad \forall k\neq k'\in \Lambda_v,\\
&W^v_{(k)}(\Phi_i(x,t))W^{\theta}_{{(k')}}(\Phi_i(x,t))=0, \quad \forall k\in \Lambda_v, k'\in\Lambda_{\theta}.
\end{alignedat}\right.
\end{equation}
\end{lem}
\begin{proof}By the definition of $\Phi_i$ in \eqref{e:phi_i-defn}, we deduce that
$\Phi_i(x,t)=\psi^{-1}_t(x)$,  the inverse mapping for the flow $\psi_{t}(x)$, where
\begin{equation}\label{eq-psi}
    \frac{\mathrm{d}}{\mathrm{d} t}\psi_{t}(x)=\bar{v}_q(\psi_t(x),t), \qquad \psi_{t}(x)|_{t=t_i}=x.
\end{equation} For any closed set $A, B\subseteq \mathbb{T}^3$ satisfying $A\cap B=\emptyset$, we have $|x-y|>0$, $\forall x\in A, \forall y\in B$. Next, we consider $d(t)=|\psi_t(x)-\psi_t(y)|$. From the equation \eqref{eq-psi}, we have
\begin{align*}
\frac{\mathrm{d}}{\mathrm{d} t}(\psi_t(x)-\psi_t(y))=\bar{v}_q(\psi_t(x),t)-\bar{v}_q(\psi_t(y),t).
\end{align*}
This together with $\|\nabla\bar{v}_q\|_{L^\infty}\le C\delta^{1/2}_{q}\lambda_q$  yields  that
\begin{align*}
\Big|\frac{\mathrm{d}}{\mathrm{d} t}d(t)\Big|\le \|\nabla \bar{v}_q\|_{L^\infty} d(t)\le C\delta^{1/2}_{q}\lambda_q d(t),
\end{align*}
which implies that $\frac{\mathrm{d}}{\mathrm{d} t}d(t)\ge - C\delta^{1/2}_{q}\lambda_q d(t)$.
Then for $t\in \text{supp}\, \eta_i(t)=[t_i+\frac{\tau_q}{6}, t_i+\frac{5\tau_q}{6}]$, we obtain that
\begin{align*}
d(t)\ge d(0)e^{-C\delta^{1/2}_{q}\lambda_q\tau_q}=|x-y|e^{-C\delta^{1/2}_{q}\lambda_q\tau_q}=|x-y|e^{-C_0}>0.
\end{align*}
This shows that $\psi_t(A)\cap \psi_t(B)=\emptyset$, that is, $\Phi^{-1}_i(A,t)\cap\Phi^{-1}_i(B,t)=\emptyset$. Therefore, the set of the supports of  $\big\{W^v_{(k)}(x)\big\}$ and $\big\{W^{\theta}_{(k)}(x)\big\}$ are mutual disjoint implies that the supports of  $\big\{W^v_{(k)}(\Phi_i(x,t))\big\}$ and $\big\{W^{\theta}_{(k)}(\Phi_i(x,t))\big\}$ are also mutual disjoint for $t\in  \text{supp}\, \eta_i(t)$, which implies~\eqref{WkWk'}.
\end{proof}}
\vskip 3mm
We are in position to construct the deformed intermittent cuboid flows. For each $k\in \Lambda_{\theta}$, we define
\begin{equation}\nonumber
 \begin{aligned}
W^{\theta}_{(\lambda_{q+1}k)}(\Phi_i(x,t))k:=\phi(\lambda_{q+1}N_{\Lambda}\bar{k}\cdot(\Phi_i(x,t)-x^k))\phi(\lambda_{q+1}N_{\Lambda}\bar{\bar{k}}\cdot(\Phi_i(x,t)-x^k))k.
 \end{aligned}
 \end{equation}
For each $k\in \Lambda_v$, we define
\begin{equation}\nonumber
 \begin{aligned}
W^v_{(\lambda_{q+1}k)}(\Phi_i(x,t))\bar{k}:=\phi(\lambda_{q+1}N_{\Lambda}k\cdot(\Phi_i(x,t)-x^k))\phi(\lambda_{q+1}N_{\Lambda}\bar{\bar{k}}\cdot(\Phi_i(x,t)-x^k))\bar{k}.
 \end{aligned}
 \end{equation}

 Let us rewrite $W^v_{(\lambda_{q+1}k)}(\Phi_i(x,t))$ and $W^{\theta}_{(\lambda_{q+1}k)}(\Phi_i(x,t))$ as
\begin{equation}\label{W}
\left\{ \begin{alignedat}{-1}
W^{\theta}_{(\lambda_{q+1}k)}(\Phi_i(x,t)):=\sum_{l,m\in\ZZ\backslash\{0\}}a_{k,l}a_{k,m}e^{2\pi \ii \lambda_{q+1}N_{\Lambda}(l \bar{k}+m\bar{\bar{k}})\cdot\Phi_i}:=\sum_{j\in I^{\theta}_k}b_{k,j}e^{2\pi\ii\lambda_{q+1}j\cdot\Phi_i},\\
W^{v}_{(\lambda_{q+1}k)}(\Phi_i(x,t)):=\sum_{l,m\in\ZZ\backslash\{0\}}a_{k,l}a_{k,m}e^{2\pi \ii \lambda_{q+1}N_{\Lambda}(l k+m\bar{\bar{k}})\cdot\Phi_i}:=\sum_{j\in I^v_k}b_{k,j}e^{2\pi\ii\lambda_{q+1}j\cdot\Phi_i},
 \end{alignedat}\right.
 \end{equation}
 where $I^{\theta}_k$ and  $I^v_k$ are defined by
 \begin{align*}
     &I^{\theta}_k=\big\{j=N_{\Lambda}l\bar{k}+N_{\Lambda}m\bar{\bar{k}}\,\,|\,\, l, m\in\ZZ\backslash\{0\}, k\in \Lambda_{\theta}\big\},\\
    & I^v_k=\big\{j=N_{\Lambda}lk+N_{\Lambda}m\bar{\bar{k}}\,\,|\,\, l, m\in\ZZ\backslash\{0\}, k\in \Lambda_v\big\}.
 \end{align*}
From the smoothness of $W^{\theta}_{(k)}$ and  $W^v_{(k)}$, we can infer that for $j\in I^{\theta}_k\cup  I^{v}_k$,
\begin{equation}\label{bkj}
|b_{k,j}|\lesssim|lm|^{-6}\lesssim|j|^{-6}.
\end{equation}
Furthermore, $W^v_{(\lambda_{q+1} k)}(\Phi_i(x,t))$ and  $W^{\theta}_{(\lambda_{q+1} k)}(\Phi_i(x,t))$ are $\TTT^3$-periodic function with zero mean. By \eqref{phinorm}, we have
 \begin{align}\label{WL2}
 \big\|W^v_{(k)}(\Phi_i(x,t))\big\|_{L^2(\TTT^3)} =\big\|W^{\theta}_{(k)}(\Phi_i(x,t))\big\|_{L^2(\TTT^3)}=1.
 \end{align}
\noindent\textbf{Step 2:\,\,Construction of the perturbation}\\
Now, we are focused on constructing the perturbation parts $(w_{q+1}, d_{q+1})$ by  \eqref{W}. First of all, we rewrite the vector function $\nabla\Phi_i\MMM_q\eta_i(t)$ in components form as follows.
\begin{equation}\label{Tqik}
-\nabla\Phi_i\MMM_q\eta_i(t)=\sum_{j=1}^3\big(-\nabla\Phi_i\MMM_q\big)^j \eta_i(t)e_j:=\sum_{k\in\Lambda_{\theta}}\mathcal{T}_{q,i,k}\eta_i(t)k,
\end{equation}
where we have used the fact that $\Lambda_{\theta}=\{e_1, e_2, e_3\}$.
The \emph{principal part} of the perturbation of glued velocity is defined by
\begin{align}
w^{(p)}_{q+1}=&\sum_{i=0}^{i_{\max}-1}w^{(p)}_{q+1,i}\notag\\
:=&\sum_{i=0}^{i_{\max}-1}\big(\delta^{-{1/2}}_{q+1}\ell^{-\alpha/2}_q\sum_{k\in\Lambda_{\theta}}\mathcal{T}_{q,i,k}\eta_i(t)\nabla\Phi^{-1}_iW^{\theta}_{(\lambda_{q+1}k)}(\Phi_i(x,t))k\notag\\
&+{\delta^{{1/2}}_{q+1}}\ell^{\alpha/2}_q\sum_{k\in\Lambda_{v}}a_{v,k}\big({R_v}\big)\eta_i(t)\nabla\Phi^{-1}_iW^v_{(\lambda_{q+1}k)}(\Phi_i(x,t))\bar{k}\big),\label{wp-def1}
\end{align}
where
\begin{equation}\label{Rv}
 R_v=\nabla\Phi_i\big(\mathrm{Id}-\tfrac{\RRR_q}{\delta_{q+1}\ell^{\alpha}_q}\big)(\nabla\Phi_i)^{\TT}+\delta^{-1}_{q+1}\ell^{-\alpha}_q\sum_{k\in\Lambda_{\theta}}\mathcal{T}^2_{q,i,k}\eta^2_i(t)k\otimes k.
 \end{equation}
According to  the representation \eqref{W}, we can rewrite $w^{(p)}_{q+1}$ as follows:
 \begin{align}\label{wp_q+1}
w^{(p)}_{q+1}=&\sum_{i=0}^{i_{\max}-1}w^{(p)}_{q+1,i}\notag\\
=&\sum_{i=0}^{i_{\max}-1}\Big(\sum_{k\in\Lambda_{\theta}}\sum_{j\in I^{\theta}_k}A_{k,i,j}ke^{2\pi\ii\lambda_{q+1}j\cdot\Phi_i}+\sum_{k\in\Lambda_{v}}\sum_{j\in I^v_k}B_{k,i,j}\bar{k}e^{2\pi\ii\lambda_{q+1}j\cdot\Phi_i}\Big),
\end{align}
where
\begin{equation}\label{AB}
\left\{ \begin{alignedat}{-1}
&A_{k,i,j}=\delta^{-{1/2}}_{q+1}\ell^{-\alpha/2}_qb_{k,j}\mathcal{T}_{q,i,k}\eta_i(t)\nabla\Phi^{-1}_i,\\
&B_{k,i,j}={\delta^{{1/2}}_{q+1}}\ell^{\alpha/2}_qb_{k,j}a_{v,k}({R_v})\eta_i(t)\nabla\Phi^{-1}_i.
\end{alignedat}\right.
\end{equation}
 The  perturbation of glued temperature is constructed  as follows:
\begin{align}
d_{q+1}=\sum_{i=0}^{i_{\max}-1}d_{q+1,i}:=&\sum_{i=0}^{i_{\max}-1}\sum_{k\in\Lambda_{\theta}}\delta^{{1/2}}_{q+1}\ell^{\alpha/2}_qW^{\theta}_{(\lambda_{q+1}k)}(\Phi_i(x))\eta_i(t)\nonumber\\
=&\sum_{i=0}^{i_{\max}-1}\sum_{k\in\Lambda_{\theta}}\sum_{j\in I^{\theta}_k}\delta^{{1/2}}_{q+1}\ell^{\alpha/2}_qb_{k,j}\eta_i(t)e^{2\pi\ii\lambda_{q+1}j\cdot\Phi_i}.\label{d}
 \end{align}

 Next, we construct the \emph{incompressibility corrector} to guarantee that the perturbation of glued velocity is divergence-free. To obtain the incompressibility corrector, we firstly give an equivalent form of $w^{(p)}_{q+1}$.

For any  vectors $j$ and $k$ satisfying $j\perp k$,  one deduces that
 \begin{align*}
 {2\pi}^{-1}\lambda^{-1}_{q+1}\curl_{\xi}\Big(\Big(\frac{\ii j\times k}{|j|^2}\Big)e^{2\pi\ii\lambda_{q+1}j\cdot\xi}\Big)=\ii j\times\Big(\frac{\ii j\times k}{|j|^2}\Big)e^{2\pi\ii\lambda_{q+1}j\cdot\xi}=ke^{2\pi\ii\lambda_{q+1}j\cdot\xi},
 \end{align*}
This fact combined with
$$j\perp k, \quad\forall\,\, k\in \Lambda_{\theta},\, j\in I^{\theta}_k; \quad j\perp \bar{k}, \quad\forall\,\, k\in \Lambda_{v},\, j\in I^{v}_k$$
 allows us to obtain from \eqref{wp_q+1} that
 \begin{align*}
 w^{(p)}_{q+1,i}=&\frac{\delta^{-{1/2}}_{q+1}\ell^{-\alpha/2}_q}{2\pi\lambda_{q+1}}\sum_{k\in\Lambda_{\theta}}\sum_{j\in I^{\theta}_k}\mathcal{T}_{q,i,k}\nabla\Phi^{-1}_i\curl_{\xi}\Big(\Big(\frac{\ii j\times b_{k,j}k}{|j|^2}\Big)e^{2\pi\ii\lambda_{q+1}j\cdot\xi}\Big)(\Phi_i)\\
 &+ \frac{\delta^{{1/2}}_{q+1}\ell^{\alpha/2}_q}{2\pi\lambda_{q+1}}\sum_{k\in\Lambda_{v}}\sum_{j\in I^v_k}a_{v,k}\big({R_v}\big)\nabla\Phi^{-1}_i\curl_{\xi}\Big(\Big(\frac{\ii j\times b_{k,j}\bar{k}}{|j|^2}\Big)e^{2\pi\ii\lambda_{q+1}j\cdot\xi}\Big)(\Phi_i).
 \end{align*}
Combining with
 \[\nabla \Phi^{-1}_i(\curl U)(\Phi_i)=\curl(\nabla\Phi^{\TT}_iU(\Phi_i)),\]
 we  rewrite $w^{(p)}_{q+1,i}$ as follows:
 \begin{align*}
 w^{(p)}_{q+1,i}=&\frac{\delta^{-{1/2}}_{q+1}\ell^{-\alpha/2}_q}{2\pi\lambda_{q+1}}\sum_{k\in\Lambda_{\theta}}\sum_{j\in I^{\theta}_k}\mathcal{T}_{q,i,k}\curl\Big(\nabla\Phi^{\TT}_i\Big(\frac{\ii j\times b_{k,j}k}{|j|^2}\Big)e^{2\pi\ii\lambda_{q+1}j\cdot\Phi_i}\Big)\\
 &+\frac{\delta^{{1/2}}_{q+1}\ell^{\alpha/2}_q}{2\pi\lambda_{q+1}}\sum_{k\in\Lambda_{v}}\sum_{j\in I^v_k}a_{v,k}\big({R_v}\big)\curl\Big(\nabla\Phi^{\TT}_i\Big(\frac{\ii j\times b_{k,j}\bar{k}}{|j|^2}\Big)e^{2\pi\ii\lambda_{q+1}j\cdot\Phi_i}\Big).
 \end{align*}
 This equivalent form of $w^{(p)}_{q+1,i}$ helps us define
\begin{equation}\label{wc}
 \begin{aligned} w^{(c)}_{q+1,i}=&\sum_{k\in\Lambda_{\theta}}\sum_{j\in I^{\theta}_k}E_{k,i,j}e^{2\pi\ii\lambda_{q+1}j\cdot\Phi_i}+\sum_{k\in\Lambda_{v}}\sum_{j\in I^v_k}F_{k,i,j}e^{2\pi\ii\lambda_{q+1}j\cdot\Phi_i},
 \end{aligned}
 \end{equation}
 where
 \begin{equation}
\left\{ \begin{alignedat}{-1}\label{EF}
 &E_{k,i,j}=\frac{\delta^{-{1/2}}_{q+1}\ell^{-\alpha/2}_q}{2\pi\lambda_{q+1}}\nabla\mathcal{T}_{q,i,k}\times\Big(\nabla\Phi^{\TT}_i\Big(\frac{\ii j\times b_{k,j}k}{|j|^2}\Big)\Big),\\
 &F_{k,i,j}=\frac{\delta^{{1/2}}_{q+1}\ell^{\alpha/2}_q}{2\pi\lambda_{q+1}}\nabla \big(a_{v,k}\big({R_v})\big)\times\Big(\nabla\Phi^{\TT}_i\Big(\frac{\ii j\times b_{k,j}\bar{k}}{|j|^2}\Big)\Big).
 \end{alignedat}\right.
 \end{equation}
Since $\curl(fW)=\nabla f\times W+f\curl W$, we get that $\div(w^{(p)}_{q+1,i}+w^{(c)}_{q+1,i})=0$. Finally, we construct $w_{q+1}$ as
\[w_{q+1}=\sum_{i=0}^{i_{\max}-1}(w^{(p)}_{q+1,i}+w^{(c)}_{q+1,i}).\]
\vskip 3mm
\noindent{\textbf{Step 3:\,\,Preliminaries for the estimates for perturbation}}\\
In order to estimate $(w_{q+1}, d_{q+1})$, we give the derivative estimates of $R_v$ for preparation.
\begin{prop}[Estimates for $R_v$]For $0\le N\le N_0+1$, we have
\begin{align}
\|R_v\|_{N}&\lesssim\delta_{q+1}\ell^{-N+\frac{5}{4}\alpha}_q,\label{estimate-Rv}\\
\|\matd {\vv_q}R_v\|_{N}&\lesssim\delta_{q+1}\delta^{1/2}_q\lambda_q\ell^{-N-\frac{3}{4}\alpha}_q \label{estimate-DRv}.
\end{align}
\end{prop}
\begin{proof}
By the definition of $\mathcal{T}_{q,i,k}$ in \eqref{Tqik}, we obtain from \eqref{e:RRR_q-N+alpha-bd} and \eqref{e:nabla-phi-i-CN}  that
\begin{align}\label{TqikN}
&\|\mathcal{T}_{q,i,k}\|_{N}\le\|\MMM_q\|_N\|\nabla\Phi_i\|_0+\|\MMM_q\|_0\|\nabla\Phi_i\|_N\lesssim\delta_{q+1} \ell_q^{-N+\frac{5}{4}\alpha}.
\end{align}
Combining \eqref{TqikN} with \eqref{e:RRR_q-N+alpha-bd} and \eqref{e:nabla-phi-i-CN}, we infer from the definition of $R_v$ in \eqref{Rv} that
\begin{align*}
\|R_v\|_{N}&\lesssim\|\nabla\Phi_i\|_N\|\RRR_q\|_0+\|\RRR_q\|_N+\delta^{-1}_{q+1}\ell^{-\alpha}_q\|\mathcal{T}_{q,i,k}\|_N\|\mathcal{T}_{q,i,k}\|_0\\
&\lesssim\delta_{q+1} \ell_q^{-N+\frac{5}{4}\alpha}+\delta^{-1}_{q+1}\cdot\delta^2_{q+1} \ell_q^{-N+\frac{5}{4}\alpha}\lesssim\delta_{q+1}\ell^{-N+\frac{5}{4}\alpha}_q.
\end{align*}
Taking advantage of \eqref{e:RRR_q-N+alpha-bd}, \eqref{e:matd-RRR_q}, \eqref{e:nabla-phi-i-CN} and \eqref{e:nabla-phi-i-matd}, one deduces that
\begin{align}
\|\matd {\vv_q}\mathcal{T}_{q,i,k}\|_{N}\le&\|\matd {\vv_q}\MMM_q\|_N\|\eta_i(t)\nabla\Phi_i\|_0+\|\matd {\vv_q}\MMM_q\|_0\|\eta_i(t)\nabla\Phi_i\|_N\nonumber\\
&+\|\eta_i(t)\matd {\vv_q}\nabla\Phi_i\|_N\|\MMM_q\|_0+\|\eta_i(t)\matd {\vv_q}\nabla\Phi_i\|_0\|\MMM_q\|_N\nonumber\\
&+\|\eta'(t)\nabla\Phi_i\|_0\|\MMM_q\|_N+\|\eta'(t)\nabla\Phi_i\|_N\|\MMM_q\|_0\nonumber\\
\lesssim&\delta_{q+1}\delta^{1/2}_q\lambda_q\ell^{-N-\frac{3}{4}\alpha}_q.\label{matdT}
\end{align}
This estimate together with \eqref{e:matd-RRR_q}, \eqref{e:nabla-phi-i-matd} and \eqref{TqikN} allow us to infer that
\begin{align*}
\|\matd {\vv_q}R_v\|_{N}&\lesssim\|(\matd {\vv_q}\nabla\Phi_i) \RRR_q\|_N+\|\nabla\Phi_i \matd {\vv_q}\RRR_q\|_N\\
&\quad+\delta^{-1}_{q+1}\ell^{-\alpha}_q({\|\matd {\vv_q}\mathcal{T}_{q,i,k}\|_{N}\|\mathcal{T}_{q,i,k}\|_0+\|\matd {\vv_q}\mathcal{T}_{q,i,k}\|_{0}\|\mathcal{T}_{q,i,k}\|_N})\\
&\lesssim\delta_{q+1}\delta_q^{1/2} \lambda_q \ell_q^{-N+\frac{5}{4}\alpha} +\delta_{q+1} \delta_{q}^{1/2}\lambda_q \ell_q^{-N-\frac{3}{4}\alpha} +\delta_{q+1} \delta_{q}^{1/2}\lambda_q \ell_q^{-N-\frac{1}{2}\alpha}\\
&\lesssim\delta_{q+1}\delta^{1/2}_q\lambda_q\ell^{-N-\frac{3}{4}\alpha}_q.
\end{align*}
\end{proof}
 Next, we give the derivative estimates for these coefficients functions $A_{k,i,j}, B_{k,i,j}$, $E_{k,i,j}$, $F_{k,i,j}$ which are defined in \eqref{AB} and \eqref{EF}.
 \begin{prop}[$C^N$ estimates for $A_{k,i,j}$, $B_{k,i,j}$, $E_{k,i,j}$ and $F_{k,i,j}$]For $j\in I^{\theta}_k\cup I^v_k$ and $0\le N\le N_0$, we have
 \begin{align}
 &\|A_{k,i,j}\|_{N}
\lesssim\delta^{1/2}_{q+1}\ell^{-N+\frac{1}{4}\alpha}_q|j|^{-6}, \label{AN}\\
&\|B_{k,i,j}\|_{N}\lesssim\delta^{1/2}_{q+1}\ell^{-N+\frac{1}{4}\alpha }_q|j|^{-6}, \label{BN}\\
 &\|E_{k,i,j}\|_{N}
\lesssim\delta^{1/2}_{q+1}\lambda^{-1}_{q+1}\ell^{-(N+1)+\frac{1}{4}\alpha}_q|j|^{-7},\label{EN}\\
&\|F_{k,i,j}\|_{N}
\lesssim\delta^{1/2}_{q+1}\lambda^{-1}_{q+1}\ell^{-(N+1)+\frac{1}{4}\alpha}_q|j|^{-7}.\label{FN}
 \end{align}
\end{prop}
\begin{proof}
Using \eqref{TqikN} and \eqref{e:nabla-phi-i-CN},  we directly obtain from \eqref{AB} that
\begin{align*}
\|A_{k,i,j}\|_{N}& \lesssim \delta^{-1/2}_{q+1}\ell^{-\alpha/2}_q|b_{k,j}|\big(\|\mathcal{T}_{q,i,k}\|_{N}\|\nabla\Phi_i\|_{0}+\|\mathcal{T}_{q,i,k}\|_{0}\|\nabla\Phi_i\|_{N}\big)\\
&\lesssim\delta^{1/2}_{q+1}\ell^{-N+\frac{3}{4}\alpha}_q|j|^{-6}.
\end{align*}
With the aid of \eqref{M} and \eqref{estimate-Rv}, for $0\le N\le N_0+1$, we have
\begin{equation}\label{avkN}
\big\|a_{(v,k)}(R_v)\big\|_{N}\lesssim_N M\big(\|R_v\|^N_1+\|R_v\|_N+1\big)\lesssim_N M (\ell^{-N+\frac{\alpha}{4}}_q+1).
\end{equation}
This inequality  together with \eqref{e:nabla-phi-i-CN}  and  \eqref{bkj} yields that
 \begin{align*}
 \|B_{k,i,j}\|_{N}&\lesssim\delta^{1/2}_{q+1}\ell^{\alpha/2}_q|b_{k,j}|(\|a_{v,k}(R_v)\|_N\|\eta_i(t)\nabla\Phi^{-1}_i\|_0+M\|\eta_i(t)\nabla\Phi^{-1}_i\|_N\big)\\
 &\lesssim_N M\delta^{1/2}_{q+1}(\ell^{-N+\frac{3}{4}\alpha}_q+\ell^{\alpha/2}_q+\ell^{-N+\frac{1}{2}\alpha}_q)|j|^{-6}\lesssim\delta^{1/2}_{q+1}\ell^{-N+\frac{\alpha}{4}}_q|j|^{-6}.
 \end{align*}
 Taking advantage of  \eqref{e:nabla-phi-i-CN}, \eqref{bkj} and \eqref{TqikN},  we get from the definition of $E_{k,i,j}$ in \eqref{EF} that, for any $0\le N\le N_0$,
\begin{align*}
\|E_{k,i,j}\|_{N}&\lesssim\tfrac{\delta^{-1/2}_{q+1}\ell^{-\alpha/2}_q}{\lambda_{q+1}}|b_{k,j}||j|^{-1}(\|\nabla\mathcal{T}_{q,i,k}\|_N\|\nabla \Phi^{\TT}_i\|_0+\|\nabla\mathcal{T}_{q,i,k}\|_0\|\nabla\Phi_i^{\TT}\|_N)\\
&\lesssim\delta^{1/2}_{q+1}\lambda^{-1}_{q+1}\ell^{-(N+1)+\frac{3}{4}\alpha}_q|j|^{-7}.
\end{align*}
By the definition of $F_{k,i,j}$ in \eqref{EF}, one deduces from \eqref{avkN} that, for $0\le N\le N_0$,
\begin{align*}
\|F_{k,i,j}\|_{N}&\lesssim\tfrac{\delta^{1/2}_{q+1}\ell^{\alpha/2}_q}{\lambda_{q+1}}|b_{k,j}||j|^{-1}\big(\big\|a_{v,k}({R_v})\big\|_{N+1}+\big\|a_{(v,k)}(R_v)\big\|_{1}\|\nabla\Phi_i\|_N\big)\\
&\lesssim M\delta^{1/2}_{q+1}\lambda^{-1}_{q+1}\ell^{-(N+1)+\frac{3\alpha}{4}}_q|j|^{-7}\lesssim \delta^{1/2}_{q+1}\lambda^{-1}_{q+1}\ell^{-(N+1)+\frac{5}{8}\alpha}_q|j|^{-7}.
\end{align*}
For the sake of estimating the transports errors in the new Reynolds and temperature stresses, we need the following estimates of material derivative.
\end{proof}
\begin{prop}[Estimates for material derivative of $A_{k,i,j}$, $B_{k,i,j}$, $E_{k,i,j}$ and $F_{k,i,j}$]For $j\in I^{\theta}_k\cup I^v_k$ and $0\le N\le N_0$, we have
 \begin{align}
 &\|\matd {\vv_q}A_{k,i,j}\|_{N}
\lesssim\delta^{1/2}_{q+1} \delta^{1/2}_q\lambda_q\ell_q^{-N-\frac{3}{2}\alpha}|j|^{-6}, \label{DAN}\\
&\|\matd {\vv_q}B_{k,i,j}\|_{N}\lesssim\delta^{1/2}_{q+1} \delta^{1/2}_q \lambda_q\ell^{-N-\frac{3}{2}\alpha}_q|j|^{-6},\label{DBN}\\
 &\|\matd {\vv_q}E_{k,i,j}\|_{N}
\lesssim\delta^{1/2}_{q+1} \delta^{1/2}_q\lambda_q \lambda^{-1}_{q+1} \ell_q^{-N-1-\frac{3}{2}\alpha}|j|^{-7},\label{DEN}\\
&\|\matd {\vv_q}F_{k,i,j}\|_{N}
\lesssim\delta^{1/2}_{q+1} \delta^{1/2}_q\lambda_q \lambda^{-1}_{q+1} \ell_q^{-N-1-\frac{3}{2}\alpha}|j|^{-7}.\label{DFN}
 \end{align}
\end{prop}
\begin{proof}Making use of \eqref{e:matd-RRR_q} and \eqref{e:nabla-phi-i-matd}, we obtain that
\begin{equation}\label{matdvvq}
\begin{aligned}
\|\matd {\vv_q}\mathcal{T}_{k,i,j}\|_N&\le\|(\matd {\vv_q}\nabla\Phi_i )\MMM_q\eta_i(t)\|_N+\|(\matd {\vv_q}\MMM_q)\nabla\Phi_i \eta_i(t)\|_N+\|\MMM_q\nabla\Phi_i \eta_i'(t)\|_N\\
&\lesssim\delta_{q+1} \delta_{q}^{1/2}\lambda_q \ell_q^{-N-\frac{5}{4}\alpha}+\tau^{-1}_q\delta_{q+1} \ell_q^{-N+\frac{3}{4}\alpha}\lesssim\delta_{q+1} \delta_{q}^{1/2}\lambda_q \ell_q^{-N-\frac{3}{4}\alpha}.
\end{aligned}
\end{equation}
The above inequality together with \eqref{e:vv_q-bound},  \eqref{e:nabla-phi-i-matd}, \eqref{bkj} and \eqref{matdT} leads to
\begin{align*}
\|\matd {\vv_q} A_{k,i,j}\|_{N}&\le \delta^{-1/2}_{q+1}\ell^{-\alpha/2}_q|b_{k,j}|(\|(\matd {\vv_q}\mathcal{T}_{k,i,j})\eta_i(t)\nabla\Phi_i\|_{N}\\
&\quad+\|(\matd {\vv_q}\nabla\Phi_i)  \mathcal{T}_{k,i,j}\eta_i(t)\|_{N}+\|\mathcal{T}_{k,i,j}\eta'_i(t)\nabla\Phi_i\|_{N})\\
&\lesssim\delta^{1/2}_{q+1} \delta_{q}^{1/2}\lambda_q \ell_q^{-N-\frac{3\alpha}{2}}|j|^{-6}.
\end{align*}
Using \eqref{estimate-DRv} and \eqref{avkN} leads to
\begin{equation}\label{Davk}
\begin{aligned}
&\big{\|}\matd {\vv_q} \big(a_{v,k}(R_v)\big)\big{\|}_N=\big\|\big((a'_{v,k})(R_v)\big)\matd {\vv_q}\big(\tfrac{ R_{v}}{\delta_{q+1}\ell^{\alpha}_q}\big)\big\|_N\\
\le&\big\|(a'_{v,k})(R_v)\|_0\big\|\matd {\vv_q}\big(\tfrac{ R_{v}}{\delta_{q+1}\ell^{\alpha}_q}\big)\big\|_N
+\big\|(a'_{v,k})(R_v)\|_N\|\matd {\vv_q}\big(\tfrac{ R_{v}}{\delta_{q+1}\ell^{\alpha}_q}\big)\big\|_0\\
\lesssim&M\delta^{1/2}_q\lambda_q\ell^{-N-\frac{3}{2}\alpha}_q\lesssim\delta^{1/2}_q\lambda_q\ell^{-N-2\alpha}_q.
\end{aligned}
\end{equation}
Collecting this and \eqref{estimate-DRv}, one can infer that for $0\le N\le N_0$,
\begin{align*}
&\|\matd {\vv_q} B_{k,i,j}\|_{N}\\
\le&\delta^{1/2}_{q+1}\ell^{\alpha/2}_q|b_{k,j}|\big(\|\matd {\vv_q} \big(a_{v,k}(R_v)\big)\eta_i\nabla \Phi^{-1}_i\|_N+\|a_{v,k}(R_v)\eta_i\matd {\vv_q}\nabla \Phi^{-1}_i\|_N\\
&+\|\eta'_i(t)a_{v,k}(R_v)\nabla\Phi^{-1}_i\|_N\big)\\
\lesssim&\delta^{1/2}_{q+1}\ell^{\alpha/2}_q|j|^{-6}(\delta^{1/2}_q\lambda_q\ell^{-N-2\alpha}_q+M\delta^{1/2}_{q+1}\lambda_q\ell^{-N}_q+M\tau^{-1}_q\ell^{-N}_q)\lesssim\delta^{1/2}_{q+1}\delta^{1/2}_q\lambda_q\ell^{-N-\frac{3}{2}\alpha}_q|j|^{-6}.
\end{align*} With aid of \eqref{e:vv_q-bound}, \eqref{TqikN} and \eqref{matdvvq}, we have
\begin{align*}
\|\matd {\vv_q}\nabla\mathcal{T}_{k,i,j}\|_N\le&\|\matd {\vv_q}\mathcal{T}_{k,i,j}\|_{N+1}+\|\nabla\vv_q\cdot\nabla\mathcal{T}_{k,i,j}\|_{N}\lesssim\delta_{q+1} \delta_{q}^{1/2}\lambda_q \ell_q^{-N-1-\frac{3}{4}\alpha}.
\end{align*}
This inequality helps us obtain that
\begin{align*}
\|E_{k,i,j}\|_{N}&\le\delta^{-1/2}_{q+1}\lambda^{-1}_{q+1}\ell^{-\alpha/2}_q|j|^{-1}|b_{k,j}|(\|\eta_i(t)(\matd {\vv_q}\nabla\mathcal{T}_{k,i,j}) \nabla\Phi^{\TT}_i\|_N\\
&\qquad+\|\eta_i(t)\nabla\mathcal{T}_{k,i,j}( \matd {\vv_q}\nabla\Phi^{\TT}_i)\|_N+\|\nabla\mathcal{T}_{k,i,j}\nabla\Phi^{\TT}_i\eta'_i(t)\|_N)\\
&\lesssim\delta^{1/2}_{q+1} \delta_{q}^{1/2}\lambda_q \lambda^{-1}_{q+1}\ell_q^{-N-1-\frac{3}{2}\alpha}|j|^{-7}.
\end{align*}
By  \eqref{e:nabla-phi-i-CN}, \eqref{e:nabla-phi-i-matd}, \eqref{avkN} and \eqref{Davk}, we obtain from the definition of $F_{k,i,j}$ in \eqref{EF} that for $0\le N\le N_0$,
\begin{align*}
\big\|\matd {\vv_q}F_{k,i,j}\big\|_N\le&\delta^{1/2}_{q+1}\lambda^{-1}_{q+1}\ell^{\alpha/2}_q|j|^{-1}|b_{k,j}|(\|\eta_i(t)\matd {\vv_q}\big(\nabla a_{v,k}(R_v))\nabla\Phi^{\TT}_i\|_N\\
&+\|\nabla \big(a_{v,k}(R_v)\big)\nabla\Phi^{\TT}_i\eta'_i(t)\|_N\\
&+\|\eta_i(t)\nabla \big(a_{v,k}(R_v)\big)( \matd {\vv_q}\nabla\Phi^{\TT}_i)\|_N)\\
\lesssim& \delta^{1/2}_{q+1}\delta_{q}^{1/2}\lambda_q\lambda^{-1}_{q+1} \ell_q^{-N-1-\frac{3\alpha}{2}}|j|^{-7}.
\end{align*}
\end{proof}
\vskip 3mm
\noindent{\textbf{Step 4:\,\,Estimates for the perturbation}}\\
Now we can bound the perturbation $(w_{q+1}, d_{q+1})$ based on the estimates of  these coefficients.
\begin{prop}[Estimates for $w_{q+1}$ and $d_{q+1}$]\label{estimate-wq+1}There exists a universal constant $M$ satisfying
\begin{align}
&\|w^{(p)}_{q+1}\|_0+\tfrac{1}{\lambda_{q+1}}\|w^{(p)}_{q+1}\|_1\le \tfrac{1}{8}\delta^{1/2}_{q+1},\label{estimate-wp}\\
&\|w^{(c)}_{q+1}\|_0+\tfrac{1}{\lambda_{q+1}}\|w^{(c)}_{q+1}\|_1\le\lambda^{-1}_{q+1}\delta^{1/2}_{q+1}\ell^{-1}_q,\label{estimate-wc}\\
&\|w_{q+1}\|_0+\tfrac{1}{\lambda_{q+1}}\|w_{q+1}\|_1\le \tfrac{1}{4}\delta^{1/2}_{q+1},\label{estimate-w}\\
&\|d_{q+1}\|_0+\tfrac{1}{\lambda_{q+1}}\|d_{q+1}\|_1\le\tfrac{1}{8}\delta^{1/2}_{q+1}.  \label{estimate-d}
\end{align}
Moreover, by interpolation, we obtain from the above estimates that
\begin{equation}\label{estimate-alpha}
\|w^{(p)}_{q+1}\|_{\alpha}\le \lambda^{\alpha}_{q+1}\delta^{1/2}_{q+1},\quad\|w^{(c)}_{q+1}\|_{\alpha}\le \lambda^{-1+\alpha}_{q+1}\ell^{-1}_q\delta^{1/2}_{q+1},\quad \|d_{q+1}\|_{\alpha}\le \lambda^{\alpha}_{q+1}\delta^{1/2}_{q+1}.
\end{equation}
\end{prop}
\begin{proof}By \eqref{AN} and \eqref{BN}, we deduce from definition of $w^{(p)}_{q+1}$ in \eqref{wp_q+1} that
\begin{align}
\|w^{(p)}_{q+1}\|_0&\le\max_{i}\sum_{k\in\Lambda_{\theta}}\sum_{j\in I^{\theta}_k}\|A_{k,i,j}\|_0+\max_{i}\sum_{k\in\Lambda_{v}}\sum_{j\in I^v_k}\|B_{k,i,j}\|_0\nonumber\\
&\lesssim\delta^{1/2}_{q+1}\ell^{\frac{\alpha}{4}}_q\sum_{j\in I^{\theta}_k\cup I^{v}_k}|j|^{-6}\le\frac{1}{16}\delta^{1/2}_{q+1}.\label{wpq+10}
\end{align}
With respect to $C^1$ norm, we obtain
\begin{align*}
\|w^{(p)}_{q+1}\|_1\le &\max_i\sum_{k\in\Lambda_{\theta}}\sum_{j\in I^{\theta}_k}(\|A_{k,i,j}\|_1+2\pi\lambda_{q+1}|j|\|\nabla\Phi_i\|_0\|A_{k,i,j}\|_0)\\
&+\max_i\sum_{k\in\Lambda_{v}}\sum_{j\in I^v_k}(\|B_{k,i,j}\|_1+2\pi\lambda_{q+1}|j|\|\nabla\Phi_i\|_0\|B_{k,i,j}\|_0)\\
\lesssim&\delta^{1/2}_{q+1}\ell^{\frac{\alpha}{4}}_q\Big(\sum_{j\in I^{\theta}_k\cup I^{v}_k}\ell^{-1}_q|j|^{-6}+2\pi\lambda_{q+1}\sum_{j\in I^{\theta}_k\cup I^v_k}|j|^{-5}
\Big).
\end{align*}
Noting  that $\sum_{j\in I^{\theta}_k\cup I^v_k}|j|^{-5}$ is finite, we have from this inequality that
\begin{equation}\label{wpq+11}
\begin{aligned}
\lambda^{-1}_{q+1}\|w^{(p)}_{q+1}\|_1\lesssim \delta^{1/2}_{q+1}\ell^{\frac{\alpha}{4}}_q\le\frac{1}{16}\delta^{1/2}_{q+1}.
\end{aligned}
\end{equation}
Collecting \eqref{wpq+10} and \eqref{wpq+11} together yields \eqref{estimate-wp}. Plugging the estimates \eqref{EN} and \eqref{FN} into \eqref{wc} yields that
\begin{align*}
\|w^{(c)}_{q+1}\|_0\le &{(2\pi\lambda_{q+1})}^{-1}(\max_i\sum_{k\in\Lambda_{\theta}}\sum_{j\in I^{\theta}_k}\|E_{k,i,j}\|_0+\max_i\sum_{k\in\Lambda_{v}}\sum_{j\in I^v_k}\|F_{k,i,j}\|_0)\\
\lesssim&\delta^{1/2}_{q+1}{\lambda^{-1}_{q+1}}\ell^{-1+\frac{\alpha}{4}}_q\sum_{j\in I^{\theta}_k\cup I^{v}_k}|j|^{-7}\le \frac{1}{2} \delta^{1/2}_{q+1}{\lambda^{-1}_{q+1}}\ell^{-1}_q,
\end{align*}
where the constant can be absorbed by $\ell^{\alpha}_q$ for large enough $a$. Furthermore,
\begin{align*}
\|w^{(c)}_{q+1}\|_1\le &\max_i\sum_{k\in\Lambda_{\theta}}\sum_{j\in I^{\theta}_k}(\|E_{k,i,j}\|_1+2\pi\lambda_{q+1}|j|\|\nabla\Phi_i\|_0\|E_{k,i,j}\|_0)\\
&\qquad+\max_i\sum_{k\in\Lambda_{v}}\sum_{j\in I^v_k}(\|F_{k,i,j}\|_1+2\pi\lambda_{q+1}|j|\|\nabla\Phi_i\|_0\|F_{k,i,j}\|_0)\\
\lesssim&\delta^{1/2}_{q+1}{\lambda^{-1}_{q+1}}\ell^{\frac{\alpha}{4}}_q\sum_{j\in I^{\theta}_k\cup I^{v}_k}(\ell^{-2}_q|j|^{-7}+\ell^{-1}_q\lambda_{q+1}|j|^{-6})\\
\le&  \frac{1}{2}\delta^{1/2}_{q+1}\ell^{-1}_q.
\end{align*}
The above two inequalities show \eqref{estimate-wc}. By the definition of \eqref{bkj} and \eqref{d}, it is easy to verify  that
\begin{align*}
\|d_{q+1}\|_0\le& \delta^{1/2}_{q+1}\ell^{\alpha/2}_q\max_i\sum_{k\in\Lambda_{\theta}}\sum_{j\in I_k}|b_{k,j}| \lesssim\delta^{1/2}_{q+1}\ell^{\alpha/2}_q\sum_{j\in I^{\theta}_k}|j|^{-6}\le\frac{1}{16}\delta^{1/2}_{q+1},\\
\|d_{q+1}\|_1\le& \delta^{1/2}_{q+1}\ell^{\alpha/2}_q\max_i\sum_{k\in\Lambda_{\theta}}\sum_{j\in I^{\theta}_k}2\pi\lambda_{q+1}\|\nabla\Phi_i\|_0|b_{k,j}||j|\\
\lesssim&\delta^{1/2}_{q+1}\ell^{\alpha/2}_q\lambda_{q+1}\le\frac{1}{16}\delta^{1/2}_{q+1}\lambda_{q+1},
\end{align*}
which imply \eqref{estimate-d}. Hence we complete the proof of Proposition \ref{estimate-wq+1}.
\end{proof}
\subsection{Estimates for Reynolds and temperature stresses }
Define $(v_{q+1},\theta_{q+1})=(\vv_q+w_{q+1}, \tb_{q}+d_{q+1})$, one verifies that $(v_{q+1}, \theta_{q+1}, p_{q+1}, \RR_{q+1}, \MM_{q+1})$ satisfies \eqref{e:subsol-B} with replacing $q$ by $q+1$, where
\begin{align*}
\RR_{q+1}
=&\underbrace{\mathcal{R}( \del_t w_{q+1} + \vv_q \cdot \nabla w_{q+1})}_{\Rtransport}+\underbrace{\mathcal R (w_{q+1}\cdot\nabla \vv_q -d_{q+1}e_3)}_{\Rnash}\\
&+\underbrace{\mathcal{R}\big( \div (w_{q+1}\otimes w_{q+1})+\div\RRR_q\big)}_{\Rosc},
\end{align*}
\begin{align*}
\MM_{q+1}
=&\underbrace{\mathcal R_{\vex}(\del_t d_{q+1} + \vv_q \cdot \nabla d_{q+1} )}_{\Mtrans}+\underbrace{ \mathcal R_{\vex}(w_{q+1}\cdot\nabla\tb_q)}_{\Mnash}+\underbrace{\mathcal R_{\vex}(\div \MMM_q + \div(w_{q+1} d_{q+1}) )}_{\Mosc},
\end{align*}
and $p_{q+1}=\pp_q$.

\begin{prop}[Estimates for $\Mtrans$ and $\Rtransport$]\label{proptrans}
\begin{equation}\label{Trans}
\begin{aligned}
\|\Mtrans\|_{\alpha}\lesssim \frac{\delta_{q+1}^{1/2}\delta_q^{1/2} \lambda_q}{\lambda_{q+1}^{1-3\alpha}},\qquad
\|\Rtransport\|_{\alpha}\lesssim \frac{\delta_{q+1}^{1/2}\delta_q^{1/2} \lambda_q}{\lambda_{q+1}^{1-3\alpha}}.
\end{aligned}
\end{equation}
\end{prop}
\begin{proof}Notice that $\matd {\vv_q}e^{2\pi\ii\lambda_{q+1}j\cdot\Phi_i}=0$ due to $\matd {\vv_q}\Phi_i=0$,  from the definition of~$d_{q+1}$ in~\eqref{d}, we can express $\Mtrans$ as follows:
\begin{align*}
\Mtrans=&\mathcal R_{\vex}\Big(\sum_{i=0}^{i_{\max}-1}\sum_{k\in\Lambda_{\theta}}\sum_{j\in I^{\theta}_k}\delta^{{1/2}}_{q+1}\ell^{\alpha/2}_qb_{k,j}\matd {\vv_q}\big(\eta_i(t)e^{2\pi\ii\lambda_{q+1}j\cdot\Phi_i}\big) \Big)\\
=&\sum_{i=0}^{i_{\max}-1}\sum_{k\in\Lambda_{\theta}}\sum_{j\in I^{\theta}_k}\mathcal R_{\vex}\Big(\delta^{{1/2}}_{q+1}\ell^{\alpha/2}_qb_{k,j}\eta'_i(t)e^{2\pi\ii\lambda_{q+1}j\cdot\Phi_i} \Big).
\end{align*}
By Lemma \ref{l:non-stationary-phase}, we have
\begin{align*}
\|\Mtrans\|_{\alpha}\lesssim &\delta^{1/2}_{q+1}\ell^{\alpha/2}_q\lambda^{\alpha-1}_{q+1}\max_i\sum_{j\in I^{\theta}_k}\|\eta'_i(t)b_{k,j}\|_0\\
&+\delta^{1/2}_{q+1}\ell^{\alpha/2}_q\lambda^{\alpha-N_0}_{q+1}\max_i\sum_{j\in I^{\theta}_k}(\|\eta'_i(t)b_{k,j}\|_{N_0+\alpha}+\|\eta'_i(t)b_{k,j}\|_{0}\|\Phi_i\|_{N_0+\alpha}).
\end{align*}
Because $\|\eta'_i(t)\|_{0}\lesssim \tau^{-1}_q$, by \eqref{lambdaN}, \eqref{e:nabla-phi-i-CN} and \eqref{bkj}, we infer from this inequality  that
\begin{align*}
\|\Mtrans\|_{\alpha}\lesssim&\delta^{1/2}_{q+1}\ell^{\alpha/2}_q\tau^{-1}_q(\lambda^{\alpha-1}_{q+1}+\lambda^{\alpha-N_0}_{q+1}\ell^{-N_0-\alpha+1}_q)\sum_{j\in I^{\theta}_k}|j|^{-6}\\
\lesssim&\frac{\delta^{1/2}_{q+1}\delta^{1/2}_q\lambda_q}{\lambda^{1-\alpha}_{q+1}\ell^{{3\alpha}/{2}}_q}+\frac{\delta^{1/2}_{q+1}\delta^{1/2}_q\lambda_q}{\lambda^{N_0-\alpha}_{q+1}\ell^{N_0-1+\frac{5\alpha}{2}}_q}
\lesssim\frac{\delta^{1/2}_{q+1}\delta^{1/2}_q\lambda_q}{\lambda^{1-3\alpha}_{q+1}}+\frac{\delta^{1/2}_{q+1}\delta^{1/2}_q\lambda_q}{\lambda^{N_0-\alpha}_{q+1}\ell^{N_0+\alpha}_q},
\end{align*}
where we have used $\ell^{-2\alpha}_q\le\lambda^{2\alpha}_{q+1}$ in the last inequality.
Combining with \eqref{lambdaN}, we complete the estimate for $\Mtrans$ in \eqref{Trans}.

Next, we turn to estimate $\Rtransport$. By decomposition $w_{q+1}=w^{(p)}_{q+1}+w^{(c)}_{q+1}$, we bound $\|\Rtransport\|_{\alpha}$ by two parts:
\begin{equation}\label{Rtrans}
\begin{aligned}
\|\Rtransport\|_{\alpha}=\|\mathcal{R}(\matd {\vv_q} w_{q+1})\|_{\alpha}\le \|\mathcal{R}(\matd {\vv_q} w^{(p)}_{q+1})\|_{\alpha}+\|\mathcal{R}(\matd {\vv_q} w^{(c)}_{q+1})\|_{\alpha}.
\end{aligned}
\end{equation}
By the definition of $w^{(p)}_{q+1}$ in \eqref{wp_q+1} and $\matd {\vv_q}e^{2\pi\ii\lambda_{q+1}j\cdot\Phi_i}=0$, we rewrite $\mathcal{R}(\matd {\vv_q} w^{(p)}_{q+1})$ as follows:
\begin{align*}
\mathcal{R}(\matd {\vv_q} w^{(p)}_{q+1})=\sum_{i=0}^{i_{\max}-1}\mathcal{R}\Big(\sum_{k\in\Lambda_{\theta}}\sum_{j\in I^{\theta}_k}(\matd {\vv_q}A_{k,i,j})ke^{2\pi\ii\lambda_{q+1}j\cdot\Phi_i}+\sum_{k\in\Lambda_{v}}\sum_{j\in I^v_k}(\matd {\vv_q}B_{k,i,j})\bar{k}e^{2\pi\ii\lambda_{q+1}j\cdot\Phi_i}\Big).
\end{align*}
Applying Lemma \ref{l:non-stationary-phase} to the above equality, we obtain that
\begin{align*}
\|\mathcal{R}(\matd {\vv_q} w^{(p)}_{q+1})\|_{\alpha}\lesssim&\lambda^{\alpha-1}_{q+1}\max_i\big(\sum_{j\in I^{\theta}_k}\|\matd {\vv_q}A_{k,i,j}\|_0+\sum_{j\in I^v_k}\|\matd {\vv_q}B_{k,i,j}\|_0\big)\\
&+\lambda^{\alpha-N_0}_{q+1}\max_i\big(\sum_{j\in I^{\theta}_k}\|\matd {\vv_q}A_{k,i,j}\|_{N_0+\alpha}+\sum_{j\in I^{v}_k}\|\matd {\vv_q}B_{k,i,j}\|_{N_0+\alpha}\big)\\
&+\lambda^{\alpha-N_0}_{q+1}\max_i\big(\sum_{j\in I^{\theta}_k}\|\matd {\vv_q}A_{k,i,j}\|_{0}+\sum_{j\in I^{v}_k}\|\matd {\vv_q}B_{k,i,j}\|_{0}\big)\|\Phi_i\|_{N_0+\alpha}.
\end{align*}
The estimates \eqref{e:nabla-phi-i-CN}, \eqref{DAN} and \eqref{DBN} give that
\begin{equation}\label{Rtransp}
 \begin{aligned}
 \|\mathcal{R}(\matd {\vv_q} w^{(p)}_{q+1})\|_{\alpha}\lesssim\frac{\delta^{1/2}_{q+1}\delta^{1/2}_{q}\lambda_q\ell^{-\frac{3}{2}\alpha}_q}{\lambda^{1-\alpha}_{q+1}}
 +\frac{\delta^{1/2}_{q+1}\delta^{1/2}_{q}\lambda_q}{\lambda^{N_0-\alpha}_{q+1}\ell^{N_0+3\alpha}_q}\lesssim \frac{\delta^{1/2}_{q+1}\delta^{1/2}_{q}\lambda_q}{\lambda^{1-3\alpha}_{q+1}}.
 \end{aligned}
 \end{equation}
Similarly, from the definition of $w^{(c)}_{q+1}$ in \eqref{wc}, one has
\begin{align*}
\mathcal{R}(\matd {\vv_q} w^{(c)}_{q+1})=\sum_{i=0}^{i_{\max}-1}\mathcal{R}\Big(\sum_{k\in\Lambda_{\theta}}\sum_{j\in I^{\theta}_k}(\matd {\vv_q}E_{k,i,j})ke^{2\pi\ii\lambda_{q+1}j\cdot\Phi_i}+\sum_{k\in\Lambda_{v}}\sum_{j\in I^v_k}(\matd {\vv_q}F_{k,i,j})\bar{k}e^{2\pi\ii\lambda_{q+1}j\cdot\Phi_i}\Big).
\end{align*}
With the aid of Lemma \ref{l:non-stationary-phase}, by  estimates \eqref{e:nabla-phi-i-CN}, \eqref{DEN} and \eqref{DFN}, we deduce that
\begin{align}
\|\mathcal{R}(\matd {\vv_q} w^{(c)}_{q+1})\|_{\alpha}\lesssim&\lambda^{\alpha-1}_{q+1}\max_i\big(\sum_{j\in I^{\theta}_k}\|\matd {\vv_q}E_{k,i,j}\|_0+\sum_{j\in I^{v}_k}\|\matd {\vv_q}F_{k,i,j}\|_0\big)\nonumber\\
&+\lambda^{\alpha-N_0}_{q+1}\max_i\big(\sum_{j\in I^{\theta}_k}\|\matd {\vv_q}E_{k,i,j}\|_{N_0+\alpha}+\sum_{j\in I^{v}_k}\|\matd {\vv_q}F_{k,i,j}\|_{N_0+\alpha}\big)\nonumber\\
&+\lambda^{\alpha-N_0}_{q+1}\max_i\big(\sum_{j\in I^{\theta}_k}\|\matd {\vv_q}E_{k,i,j}\|_{0}+\sum_{j\in I^{v}_k}\|\matd {\vv_q}F_{k,i,j}\|_{0}\big)\|\Phi_i\|_{N_0+\alpha}\nonumber\\
\lesssim&\frac{\delta^{1/2}_{q+1}\delta^{1/2}_{q}\lambda_q\ell^{-1-\frac{3}{2}\alpha}_q}{\lambda^{2-\alpha}_{q+1}}+\frac{\delta^{1/2}_{q+1}\delta^{1/2}_{q}\lambda_q}{\lambda^{N_0+1-\alpha}_{q+1}\ell^{N_0+1+3\alpha}_q}
\lesssim\frac{\delta^{1/2}_{q+1}\delta^{1/2}_{q}\lambda_q}{\lambda^{1-3\alpha}_{q+1}}.\label{Rtransc}
\end{align}
Substituting inequalities \eqref{Rtransp} and \eqref{Rtransc} into \eqref{Rtrans}, we obtain the second estimate in~\eqref{Trans}.
\end{proof}
\begin{prop}[Estimates for $\Mnash$ and $\Rnash$]
\begin{equation}\label{nash}
\begin{aligned}
\|\Mnash\|_{\alpha}\lesssim\frac{\delta^{1/2}_{q+1}\delta^{1/2}_{q}\lambda_q}{\lambda^{1-\alpha}_{q+1}},\qquad
\|\Rnash\|_{\alpha}\lesssim\frac{\delta^{1/2}_{q+1}\delta^{1/2}_{q}\lambda_q}{\lambda^{1-3\alpha}_{q+1}}.
\end{aligned}
\end{equation}
\end{prop}
\begin{proof}Noting that
\begin{align*}
\Mnash=\mathcal{R}_{\vex}(w^{(p)}_{q+1}\cdot\nabla \tb_q)+\mathcal{R}_{\vex}(w^{(c)}_{q+1}\cdot\nabla \tb_q),
\end{align*}
we estimate $\mathcal{R}_{\vex}(w^{(p)}_{q+1}\cdot\nabla \tb_q)$ and $\mathcal{R}_{\vex}(w^{(c)}_{q+1}\cdot\nabla \tb_q)$ separately. By the definition of $w^{(p)}_{q+1}$ in \eqref{wp_q+1}, we have
\begin{align}
\|\mathcal{R}_{\vex}(w^{(p)}_{q+1}\cdot\nabla \tb_q)\|_{\alpha}\le& \max_i\sum_{k\in \Lambda_{\theta}}\sum_{j\in I^{\theta}_k} \|\mathcal{R}_{\vex}(A_{k,i,j}ke^{\ii\lambda_{q+1}j\cdot\Phi_i}\cdot\nabla\tb_q )\|_{\alpha}\nonumber\\
&+\max_i\sum_{k\in \Lambda_{v}}\sum_{j\in I^v_k} \|\mathcal{R}_{\vex}(B_{k,i,j}\bar{k} e^{\ii\lambda_{q+1}j\cdot\Phi_i}\cdot\nabla\tb_q)\|_{\alpha}.\label{Nash}
\end{align}
Employing Lemma \ref{l:non-stationary-phase} to \eqref{Nash} yields that
\begin{align}
&\|\mathcal{R}_{\vex}(w^{(p)}_{q+1}\cdot\nabla \tb_q)\|_{\alpha}\nonumber\\
\lesssim  & \lambda^{\alpha-1}_{q+1}\max_i\big(\sum_{k\in \Lambda_{\theta}}\sum_{j\in I^{\theta}_k}{\|A_{k,i,j}\|_0\|\nabla\tb_q\|_0+\sum_{k\in \Lambda_{\theta}}\sum_{j\in I^{v}_k}\|B_{k,i,j}\|_0\|\nabla\tb_q\|_0}\big)\nonumber\\
&+ \lambda^{\alpha-N_0}_{q+1}\max_i\sum_{k\in \Lambda_{\theta}}\sum_{j\in I^{\theta}_k}\big(\|A_{k,i,j}k \cdot\nabla\tb_q\|_{N_0+\alpha}+\|A_{k,i,j}\|_0\|\nabla\tb_q\|_0\|\Phi_i\|_{N_0+\alpha}\big) \nonumber\\
&+\lambda^{\alpha-N_0}_{q+1}\max_i\sum_{k\in \Lambda_{\theta}}\sum_{j\in I^{v}_k}\big(\|B_{k,i,j}k \cdot\nabla\tb_q\|_{N_0+\alpha}+\|B_{k,i,j}\|_0\|\nabla\tb_q\|_0\|\Phi_i\|_{N_0+\alpha}\big).\label{Nashp}
\end{align}
With aid of \eqref{e:vv_q-bound}, \eqref{e:nabla-phi-i-CN}, \eqref{AN} and  \eqref{BN}, we deduce
\begin{align*}
&\|\mathcal{R}_{\vex}(w^{(p)}_{q+1}\cdot\nabla \tb_q)\|_{\alpha}
\lesssim \Big(\frac{\delta^{1/2}_{q+1}\delta_{q}^{1 / 2} \lambda_{q}\ell^{\alpha/4}_q }{\lambda^{1-\alpha}_{q+1}}+
\frac{\delta^{1/2}_{q+1}\delta_{q}^{1 / 2}  \lambda_q}{\lambda^{N_0-\alpha}_{q+1}\ell^{N_0+\alpha}_q}\Big)
\sum_{j\in I^{\theta}_k\cup I^{v}_k}|j|^{-6}\lesssim \frac{\delta^{1/2}_{q+1}\delta_{q}^{1 / 2} \lambda_{q} }{\lambda^{1-\alpha}_{q+1}}.
\end{align*}
Using the definition of $w^{(c)}_{q+1}$ in \eqref{wc}, by Lemma \ref{l:non-stationary-phase}, one gets
\begin{align*}
&\|\mathcal{R}_{\vex}(w^{(c)}_{q+1}\cdot\nabla \tb_q)\|_{\alpha}\\
\lesssim  & \lambda^{\alpha-1}_{q+1}\max_i\big(\sum_{k\in \Lambda_{\theta}}\sum_{j\in I^{\theta}_k}{\|E_{k,i,j}\|_0\|\nabla\tb_q\|_0+\sum_{k\in \Lambda_{\theta}}\sum_{j\in I^{v}_k}\|F_{k,i,j}\|_0\|\nabla\tb_q\|_0}\big)\\
&+ \lambda^{\alpha-N_0}_{q+1}\max_i\sum_{k\in \Lambda_{\theta}}\sum_{j\in I^{\theta}_k}\big(\|E_{k,i,j}k \cdot\nabla\tb_q\|_{N_0+\alpha}+\|E_{k,i,j}\|_0\|\nabla\tb_q\|_0\|\Phi_i\|_{N_0+\alpha}\big)\\
&+ \lambda^{\alpha-N_0}_{q+1}\max_i\sum_{k\in \Lambda_{\theta}}\sum_{j\in I^v_k}\big(\|F_{k,i,j}k \cdot\nabla\tb_q\|_{N_0+\alpha}+\|F_{k,i,j}\|_0\|\nabla\tb_q\|_0\|\Phi_i\|_{N_0+\alpha}\big).
\end{align*}
Applying  \eqref{e:vv_q-bound}, \eqref{e:nabla-phi-i-CN}, \eqref{EN} and \eqref{FN} to this inequality, we have
\begin{align}
\|\mathcal{R}_{\vex}(w^{(c)}_{q+1}\cdot\nabla \tb_q)\|_{\alpha}
\lesssim & \Big(\frac{\delta^{1/2}_{q+1}\delta_{q}^{1 / 2} \lambda_{q} \ell^{-1+\frac{1}{4}\alpha}_q}{\lambda^{2-\alpha}_{q+1}}+
\frac{\delta^{1/2}_{q+1}\delta_{q}^{1 / 2}  \lambda_q}{\lambda^{N_0+1-\alpha}_{q+1}\ell^{N_0+1}_q}\Big)
\sum_{j\in I^{\theta}_k \cup I^v_k}|j|^{-7}\nonumber\\
\lesssim&\frac{\delta^{1/2}_{q+1}\delta_{q}^{1 / 2} \lambda_{q} }{\lambda^{1-\alpha}_{q+1}},\label{Nashc}
\end{align}
where we have used the fact that $\ell^{-1}_q<\lambda_{q+1}$. Plugging \eqref{Nashp} and \eqref{Nashc} into \eqref{Nash}, we obtain  the estimate for $\Mnash$ in \eqref{nash}. According to the definition of $\Rnash$, we have
\begin{equation}\label{Rnash}
\begin{aligned}
\|\Rnash\|_{\alpha}\le\|\mathcal R (w_{q+1}\cdot\nabla \vv_q )\|_{\alpha}+\|\mathcal R (d_{q+1}e_3)\|_{\alpha}
\end{aligned}
\end{equation}
Thanks to  the same upper bound of $(\vv_q, \tb_q)$ in \eqref{e:vv_q-bound}, we immediately infer from the estimate for $\mathcal{R}_{\vex}(w_{q+1}\cdot\nabla \tb_q)$ that
\begin{equation}\label{Rnash1}
\|\mathcal{R}(w_{q+1}\cdot\nabla \vv_q)\|_{\alpha}\lesssim \frac{\delta^{1/2}_{q+1}\delta_{q}^{1 / 2} \lambda_{q} }{\lambda^{1-\alpha}_{q+1}}.
\end{equation}
From  \eqref{d}, we have
\[\mathcal{R}(d_{q+1}e_3)=\sum_{i=0}^{i_{\max}-1}\sum_{k\in\Lambda_{\theta}}\sum_{j\in I^{\theta}_k}\mathcal{R}\big(\delta^{{1/2}}_{q+1}\ell^{\alpha/2}_q\eta_i(t)b_{k,j}e^{2\pi\ii\lambda_{q+1}j\cdot\Phi_i}e_3\big).\]
By Lemma  \ref{l:non-stationary-phase}, \eqref{e:nabla-phi-i-CN} and \eqref{bkj}, we get by \eqref{lambdaN} that
\begin{align}
\|\mathcal{R}(d_{q+1}e_3)\|_{\alpha}
\lesssim  & \lambda^{\alpha-1}_{q+1}\sum_{j\in I^{\theta}_k}\delta^{{1/2}}_{q+1}\ell^{\alpha/2}_q|b_{k,j}|+ \lambda^{\alpha-N_0}_{q+1}\sum_{j\in I^{\theta}_k}\delta^{{1/2}}_{q+1}\ell^{\alpha/2}_q|b_{k,j}|\|\Phi_i\|_{N_0+\alpha}\nonumber\\
\lesssim  &\lambda^{\alpha-1}_{q+1}\sum_{j\in I^{\theta}_k}\delta^{{1/2}}_{q+1}\ell^{\alpha/2}_q|j|^{-6}+ \lambda^{\alpha-N_0}_{q+1}\sum_{j\in I^{\theta}_k}\delta^{{1/2}}_{q+1}\ell^{\alpha/2}_q|j|^{-6}\ell^{-N_0+1-\alpha}_q\nonumber\\
\lesssim &\frac{\delta^{1/2}_{q+1}}{\lambda^{1-\alpha}_{q+1}}+\frac{\delta^{1/2}_{q+1}}{\lambda^{N_0-\alpha}_{q+1}\ell^{N_0+\alpha}_q}\lesssim\frac{\delta^{1/2}_{q+1}}{\lambda^{1-\alpha}_{q+1}}.\label{Rnash2}
\end{align}
Putting \eqref{Rnash1} and \eqref{Rnash2} into \eqref{Rnash}, we prove \eqref{nash}.
\end{proof}

\begin{prop}[Estimates for $\Mosc$]\label{propMosc}
\begin{align*}
\|\Mosc\|_{\alpha}\le\frac{\delta^{1/2}_{q+1}\delta^{1/2}_q\lambda_q}{\lambda^{1-3\alpha}_{q+1}}.
\end{align*}
\end{prop}
\begin{proof}Owning to $w_{q+1}=w^{(p)}_{q+1}+w^{(c)}_{q+1}$, we split  $\Mosc$ into the following two parts:
\begin{equation}\label{Mosc}
\begin{aligned}
\|\Mosc\|_{\alpha}\le \|\mathcal{R}_{\vex}\div (\MMM_q+w^{(p)}_{q+1}d_{q+1})\|_{\alpha}+\|\mathcal{R}_{\vex}\div (w^{(c)}_{q+1}d_{q+1})\|_{\alpha}.
\end{aligned}
\end{equation}
Taking advantage of  \eqref{estimate-w}--\eqref{estimate-alpha}, one has
\begin{align}
\|\mathcal{R}_{\vex}\div (w^{(c)}_{q+1}d_{q+1})\|_{\alpha}\lesssim& \|w^{(c)}_{q+1}\|_{\alpha}\|d_{q+1}\|_0+\|w^{(c)}_{q+1}\|_{0}\|d_{q+1}\|_{\alpha}\nonumber\\
\lesssim&\lambda^{-1+\alpha}_{q+1}\ell^{-1}_q\delta_{q+1}\le \frac{\delta^{{1}/{2}}_{q+1}\delta^{{1}/{2}}_q\lambda_q}{\lambda^{1-3\alpha}_{q+1}}.\label{Mosc1}
\end{align}
where we have used the definition of $\ell_q$ in \eqref{e:ell}. Next, we aim to estimate the first term on the right-hand side of \eqref{Mosc}.
Note that $\eta_i$ possess disjoint supports and \eqref{WkWk'},
\begin{align*}
\sum_{k\in\Lambda_{\theta}}W^{\theta}_{(\lambda_{q+1}k)}(\Phi_i(x,t))\sum_{k\in\Lambda_{v}}W^v_{(\lambda_{q+1}k)}(\Phi_i(x,t))=0,
\end{align*}
we have by \eqref{wp_q+1} and \eqref{d} that
\begin{align}
&\mathcal{R}_{\vex}\div \big(\MMM_q+w^{(p)}_{q+1}d_{q+1}\big)\nonumber\\
=&\mathcal{R}_{\vex}\div \big(\MMM_q+\sum^{i_{\max-1}}_{i=0}\sum_{k\in\Lambda_{\theta}}\sum_{{j,j'\in I^{\theta}_k}}\mathcal{T}_{q,i,k}b_{k,j}b_{k,j'}\eta^2_i(t)\nabla\Phi^{-1}_ike^{2\pi \ii\lambda_{q+1}(j+j')\cdot\Phi_i}\big)\nonumber\\
=&\mathcal{R}_{\vex}\div \big(\MMM_q+\sum^{i_{\max-1}}_{i=0}\eta_i(t)\nabla\Phi^{-1}_i\sum_{k\in\Lambda_{\theta}}\mathcal{T}_{q,i,k}\eta_i(t)k\big)\nonumber\\
&+\mathcal{R}_{\vex}\div \Big(\sum^{i_{\max-1}}_{i=0}\sum_{k\in\Lambda_{\theta}}\sum_{\substack{j,j'\in I^{\theta}_k\\j+j'\neq 0}}\mathcal{T}_{q,i,k}\eta^2_i(t)\nabla\Phi^{-1}_ikb_{k,j}b_{k,j'}e^{2\pi\ii\lambda_{q+1}(j+j')\cdot\Phi_i}\Big),\label{mainMosc}
\end{align}
where we have used \eqref{WL2}. From the definition of $\mathcal{T}_{q,i,k}$  in \eqref{Tqik}, it is easy to verify that
\begin{equation}\label{Mosc=0}
\mathcal{R}_{\vex}\div \Big(\MMM_q+\sum^{i_{\max-1}}_{i=0}\eta_i(t)\nabla\Phi^{-1}_i\sum_{k\in\Lambda_{\theta}}\mathcal{T}_{q,i,k}\eta_i(t)k\Big)=\mathcal{R}_{\vex}\div \big(\MMM_q-\sum^{i_{\max-1}}_{i=0}\MMM_q\eta^2_i(t)\big)=0.
\end{equation}
Therefore, we only need to estimate $\|\mathcal{R}_{\vex}\div(\mathcal{T}_{q,i,k}\eta^2_i(t)\nabla\Phi^{-1}_ikb_{k,j}b_{k,j'}e^{2\pi\ii\lambda_{q+1}(j+j')\cdot\Phi_i})\|_{\alpha}$.   Note that $j+j'=l\bar{k}+m\bar{\bar{k}}(\forall\,l,m\neq 0)$, where $k\perp \bar{k}\perp \bar{\bar{k}}$, this shows $(j+j')\perp k$. Hence,
\begin{equation}\label{jk}
(\nabla\Phi^{-1}_ik)\cdot (j\cdot\nabla\Phi_i)=(\nabla\Phi^{-1}_ik)^{\TT}((j+j')^{\TT}\nabla\Phi_i)^{\TT}=k^{\TT}\nabla\Phi^{-{\TT}}\nabla\Phi^{\TT}_i(j+j')=0,
\end{equation}
from which it follows that
\begin{align*}
&\mathcal{R}_{\vex}\div(\mathcal{T}_{q,i,k}\eta^2_i(t)\nabla\Phi^{-1}_ikb_{k,j}b_{k,j'}e^{2\pi\ii\lambda_{q+1}(j+j')\cdot\Phi_i})\\
=&\mathcal{R}_{\vex}\big(\div\big(\mathcal{T}_{q,i,k}\eta^2_i(t)\nabla\Phi^{-1}_ikb_{k,j}b_{k,j'}\big)e^{2\pi\ii\lambda_{q+1}(j+j')\cdot\Phi_i}\big)\\
&+\mathcal{R}_{\vex}\big(2\pi\ii\lambda_{q+1}b_{k,j}b_{k,j'}\mathcal{T}_{q,i,k}\eta^2_i(t)(\nabla\Phi^{-1}_ik)\cdot ((j+j')\cdot\nabla\Phi_i)e^{2\pi\ii\lambda_{q+1}(j+j')\cdot\Phi_i}\big)\\
=&\mathcal{R}_{\vex}\big(\div\big(\mathcal{T}_{q,i,k}\eta^2_i(t)\nabla\Phi^{-1}_ikb_{k,j}b_{k,j'}\big)e^{2\pi\ii\lambda_{q+1}(j+j')\cdot\Phi_i}\big).
\end{align*}
Applying Lemma \ref{l:non-stationary-phase} to this  equality shows that
\begin{align*}
&\Big\|\mathcal{R}_{\vex}\div(\mathcal{T}_{q,i,k}\eta^2_i(t)\nabla\Phi^{-1}_ikb_{k,j}b_{k,j'}e^{2\pi\ii\lambda_{q+1}(j+j')\cdot\Phi_i})\Big\|_{\alpha}\\
\lesssim &\lambda^{-1+\alpha}_{q+1}\|\div\big(\mathcal{T}_{q,i,k}\eta^2_i(t)\nabla\Phi^{-1}_ikb_{k,j}b_{k,j'}\big)\|_0\\
&+\lambda^{-N_0+\alpha}_{q+1}(\|\div\big(\mathcal{T}_{q,i,k}\eta^2_i(t)\nabla\Phi^{-1}_ikb_{k,j}b_{k,j'}\big)\|_{N_0+\alpha}\\
&+\|\div\big(\mathcal{T}_{q,i,k}\eta^2_i(t)\nabla\Phi^{-1}_ikb_{k,j}b_{k,j'}\big)\|_0
\|\Phi_i\|_{N+\alpha}).
\end{align*}
By \eqref{e:nabla-phi-i-CN}, \eqref{bkj} and  \eqref{TqikN}, one can obtain that
\begin{align*}
&\Big\|\mathcal{R}_{\vex}\div(\mathcal{T}_{q,i,k}\eta^2_i(t)\nabla\Phi^{-1}_ikb_{k,j}b_{k,j'}e^{2\pi\ii\lambda_{q+1}j\cdot\Phi_i})\Big\|_{\alpha}\\
\lesssim &\frac{(\|\mathcal{T}_{q,i,k}\|_1+\|\mathcal{T}_{q,i,k}\|_0\|\eta^2_i(t)\nabla\Phi^{-1}_i\|_1)|b_{k,j}b_{k,j'}|}{\lambda^{1-\alpha}_{q+1}}\\
&+
\frac{(\|\mathcal{T}_{q,i,k}\|_1+\|\mathcal{T}_{q,i,k}\|_0\|\eta^2_i(t)\nabla\Phi^{-1}_i\|_1)|b_{k,j}b_{k,j'}|\|\Phi_i\|_{N_0+\alpha}}{\lambda^{N_0-\alpha}_{q+1}}\\
&+\frac{(\|\mathcal{T}_{q,i,k}\|_{N_0+1+\alpha}+\|\mathcal{T}_{q,i,k}\|_0\ \|\eta^2_i(t)\nabla\Phi^{-1}_i\|_{N_0+1+\alpha})|b_{k,j}b_{k,j'}|}{\lambda^{N_0-\alpha}_{q+1}}\\
\lesssim&
\big(\frac{\delta^{1/2}_{q+1}\delta^{1/2}_{q}\lambda^{1+3\alpha/2}_q\ell^{\alpha}_q}{\lambda^{1-\alpha}_{q+1}}+\frac{\delta^{1/2}_{q+1}\delta^{1/2}_{q}\lambda^{1+3\alpha/2}_q}{\lambda^{N_0-\alpha}_{q+1}\ell^{N_0}_q}\big)|jj'|^{-6}.
\end{align*}
Because $\sum_{j,j'\in I^{\theta}_k}|jj'|^{-6}$ is finite and $\lambda^{3\alpha/2}_{q}\ell^{\alpha}_q\le \lambda^{\alpha/2}_{q}\le  \lambda^{\alpha/2}_{q+1}$, we deduce from the above inequality that
\begin{align}
&\Big\|\sum^{i_{\max-1}}_{i=0}\sum_{k\in\Lambda_{\theta}}\sum_{\substack{j,j\in I^{\theta}_k\nonumber\\j+j'\neq 0}}\mathcal{R}_{\vex}\div(\mathcal{T}_{q,i,k}\eta^2_i(t)\nabla\Phi^{-1}_ikb_{k,j}b_{k,j'}e^{2\pi\ii\lambda_{q+1}(j+j')\cdot\Phi_i})\Big\|_{\alpha}\nonumber\\
\lesssim&\frac{\delta^{1/2}_{q+1}\delta^{1/2}_{q}\lambda^{1+3\alpha/2}_q\ell^{\alpha}_q}{\lambda^{1-\alpha}_{q+1}}+\frac{\delta^{1/2}_{q+1}\delta^{1/2}_{q}\lambda^{1+3\alpha/2}_q}{\lambda^{N_0-\alpha}_{q+1}\ell^{N_0}_q}
\lesssim\frac{\delta^{1/2}_{q+1}\delta^{1/2}_{q}\lambda_q}{\lambda^{1-3\alpha}_{q+1}}.\label{Moscmain1}
\end{align}
Plugging \eqref{Mosc=0} and \eqref{Moscmain1} into \eqref{mainMosc}, one obtains that
\begin{align*}
\|\mathcal{R}_{\vex}\div \big(\MMM_q+w^{(p)}_{q+1}d_{q+1}\big)\|_{\alpha}\lesssim\frac{\delta^{1/2}_{q+1}\delta^{1/2}_{q}\lambda_q}{\lambda^{1-3\alpha}_{q+1}}.
\end{align*}
Collecting this inequality together with \eqref{Mosc1} and \eqref{Mosc}, we complete  this proof of Proposition \ref{propMosc}.
\end{proof}
\begin{prop}[Estimates for $\Rosc$]\label{propRosc}
\begin{align*}
\|\Rosc\|_{\alpha}\le \frac{\delta^{1/2}_{q+1}\delta^{1/2}_{q}\lambda_{q}}{\lambda^{1-3\alpha}_{q+1}}.
\end{align*}
\end{prop}
\begin{proof}By the decomposition $w_{q+1}=w^{(p)}_{q+1}+w^{(c)}_{q+1}$, we divide $\Rosc$ into the four parts:
\begin{align*}
\Rosc=&\mathcal{R}\div(w^{(p)}_{q+1}\otimes w^{(p)}_{q+1}+\RRR_q)+\mathcal{R}\div(w^{(c)}_{q+1}\otimes w^{(c)}_{q+1})\\
&+\mathcal{R}\div(w^{(p)}_{q+1}\otimes w^{(c)}_{q+1})+\mathcal{R}\div(w^{(c)}_{q+1}\otimes w^{(p)}_{q+1})
\end{align*}
By  \eqref{estimate-wp}, \eqref{estimate-wc} and \eqref{estimate-alpha}, one immediately shows that
\begin{align}
&\|\Rosc-\mathcal{R}\div(w^{(p)}_{q+1}\otimes w^{(p)}_{q+1}+\RRR_q)\|_{\alpha}\nonumber\\
\lesssim&
\|w^{(c)}_{q+1}\|_{\alpha}\|w^{(c)}_{q+1}\|_{0}+\|w^{(p)}_{q+1}\|_{\alpha}\|w^{(c)}_{q+1}\|_{0}+\|w^{(p)}_{q+1}\|_{0}\|w^{(c)}_{q+1}\|_{\alpha}\nonumber\\
\lesssim&\frac{\delta_{q+1}\ell^{-1}_q}{\lambda^{1-\alpha}_{q+1}}=\frac{\delta^{1/2}_{q+1}\delta^{1/2}_{q}\lambda^{1+3\alpha/2}_q}{\lambda^{1-\alpha}_{q+1}}\le \frac{\delta^{1/2}_{q+1}\delta^{1/2}_{q}\lambda_q}{\lambda^{1-3\alpha}_{q+1}}.\label{Roscs}
\end{align}
Noting that  $\eta_i$ have disjoint supports and \eqref{WkWk'}, we deduce from the definition of $w^{(p)}_{q+1}$ in~\eqref{wp-def1} that
\begin{align*}
&\mathcal{R}\div(w^{(p)}_{q+1}\otimes w^{(p)}_{q+1}+\RRR_q)\\
=&\mathcal{R}\div\big(\delta^{-1}_{q+1}\sum_i\sum_{k\in\Lambda_{\theta}}\mathcal{T}^2_{q,i,k}\big(W^{\theta}_{(\lambda_{q+1}k)} (\Phi_i)\big)^2\eta^2_i(t)\nabla\Phi^{-1}_i(k\otimes k)\nabla\Phi^{-{\TT}}_i+\RRR_q\\
&+\delta_{q+1}\sum_i\sum_{k\in\Lambda_{v}}a^2_{v,k}(R_v)\big(W^v_{(\lambda_{q+1}k)} (\Phi_i)\big)^2\eta^2_i(t)\nabla\Phi^{-1}_i(\bar{k}\otimes \bar{k})\nabla\Phi^{-{\TT}}_i\big).
\end{align*}
By \eqref{WL2} and \eqref{W},  we  rewrite the above equality as follows:
\begin{align}
&\mathcal{R}\div(w^{(p)}_{q+1}\otimes w^{(p)}_{q+1}+\RRR_q)\nonumber\\
=&\delta^{-1}_{q+1}\ell^{-\alpha}_q\sum_i\sum_{k\in\Lambda_{\theta}}\mathcal{R}\div\big(\eta^2_i(t)\mathcal{T}^2_{q,i,k}\nabla\Phi^{-1}_i(k\otimes k)\nabla\Phi^{-{\TT}}_i+\RRR_q\big)\nonumber\\
&+\delta_{q+1}\ell^{\alpha}_q\sum_i\sum_{k\in\Lambda_{v}}\mathcal{R}\div\big(\eta^2_i(t)\nabla\Phi^{-1}_ia^2_{v,k}(R_v)(\bar{k}\otimes \bar{k})\nabla\Phi^{-{\TT}}_i\big)\nonumber\\
&+\delta^{-1}_{q+1}\ell^{-\alpha}_q\sum_i\sum_{k\in\Lambda_{\theta}}\sum_{\substack{j,j'\in I^{\theta}_k\\ j+j'\neq 0}}\mathcal{R}\div\big(\eta^2_i(t)\mathcal{T}^2_{q,i,k}\nabla\Phi^{-1}_i(k\otimes k)\nabla\Phi^{-{\TT}}_ib_{k,j}b_{k,j'}e^{2\pi\ii\lambda_{q+1}(j+j')\cdot\Phi_i}\big)\nonumber\\
&+\delta_{q+1}\ell^{\alpha}_q\sum_i\sum_{k\in\Lambda_{v}}\sum_{\substack{j,j'\in I^v_k\\ j+j'\neq 0}}\mathcal{R}\div\big(a^2_{v,k}(R_v)\nabla\Phi^{-1}_i(\bar{k}\otimes \bar{k})\nabla\Phi^{-{\TT}}_i b_{k,j}b_{k,j'}e^{2\pi\ii\lambda_{q+1}(j+j')\cdot\Phi_i}\big),\label{Roscdecom}
\end{align}
where we have used \eqref{WL2}.
According to Lemma \ref{first S}, the definition of $R_v$ in \eqref{Rv}, the first two terms on the right-hand side of \eqref{Roscdecom}
equal to
\begin{align}
&\mathcal{R}\div\big(\delta_{q+1}\ell^{\alpha}_q{\rm Id}\sum_i\eta^2_i(t) -\sum_i\eta^2_i(t) \nabla\Phi^{-1}_iR_v\nabla\Phi^{-{\TT}}_i\nonumber\\
&+\delta^{-1}_{q+1}\ell^{-\alpha}_q\sum_i\sum_{k\in\Lambda_{\theta}}\eta^2_i(t)\mathcal{T}^2_{q,i,k}\nabla\Phi^{-1}_i(k\otimes k)\nabla\Phi^{-{\TT}}_i+\RRR_q\big)=0.\label{Rosc=0}
\end{align}
Hence, it suffices to estimate the last two terms on the right-hand side of \eqref{Roscdecom}.
Using \eqref{jk}, for $(j+j')\perp k$, we have
\begin{align*}
&\div(\eta^2_i(t)\mathcal{T}^2_{q,i,k}\nabla\Phi^{-1}_i(k\otimes k)\nabla\Phi^{-{\TT}}_ib_{k,j}b_{k,j'}e^{2\pi\ii\lambda_{q+1}(j+j')\cdot\Phi_i})\\
=&\div(\eta^2_i(t)\mathcal{T}^2_{q,i,k}\nabla\Phi^{-1}_i(k\otimes k)\nabla\Phi^{-{\TT}}_ib_{k,j}b_{k,j'})e^{2\pi\ii\lambda_{q+1}(j+j')\cdot\Phi_i},
\end{align*}
and
\begin{align*}
&\div(a^2_{v,k}(R_v)\nabla\Phi^{-1}_i(\bar{k}\otimes \bar{k})\nabla\Phi^{-{\TT}}_i b_{k,j}b_{k,j'}e^{2\pi\ii\lambda_{q+1}(j+j')\cdot\Phi_i})\\
=&\div(a^2_{v,k}(R_v)\nabla\Phi^{-1}_i(\bar{k}\otimes \bar{k})\nabla\Phi^{-{\TT}}_i b_{k,j}b_{k,j'})e^{2\pi\ii\lambda_{q+1}(j+j')\cdot\Phi_i}.
\end{align*}
With the aid of Lemma \ref{l:non-stationary-phase}, we get from the above two equalities   that
\begin{align}
&\|\mathcal{R}\div(\eta^2_i(t)\mathcal{T}^2_{q,i,k}\nabla\Phi^{-1}_i(k\otimes k)\nabla\Phi^{-{\TT}}_ib_{k,j}b_{k,j'}e^{2\pi\ii\lambda_{q+1}(j+j')\cdot\Phi_i})\|_{\alpha}\nonumber\\
\lesssim&\frac{\|\div(\eta^2_i(t)\mathcal{T}^2_{q,i,k}\nabla\Phi^{-1}_i(k\otimes k)\nabla\Phi^{-{\TT}}_ib_{k,j}b_{k,j'})\|_0}{\lambda^{1-\alpha}_{q+1}}\nonumber\\
&+\frac{\|\div(\eta^2_i(t)\mathcal{T}^2_{q,i,k}\nabla\Phi^{-1}_i(k\otimes k)\nabla\Phi^{-{\TT}}_ib_{k,j}b_{k,j'})\|_{N_0+\alpha}}{\lambda^{N_0-\alpha}_{q+1}}\nonumber\\
&+\frac{\|\div(\eta^2_i(t)\mathcal{T}^2_{q,i,k}\nabla\Phi^{-1}_i(k\otimes k)\nabla\Phi^{-{\TT}}_ib_{k,j}b_{k,j'})\|_{0}\|\Phi_i\|_{N_0+\alpha}}{\lambda^{N_0-\alpha}_{q+1}}\nonumber\\
\lesssim&\Big(\frac{\delta^2_{q+1}\ell^{-1+\frac{5\alpha}{2}}_q}{\lambda^{1-\alpha}_{q+1}}+\frac{\delta^2_{q+1}\ell^{-N_0-1+\frac{3\alpha}{2}}_q}{\lambda^{N_0-\alpha}_{q+1}}\Big)|jj'|^{-6},\label{Rosc1}
\end{align}
where the last inequality holds by virtue of \eqref{e:nabla-phi-i-CN},  \eqref{bkj} and \eqref{TqikN}. In the same way,  we obtain
\begin{align}
&\|\mathcal{R}\div(a^2_{v,k}(R_v)\nabla\Phi^{-1}_i(\bar{k}\otimes \bar{k})\nabla\Phi^{-{\TT}}_i b_{k,j}b_{k,j'}e^{2\pi\ii\lambda_{q+1}(j+j')\cdot\Phi_i})\|_{\alpha}\nonumber\\
\lesssim&\frac{\|\div(a^2_{v,k}(R_v)\nabla\Phi^{-1}_i(\bar{k}\otimes \bar{k})\nabla\Phi^{-{\TT}}_i b_{k,j}b_{k,j'})\|_0}{\lambda^{1-\alpha}_{q+1}}\nonumber\\
&+\frac{\|\div(a^2_{v,k}(R_v)\nabla\Phi^{-1}_i(\bar{k}\otimes \bar{k})\nabla\Phi^{-{\TT}}_i b_{k,j}b_{k,j'})\|_{N_0+\alpha}}{\lambda^{N_0-\alpha}_{q+1}}\nonumber\\
&+\frac{\|\div(a^2_{v,k}(R_v)\nabla\Phi^{-1}_i(\bar{k}\otimes \bar{k})\nabla\Phi^{-{\TT}}_i b_{k,j}b_{k,j'}))\|_{0}\|\Phi_i\|_{N_0+\alpha}}{\lambda^{N_0-\alpha}_{q+1}}\nonumber\\
\lesssim&\Big(\frac{\ell^{-1+\frac{\alpha}{8}}_q}{\lambda^{1-\alpha}_{q+1}}+\frac{\ell^{-N_0-1-\frac{7\alpha}{8}}_q}{\lambda^{N_0-\alpha}_{q+1}}\Big)|jj'|^{-6},\label{Rosc2}
\end{align}
by using the estimates \eqref{e:nabla-phi-i-CN}, \eqref{bkj} and  \eqref{avkN}. Collecting \eqref{Rosc=0}--\eqref{Rosc2} together, we deduce from \eqref{Roscdecom} that
\begin{align*}
&\|\mathcal{R}\div(w^{(p)}_{q+1}\otimes w^{(p)}_{q+1}+\RRR_q)\|_{\alpha}\\
\lesssim&
\Big(\frac{\delta_{q+1}\ell^{-1+\frac{3\alpha}{2}}_q}{\lambda^{1-\alpha}_{q+1}}+\frac{\delta_{q+1}\ell^{-N_0-1+\frac{\alpha}{2}}_q}{\lambda^{N_0-\alpha}_{q+1}}
+\frac{\delta_{q+1}\ell^{-1+\frac{9\alpha}{8}}_q}{\lambda^{1-\alpha}_{q+1}}+\frac{\delta_{q+1}\ell^{-N_0-1+\frac{\alpha}{8}}_q}{\lambda^{N_0-\alpha}_{q+1}}\Big)\sum_{\substack{j,j'\in I^{\theta}_k\cup I^v_k\\j+j'\neq 0}}|jj'|^{-6}\\
\lesssim&\frac{\delta^{1/2}_{q+1}\delta^{1/2}_{q}\lambda^{1+\frac{3\alpha}{2}}_q}{\lambda^{1-\alpha}_{q+1}}+\frac{\delta^{1/2}_{q+1}\delta^{1/2}_{q}\lambda^{1+\frac{3\alpha}{2}}_{q}}{\lambda^{N_0-\alpha}_{q+1}\ell^{N_0}_q}
\lesssim\frac{\delta^{1/2}_{q+1}\delta^{1/2}_{q}\lambda_q}{\lambda^{1-3\alpha}_{q+1}}.
\end{align*}
Combining with \eqref{Roscs}, we complete the proof of  Proposition \ref{propRosc}.
\end{proof}
\subsection{Proof of Proposition \ref{p:main-prop} }
In order to prove Proposition \ref{p:main-prop}, we need to show that $(v_{q+1}, \theta_{q+1}, \RR_{q+1}, \MM_{q+1})$ satisfies \eqref{p_q}, \eqref{e:vq-C0}--\eqref{e:initial} with $q$ replaced by $q+1$.

First of all, we infer from \eqref{ppp} that $\int_{\TTT^3}\bar{p}_q\dd x=0$. Since $p_{q+1}=\bar{p}_q$, we get that $p_{q+1}$ satisfies \eqref{p_q} at $q+1$ level.
Next, By \eqref{e:stability-vv_q-N}, \eqref{estimate-w} and \eqref{estimate-d}, we have
\begin{align*}
\|(v_{q+1}, \theta_{q+1})\|_0\le&\|(\vv_{q}, \tb_{q})\|_0+\|(w_{q+1}, d_{q+1})\|_0\\
\le&\sum_{i=1}^q\delta^{1 / 2}_i+\delta_{q+1}^{1 / 2} \ell_q^\alpha+\frac{3}{8}\delta^{1/2}_{q+1}\le \sum_{i=1}^{q+1}\delta^{1 / 2}_i,
\end{align*}
This inequality shows that $(v_{q+1}, \theta_{q+1})$ satisfies \eqref{e:vq-C0} at $q+1$ level. By virtue of \eqref{e:vv_q-bound}, \eqref{estimate-w} and \eqref{estimate-d}, we obtain for large enough $a$,
\begin{align*}
\|(v_{q+1}, \theta_{q+1})\|_1\le &\|(\vv_{q}, \tb_{q})\|_1+\|(w_{q+1}, d_{q+1})\|_1\\
\le&C \delta_{q}^{1 / 2} \lambda_{q} \ell_q^{-N} +\tfrac{3}{8}\delta_{q+1}^{1 / 2} \lambda_{q+1}\le \delta_{q+1}^{1 / 2} \lambda_{q+1}.
\end{align*}
Collecting Proposition \ref{proptrans}--Proposition \ref{propRosc} together shows that
\begin{align*}
\|(\RR_{q+1}, \MM_{q+1})\|_{\alpha}\le  \frac{C\delta_{q+1}^{1/2}\delta_q^{1/2} \lambda_q}{\lambda_{q+1}^{1-4\alpha}}.
\end{align*}
Therefore, we have by \eqref{alpha0} that
\begin{align*}
\|(\RR_{q+1}, \MM_{q+1})\|_{\alpha}\le C  \delta_{q+2}\lambda_{q+1}^{-4\alpha}\le\delta_{q+2}\lambda_{q+1}^{-3\alpha}.
\end{align*}
Notice  that
\[\supp w_{q+1},\,\, \supp d_{q+1}\subset \cup_{i=0}^{i_{\max}-1}(I_i+\big[-\tfrac{\tau_q}{6}, \tfrac{\tau_q}{6}\big])\subset\big[2\widetilde{T}+\tfrac{\tau_1}{6}, 3\widetilde{T}-\tfrac{\tau_q}{6}\big],\]
we deduce
$$(v_{q+1}, \theta_{q+1})(t,x)=(\vv_q, \tb_q)(t,x), \quad t\in[0, 2\widetilde{T}+\tfrac{\tau_1}{6}]\cup [3\widetilde{T}-\tfrac{\tau_q}{6}, 1], \,x\in\TTT^3.$$ Combining with \eqref{vvqtbq}, we gets
\begin{align*}
(\vv_{q}, \tb_q)\equiv (v^{(1)}, b^{(1)}) \,\,{\rm on}\,\,\big[0, 2\widetilde{T}+\tfrac{\tau_q}{6}\big],\quad (\vv_{q}, \tb_q)\equiv (v^{(2)}, b^{(2)}) \,\,\text{\rm on}\,\,\big[3\widetilde{T}-\tfrac{\tau_q}{6}, 1\big],
\end{align*}
which implies that \eqref{e:initial1} and \eqref{e:initial} hold at $q+1$ level due to $20\tau_{q+1}<\tau_q$ for $b>1$ and large~$a$. Hence, we complete the proof of Proposition \ref{p:main-prop}.

As mentioned in Section 4.1, Proposition \ref{p:main-prop} implies Theorem~\ref{t:main}. Thus,  we complete the proof of Theorem~\ref{t:main}.

\section*{Acknowledgement}
{We  thank associated editor and  the anonymous referee  for their
invaluable comments  which helped to improve the paper.  This work was supported by the National Key Research and Development Program of China (No. 2022YFA1005700) and NSFC Grant No.12371095 and No.12071043. }

\end{document}